\documentstyle[11pt]{article}

\input amssym

\font\tinybbfont=msbm6
\font\scriptsizebbfont=msbm7 scaled \magstep 1
 1
\font\footnotesizebbfont=msbm9 scaled \magstep 0
\font\smallbbfont=msbm7 scaled \magstep 2
\font\bbfont=msbm9 scaled \magstep1  
\font\largebbfont=msbm10 scaled \magstep 1
\font\Largebbfont=msbm10 scaled \magstep 2
 3
 4

\def\tinyBbb#1{\hbox{\tinybbfont #1}}
\def\scriptsizeBbb#1{\hbox{\scriptsizebbfont #1}}

\def\footnotesizeBbb#1{\hbox{\footnotesizebbfont #1}}
\def\smallBbb#1{\hbox{\smallbbfont #1}}
\def\Bbb#1{\hbox{\bbfont #1}}
\def\largeBbb#1{\hbox{\largebbfont #1}}
\def\LargeBbb#1{\hbox{\Largebbfont #1}}

\newcommand{\Ann}{\mbox{\it Ann}\,}

\newcommand{\AzumayaSpace}{\mbox{\it ${\cal A}$zumaya${\cal S}$pace}\,}

\newcommand{\Brane}{\mbox{\it ${\cal B}$rane}\,}
 \newcommand{\scriptsizeBrane}{\mbox{\scriptsize\it ${\cal B}$rane}\,}

\newcommand{\Centralizer}{\mbox{\it Centralizer}\,}
\newcommand{\Centralizersheaf}{\mbox{\it ${\cal C}$entralizer}\,}
\newcommand{\Chow}{\mbox{\it Chow}\,}
 \newcommand{\scriptsizeChow}{\mbox{\scriptsize\it Chow}}

\newcommand{\Diag}{\mbox{\it Diag}\,}
 \newcommand{\scriptsizeDiag}{\mbox{\scriptsize\it Diag}\,}

\newcommand{\End}{\mbox{\it End}\,}

\newcommand{\GL}{\mbox{\it GL}}
 
\newcommand{\Gr}{\mbox{\it Gr}\,}
\newcommand{\Hilb}{\mbox{\it Hilb}\,}
 \newcommand{\scriptsizeHilb}{\mbox{\scriptsize\it Hilb}}

\newcommand{\Id}{\mbox{\it Id}}
 \newcommand{\scriptsizeId}{\mbox{\scriptsize\it Id}}

\newcommand{\Ker}{\mbox{\it Ker}\,}

\newcommand{\Lie}{\mbox{\scriptsize\it Lie}}
 \newcommand{\tinyLie}{\mbox{\tiny\it Lie}}
\newcommand{\Map}{\mbox{\it Map}\,}

\newcommand{\Mod}{\mbox{\it Mod}\,}
\newcommand{\Mor}{\mbox{\it Mor}\,}
 
 \newcommand{\footnotesizeMor}{\mbox{\footnotesize\it Mor}\,}

\newcommand{\Proj}{\mbox{\it Proj}\,}

\newcommand{\PP}{\mbox{\it P$\!$P}}

 \newcommand{\footnotesizeQFT}{\mbox{\footnotesize\it QFT}\,}

\newcommand{\Rep}{\mbox{\it Rep}\,}

\newcommand{\RingSystemCategory}{\mbox{\it ${\cal R}$ing${\cal S}ystem$}}
\newcommand{\Scheme}{\mbox{\it ${\cal S}$cheme}\,}

\newcommand{\Space}{\mbox{\it Space}\,}

 \newcommand{\footnotesizeSpace}{\mbox{\footnotesize\it Space}\,}

\newcommand{\Spec}{\mbox{\it Spec}\,}
 \newcommand{\boldSpec}{\mbox{\it\bf Spec}\,}
 \newcommand{\scriptsizeSpec}{\mbox{\scriptsize\it Spec}\,}
 \newcommand{\footnotesizeSpec}{\mbox{\footnotesize\it Spec}\,}

\newcommand{\Stab}{\mbox{\it Stab}\,}

\newcommand{\Sym}{\mbox{\it Sym}}

\newcommand{\an}{\mbox{\scriptsize\it an}}

\newcommand{\coker}{\mbox{\it coker}\:}
\newcommand{\degree}{\mbox{\it deg}\,}

\newcommand{\dimm}{\mbox{\it dim}\,}

\newcommand{\image}{\mbox{\it Im}\,}
 \newcommand{\scriptsizeimage}{\mbox{\scriptsize\it Im}\,}

\newcommand{\scriptsizemin}{\mbox{\scriptsize\rm min}}

\newcommand{\nc}{\mbox{\scriptsize\it nc}}
 \newcommand{\tinync}{\mbox{\tiny\it nc}}
\newcommand{\scriptsizenoncommutative}{\mbox{\scriptsize\rm
                                             noncommutative}}
\newcommand{\scriptsizeorb}{\mbox{\scriptsize\rm orb\,}}

\newcommand{\pt}{\mbox{\it pt}}
\newcommand{\rank}{\mbox{\it rank}\,}

\newcommand{\redscriptsize}{\mbox{\scriptsize\rm red}\,}

\newcommand{\ringsetscriptsize}{\mbox{\scriptsize\rm ring-set}\,}
 \newcommand{\ringsettiny}{\mbox{\tiny\rm ring-set}\,}
\newcommand{\scriptsizescheme}{\mbox{\scriptsize\it scheme}}
\newcommand{\scriptsizest}{\mbox{\scriptsize\it st}}

\newcommand{\type}{\mbox{\it type}}

\newcommand{\dottedrightarrow}{\;\;
                               \mbox{-}\;\mbox{-}\;\mbox{-}\rightarrow\,}
\newcommand{\doublearrow}{\begin{array}{c}
                            \rightarrow \\[-2ex]
                            \rightarrow \end{array}}

\newcommand{\pprec}{\prec\!\prec}
\newcommand{\precleftarrow}{\vspace{-.2ex}\begin{array}{c}\prec \\[-1.9ex]
                                             \leftarrow \end{array}}
\newcommand{\succrightarrow}{\vspace{-.2ex}\begin{array}{c} \succ \\[-1.9ex]
                                              \rightarrow \end{array}}

\begin{document}

\enlargethispage{23cm}

\begin{titlepage}

$ $

\vspace{-1cm}

\noindent\hspace{-1cm}
\parbox{6cm}{\small August 2007}\
   \hspace{6.2cm}\
   \parbox[t]{5cm}{math.AG/yymm.nnnn \\ D(1): NCAG, D0.}

\vspace{2cm}

\centerline{\large\bf
 Azumaya-type noncommutative spaces and morphisms therefrom:}
\vspace{1ex}
\centerline{\large\bf
 Polchinski's D-branes in string theory from Grothendieck's viewpoint}

\bigskip

\vspace{3em}
\centerline{\large
  Chien-Hao Liu
  \hspace{1ex} and \hspace{1ex}
  Shing-Tung Yau
}

\vspace{3em}

\begin{quotation}
\centerline{\bf Abstract}

\vspace{0.3cm}

\baselineskip 12pt  
{\small
 Grothendieck's equivalence of a commutative function ring and
  a local geometric space gives rise to the language of schemes and
  functor of points in 1960s that rewrote commutative algebraic geometry
 while Polchinski's identification/recognition in 1995 of D-branes
   -- studied since the second half of 1980s as boundary conditions
      for open strings --
   as the source of Ramond-Ramond fields created by closed superstrings
   in the space-time
  rewrote string theory.
 In this work, we explain how a noncommutative version of
  Grothendieck's equivalence gives rise to a prototype intrinsic
  definition of D-branes that can reproduce the key,
  originally open-string-induced, properties of D-branes
  described in Polchinski's works.
 After the discussion of Azumaya-type noncommutative spaces and
  morphisms therefrom that form the algebro-geometric foundation
  of the current work,
 basic properties of D0-branes on a smooth curve/surface or
  a quasi-projective variety, the associated Chan-Paton modules,
  the Higgsing/un-Higgsing behavior -- all under the current setting --,
  and their relation with Hilbert schemes and Chow varieties are given.
 When applied to the case of D0-branes on a (commutative) projective
  complex smooth surface, this gives also a picture in the current
  pure algebro-geometric setting that resembles gas of D0-branes
  in a work of Vafa.
 Related supplementary discussions/remarks are given in footnotes.
} 
\end{quotation}

\bigskip

\baselineskip 12pt
{\footnotesize
\noindent
{\bf Key words:} \parbox[t]{14cm}{
  D-brane, 
  Polchinski;
  noncommutative geometry, Azumaya, Grothendieck; 
  D0-brane, moduli space.}}

\bigskip

\noindent {\small MSC number 2000:
 14A22, 81T30; 14A10, 16G30, 81T75.
} 

\bigskip

\baselineskip 10pt
{\scriptsize
\noindent{\bf Acknowledgements.}
 We thank
   Andrew Strominger and Cumrun Vafa
  for lectures/discussions that influence the project;
   Duiliu-Emanuel Diaconescu, Gang Liu, Kefeng Liu, Pan Peng
   for discussions on open/closed string duality;
  Dep't of Mathematics of UCLA for hospitality while the work
   is in preparation.
 Strings/branes is a topic
  that mixes notions/techniques/art from physics/QFT/SUSY with
  languages/techniques/naturality from mathematics.
 For that, C.-H.L.\ thanks in addition
  Paul Aspinwall, Michael Douglas, Jeffrey Harvey, Joseph Polchinski,
  Ashoke Sen
   for discussions at three TASI schools at U.\ Colorado
   at Boulder and
  Orlando Alvarez, Philip Candelas, Hungwen Chang, Chong-Sun Chu,
  Xenia de la Ossa, Jacques Distler, Daniel Freed, Joe Harris,
  Pei-Ming Ho, Kentaro Hori, Shinobu Hosono, Albrecht Klemm,
  Shiraz Minwalla, Gregory Moore, Mircea Mustata, Rafael Nepomechie,
  Mihnea Popa, Lisa Randall, Margaret Symington, Zheng Yin,
  Barton Zwiebach
   for joint influence/education over the years in the background;
  John Beachy
   for a discussion on noncommutative localizations;
  Melanie Becker, Alex Maloney
   for talks;
  D.-E.D., G.L.\
   for enlightenment on subtleties of open-string world-sheet instantons;
  Lubos Motl, Li-Sheng Tseng, Ilia Zharkov
   for conversations;
  William Oxbury
   for preprint/communication [Ox];
  participants of GSS, organized by Monica Maria Guica,
   for discussions on themes in string theory;
  Nima Arkani-Hamed, John Duncan, Dennis Gaitsgory, Katrin Wehrheim
   for topic courses;
  Rev.\ Robert Campbell Willman and Betty
   for hospitality/dinner and a discussion on a book of Brian Greene;
  Ling-Miao Chou for the moral support.
 The project is supported by NSF grants DMS-9803347 and DMS-0074329.
} 

\end{titlepage}

\newpage
\begin{titlepage}

$ $

\vspace{12em} 

\centerline{\small\it
 Chien-Hao Liu dedicates this work to his teacher Ann L.\ Willman,}
\centerline{\small\it
 who is giving him yet another lesson}
\centerline{\small
 -- {\it the grace, courage, will power, and inner peace
         while in the turmoil of life} --}
\centerline{\small\it
 throughout the treatment of her cancer.}

\end{titlepage}

\newpage
$ $

\vspace{-4em}  

\centerline{\sc
 Azumaya-Type Noncommutative Spaces, Morphisms, and D-Branes}

\vspace{2em}

\baselineskip 13pt  

\begin{flushleft}
{\Large\bf 0. Introduction and outline.}
\end{flushleft}

\begin{flushleft}
{\bf Introduction.}
\end{flushleft}
A {\it D-brane} (in full name: {\it Dirichlet brane} or
 {\it Dirichlet membrane})\footnote{D-brane
             theory and open string theory are in a way counterpart
             to and interacting with each other. As a consequence,
             supersymmetric D-brane theory and open Gromov-Witten
             theory are closely related.
            In a train of communications with Duiliu-Emanuel Diaconescu
             [Dia] on a vanishing lemma in the last section of [L-Y3]
             and its comparison with [D-F], he drew our attention to
             the important distinction between pure open GW-invariants
             and open-string world-sheet instantons.
            The former depends only on the boundary condition set up
              on the stable maps by supersymmetric D-branes and
              a decoration on the brane (cf.\ [L-Y2: Sec.~7.2])
             while the latter may interact via Wilson loops with
              the {\it general} gauge fields on the D-branes as well
              (cf.\ [Wi2: Sec.~4.2 and Sec.~4.4] and
                    [D-F: introduction part of Sec.~3]).
            Thus,
             {\it D-brane theory and the field theory thereupon
             are a part in understanding open-string world-sheet
             instantons beyond the pure Gromov-Witten sector}.
            We attribute this footnote to him and thank him
             for the patient explanations of [D-F]
           to us.}
 in string theory is by definition (i.e.\ by the very word `Dirichlet')
 a boundary condition for the end-points of open strings.
From the viewpoint of the field theory on the open-string world-sheet
 aspect, it is a boundary state in the $d=2$ conformal field theory
 with boundary.
From the viewpoint of open string target space(-time) $M$, it is
 a cycle or a union of submanifolds $Z$ in $M$ with a gauge bundle
 (on $Z$) that carries the Chan-Paton index for the end-points of
 open strings.
For the second viewpoint, Polchinski recognized in 1995 in [Pol2]
 that a D-brane is indeed a source of the Ramond-Ramond fields
 on $M$ created by the oscillations of closed superstrings in $M$.
In particular, in a specific region of the Wilson's theory-space for
 D-branes, D-branes can be identified with the solitonic/black
 branes studied earlier\footnote{See
                            [D-K-L] for a review and more references.}
 in supergravity and (target) space-time aspect of superstrings.
This recognition is so fundamental that it gave rise to
 the second revolution of string theory.
When $M$ is compactified on a Calabi-Yau space $Y$, the preservation
 of supersymmetries in either the field theory on the open-string
 world-sheet or in the effective field theory after
 the compactification requires the D-brane to be supported on
 a union of Lagrangian submanifolds/subspaces or holomorphic cycles,
(cf.\ [B-B-St], [H-I-V], and [O-O-Y]).
When we focus only on the internal/compactified part of space-time,
 this gives us a preliminary mathematical definition of
 supersymmetric D-branes as a union of Lagrangian submanifolds
 with gauge bundles or a coherent (possibly torsion) sheaf on $Y$.
While such definitions of D-branes is already very convenient in
 the study of superstring theory with branes and of stringy dualities,
they are not adequate to serve as the intrinsic definition of D-branes
 as, among other issues, in general they cannot reproduce by themselves
 a key property of D-branes
 -- the Higgsing/un-Higgsing behavior of D-branes --
 in its own mathematical framework in a natural way.

This subtlety actually does not seem to bother string theorists,
 likely for two reasons:
 \begin{itemize}
  \item[(1)]
   The picture of supersymmetric D-branes as cycles in $Y$
    with a gauge bundle is generically correct/enough
    in the regime where branes are still branes.

  \item[(2)]
   Under deformations of D-branes for which the mathematical
   picture in Item (1) is not complete enough to dictate
   the details, the very definition of D-branes as where
   open strings end tells us that we can look at
   the related open string theory, particularly its induced fields
    and their effective action on the brane, to determine
   what happens to the deformed D-branes.
 \end{itemize}
Depending on one's taste/weight on such a subtlety,
one is either satisfied with this picture or not.
And if not, one is led to the following question:
 \begin{itemize}
  \item[$\cdot$]
  {\bf Q.\ [D-brane]}$\;$
  {\it What is a D-brane intrinsically?}
 \end{itemize}
In other words, what is the intrinsic definition of D-branes
 so that by itself it can produce the properties of D-branes
 that are consistent with, governed by, or originally produced by
 open strings as well?
This is the guiding question of the current work.

The answer to this question is indeed already suggested by string
 theorists: it is hinted already in the works (e.g.\ [Pol3])
 of Polchinski and later put with even more weight
 by other string theorists\footnote{See, for example,
                                     [Dou4] and [Dou5] of Douglas and
                                     [S-W2] of Seiberg and Witten
                                     for the development and
                                     more references up to 1999.}
 that D-branes have a close tie with noncommutative geometry.
One cannot expect to have a good answer to Question [D-brane]
 without bringing appropriate noncommutative geometry into
 the intrinsic definition of D-branes.
Indeed,
 Polchinski's {\it description of deformations of stacked D-branes}
  together with
 Grothendieck's {\it local equivalence of rings
  and spaces/geometries} and the notion of {\it functors of points}
  (Sec.~2.1) implies immediately (Sec.~2.2):
 \begin{itemize}
  \item[$\cdot$]
  \parbox[t]{13.6cm}{{\bf Polchinski-Grothendieck Ansatz
                          [D-brane: noncommutativity].} {\it
   The world-volume of a D-brane carries a noncommutative structure
   locally associated to a function ring of the form $M_n(R)$,
    i.e., the $n\times n$ matrix-ring over a ring $R$
          for some $n\in {\Bbb Z}_{\ge 1}$.}}
  \end{itemize}
This brings us to a technical world in mathematics:
 noncommutative geometry.
Due to the different languages
 used in differential geometry and in algebraic geometry
 for noncommutative geometry
 (though the philosophy to equate locally a space and
  a function ring in each category is in common),
 we focus now on supersymmetric D-branes of B-type,
 for which algebro-geometric language is appropriate.

From the basic properties of D-branes spelt out explicitly
 in the work of Polchinski, there are a special class of
 noncommutative spaces that are particularly related to D-branes,
 namely the Azumaya-type noncommutative spaces.
These are the noncommutative spaces that locally have their
 function ring the matrix ring $M_n(R)$ over a commutative ring $R$.
The ansatz of Grothendieck on the equivalence of a ring and a local
 geometry, when extended to the noncommutative case as well,
 enables us to directly look at rings themselves
 without having to deal with the technical subtle issue of
 the functorial construction of an associated space
 (i.e. a set of points with topology and other structures)
 to a ring as Grothendieck did in 1960s for commutative rings
 that rewrote commutative algebraic geometry.
His ansatz of the contravariant equivalence of
 morphisms-between-spaces and morphisms-between-rings-locally,
 and the ansatz of composability, which says that
 the composition of morphisms $X\rightarrow Y$, $Y\rightarrow Z$
  between spaces should be a morphism $X\rightarrow Z$,
 can then be used to give the notion of morphisms
 from an Azumaya-type noncommutative space to
 a (either commutative or noncommutative) space without having
 the spaces themselves.
In this way, an Azumaya-type noncommutative space $X$ can be phrased
  purely as a gluing system ${\cal R}$ of matrix rings and
 a morphism from $X$ can be phrased purely as a gluing system of
  ring-homomorphisms to ${\cal R}$.
A quasi-coherent sheaf on $X$ is then a gluing system of modules
 over rings in ${\cal R}$. (Sec.~1.)

Once this language is formulated precisely,
the following prototype definition of D-branes (of B-type and
 when a ``brane" is still a brane) (Definition 2.2.3):
 \begin{itemize}
  \item[$\cdot$]
   \parbox[t]{13.6cm}{{\bf Definition [D-brane].} {\rm
    A {\it D-brane} is an Azumaya-type noncommutative space $X$
     with a fundamental module (i.e.\ the {\it Chan-Paton sheaf})
     of its noncommutative structure sheaf.
    A {\it D-brane on an open-string target-space $Y$}
     is the image of a morphism from such an $X$ to $Y$
     with the push-forward Chan-Paton sheaf.}
    } 
 \end{itemize}
 alone gives a Higgsing/un-Higgsing property of D-branes
 in its own right that is consistent and originally deduced
 via open strings in the work of Polchinski;
 (Sec.~2.2 for highlights for general D-branes;
  Sec.~3.2 for the case of D0-branes; and
  Sec.~4.1 - 4.4 for D0-branes
   on a commutative quasi-projective space).
In particular, except that we have to stay on algebraic groups
 in the pure algebro-geometric setting,
D0-branes in the current setting that move on a (commutative)
 smooth complex projective surface $Y$ has the same Higgsing/un-Higgsing
 feature of gas of D0-branes in [Vafa1] of Vafa when we choose
 the morphims of the D0-brane to $Y$ appropriately;
 (the last theme in Sec.~4.4).
The anticipation (Sec.~4.5)
 that:
 \begin{itemize}
  \item[$\cdot$]
   \parbox[t]{13.6cm}{{\bf Anticipation
       [universal moduli space from D-branes].} {\it
    The moduli space of D-branes
    -- or in general of D-branes coupled with NS-branes
      when defined correctly -- on a target space
    should encompass simultaneously
    several standard moduli spaces in commutative geometry.}
   } 
 \end{itemize}
 is supported in the study of the moduli space of D0-branes;
 (Sec.~3 and Sec.~4.1- Sec.~4.4).

Finally, a word about reading the current work:
Noncommutative geometry, in the language of either differential
 geometry or algebraic geometry, is a demanding topic and
there is no way to bypass it.
Readers who already know D-branes in the string-theoretic aspect
  from [Pol3] or [Pol4] are suggested to read Sec.~4.1 first to see
 how algebraic geometry in the line of Grothendieck is used
  to implement Polchinski's picture in a most elementary case:
 D0-branes on the complex line ${\Bbb C}$.
Various general features of D-branes and their moduli space,
 following the above prototype definition, reveal themselves
 already in this example in a simplified form.


\bigskip

\noindent
{\it Remark 0.1 $[$diverse D-``branes"$]$.} {\rm
 Mathematicians should be aware that
  there are numerous string theorists whose collective contribution
  shaped the understanding of D-branes nowadays, cf.\ the limited
  ``short'' list of stringy references of the current work, which have
  influenced us and became part of the background of the project.
 Their works led to diverse meanings/roles of D-branes in various
  physical contents.
 The current work addresses D-branes when they are ``still branes",
  i.e.\ in the sense of [D-L-P], [P-C], [Pol2], [Pol3], [Pol4],
  and, e.g., [B-V-S1], [B-V-S2], [Vafa1], [Vafa2]
  that they are manifold/variety-type objects.
 The terms `Polchinski's D-brane' and `D0-brane gas'
  occasionally used in this work refer to [D-L-P], [Vafa1], and
  Polchinski's special contribution to this topic.
 Physicists use the same term `D-branes' in the various different
  physical contents with good reasons, particularly from the aspect
  of stringy dualities.
 However, this is unfortunate/inconvenient for us as these other
  types of D-``branes" are no longer branes and have/involve
  very different mathematical contents/language as well.
 Lacking an official terminology, we use above-mentioned terms
  and terms like `D-branes in the sense of Polchinski'
  to single out the particular meaning/type of D-branes studied
  in the above-quoted stringy works in the earlier years of D-branes
  for convenience.
} 

\bigskip

\noindent
{\it Remark 0.2 $[$other brane$]$.} {\rm
 It should be mentioned that, while D-branes have been a central
  object in string theory since 1995, there are other types of branes,
  (e.g.. NS-branes) in string theory as well that serve as the source
  for other types of fields created by closed strings in space-time;
 see [Pol4], [Jo], and [B-B-Sc] for a review.
 It is also worth noting that, since the work of Randall and Sundrum
  [R-S] in 1999, the use of branes has been extended outside of
  string theory and gives a new insight to the weakness of gravity
  in comparison with electro-magnetic, weak, and strong interactions
  in nature.
 That route hints at a connection of hyperbolic geometry and branes
  -- a topic in its own right.
} 

\bigskip

\noindent
{\bf Convention.}
 Standard notations, terminology, operations, facts in
  (1) (noncommutative; commutative) ring theory;
  (2) (commutative) algebraic geometry;
  (3) quantum field theory, supersymmetry; string theory
  can be found respectively in
  (1) [Jac]; [Mat];
  (2) [Ha], [E-H];
  (3) [I-Z], [P-S], [W-B]; [B-B-Sc], [G-S-W], [Jo], [Pol4], [Zw].
 \begin{itemize}
  \item[$\cdot$]
   Except the {\it zero-ring} $0$,
    all {\it rings} or {\it algebras} (over an algebraically closed field)
    $R$ in the general discussion of this work are {\it associative}
    with an {\it identity} $1$ and
    are both {\it left- and right-Noetherian}.
   The term ``{\it $R$-modules}", including ``{\it ideals}" in $R$,
    means ``{\it left} $R$-modules" (cf.\ {\it left} ideals in $R$)
    unless otherwise noted.
   $Z(R) :=$ the {\it center} of $R$.
   $M_n(R) :=$ the $n\times n$ {\it matrix ring} with entries in $R$.

  \item[$\cdot$]
   The term {\it field} has two completely different meanings:
    field in {\it quantum field theory}
     vs.\ field in the {\it theory of rings}.

  \item[$\cdot$]
   The {\it analytic space} ${\Bbb C}^n$, with the standard topology,
    of closed points in the affine space ${\Bbb A}^n$ over ${\Bbb C}$
    is constantly denoted directly by ${\Bbb A}^n$.
   Similarly for ${\Bbb P}^n$ and other varieties.
   (In this work, we use the term `varieties/schemes' mainly only
     to manifest/emphasize the fact that they arise from gluing
     of affine charts associated to rings.)
   In this way, ${\Bbb C}^n$ is kept to mean ${\Bbb C}^n$ as
    a ${\Bbb C}$-  or $M_n({{\Bbb C}})$-{\it module} as best possible.
   ${\Bbb C}^n$ as the $n$-th {\it product ring} of ${\Bbb C}$ will
    be denoted also by $\prod_n{\Bbb C}$.

  \item[$\cdot$]
   A representation (resp.\ commuting) scheme with the reduced scheme
    structure will be called {\it representation} (resp.\ {\it commuting})
    {\it variety} for simplicity.
   Irreducibility is {\it not} implied here.
   (In fact, in general they are not irreducible.)

  \item[$\cdot$]
   Omitted subscripts (resp.\ superscripts) are indicated
    by $_{\bullet}$ (resp.\ $^{\bullet}$).
 \end{itemize}

\bigskip

\begin{flushleft}
{\bf Outline.}
\end{flushleft}
{\small
\baselineskip 11pt  
\begin{itemize}
 \item[1.]
  Azumaya-type noncommutative spaces and morphisms therefrom.
  \vspace{-.6ex}
  \begin{itemize}
   \item[1.1]
    Azumaya-type noncommutative spaces and morphisms therefrom.

   \item[1.2]
    A noncommutative space as a gluing system of rings.
  \end{itemize}

 \item[2.]
  D-branes from the viewpoint of Grothendieck.
  \vspace{-.6ex}
  \begin{itemize}
   \item[2.1]
    The notion of a space(-time): functor of points vs.\ probes.

   \item[2.2]
    D-branes as Azumaya-type noncommutative spaces.
  \end{itemize}

 \item[3.]
  $\Mor(\Space M_n({\smallBbb C}),Y)$ as a coarse moduli space.
  \vspace{-.6ex}
  \begin{itemize}
   \item[3.1]
    Central localizations of Artinian rings and their modules.

   \item[3.2]
    $\Mor(\Space M_n({\smallBbb C}),Y)$ as a coarse moduli space.
  \end{itemize}

 \item[4.]
  D0-branes on a commutative quasi-projective variety.
  \vspace{-.6ex}
  \begin{itemize}
   \item[4.1]
    D0-branes on the complex affine line ${\smallBbb A}^1$.

   \item[4.2]
    D0-branes on the complex projective line ${\smallBbb P}^1$.

   \item[4.3]
    D0-branes on the complex affine plane ${\smallBbb A}^2$.

   \item[4.4]
    D0-branes on a complex quasi-projective variety.

   \item[4.5]
    A remark on D-branes and universal moduli space.
  \end{itemize}
\end{itemize}
} 

\bigskip

\newpage
\section{Azumaya-type noncommutative spaces and morphisms\\ therefrom.}

We introduce in Sec.~1.1 a class of noncommutative spaces
  that are relevant to D-branes.
Its foundation, central localizations of noncommutative rings,
 is given in Sec.~1.2.
The ring-theoretic description of a space in Sec.~1.2 allows us
 to study as well the space of morphisms between noncommutative
 spaces without having to construct the noncommutative spaces.

\bigskip

\subsection{Azumaya-type noncommutative spaces and morphisms therefrom.}

\noindent
{\bf Definition 1.1.1 [Azumaya-type noncommutative space].} {\rm
 An {\it Azumaya-type noncommutative space} is a triple
  $(X,{\cal O}_X,{\cal O}_X^{\nc})$,
  where
   $(X,{\cal O}_X)$ is a (commutative Noetherian) scheme,
    as defined in [Ha], and
   ${\cal O}_X^{\nc}$ is a coherent sheaf of noncommutative
   ${\cal O}_X$-algebras\footnote{The category of noncommutative
                                   algebras includes also commutative
                                   algebras.
                                  We will call a sheaf ${\cal G}$ of
                                   ${\cal O}_X$-algebras simply
                                   an ${\cal O}_X$-algebra.
                                  The {\it center} ${\cal Z}({\cal G})$
                                   of ${\cal G}$ is, by definition,
                                   the sheaf associated to the presheaf
                                   that assigns to each open set $U$ of
                                   $X$ the sub-${\cal O}_X(U)$-algebra
                                   $Z({\cal G}(U))$ of ${\cal G}(U)$.}
    on $X$ that contains ${\cal O}_X$
    by $1\cdot {\cal O}_X$ in its center ${\cal Z}({\cal O}_X^{\nc})$.
 We will call ${\cal O}_X$ (resp.\ ${\cal O}_X^{\nc}$)
  the {\it commutative} (resp.\ {\it noncommutative})
  {\it structure sheaf} of $X$.
 A {\it strict morphism} from $(X,{\cal O}_X,{\cal O}_X^{\nc})$ to
  $(Y,{\cal O}_Y,{\cal O}_Y^{\nc})$ is a triple
  $(f,f^{\sharp}, f^{\sharp\nc})$,
  where $(f:X\rightarrow Y\,,\,
          f^{\sharp}:{\cal O}_Y\rightarrow f_{\ast}{\cal O}_X)$
   gives a morphism of schemes from $(X,{\cal O}_X)$ to $(Y,{\cal O}_Y)$
    and
   $f^{\sharp\nc}:{\cal O}_Y^{\nc}\rightarrow f_{\ast}{\cal O}_X^{\nc}$
    is a homomorphism of ${\cal O}_Y$-algebras that extends $f^{\sharp}$.
 A {\it general morphism} from $(X,{\cal O}_X,{\cal O}_X^{\nc})$ to
   $(Y,{\cal O}_Y,{\cal O}_Y^{\nc})$ consists of the following data:
   \begin{itemize}
    \item[$\cdot$]
     an inclusion pair
      ${\cal O}_X\subset {\cal A}\subset {\cal A}^{\nc}
                                 \subset {\cal O}_X^{\nc}$
      of ${\cal O}_X$-subalgebras such that
      ${\cal A}\subset {\cal Z}({\cal A}^{\nc})\,$;

    \item[$\cdot$]
     a strict morphism $(f,f^{\sharp},f^{\sharp\nc})$ from
      $(X^{\prime},{\cal O}_{X^{\prime}},{\cal O}_{X^{\prime}}^{\nc})$
      to $(Y,{\cal O}_Y,{\cal O}_Y^{\nc})\,$,
     where
      \begin{itemize}
       \item[-]
        $X^{\prime}:=\boldSpec{\cal A}$
         is equipped with the tautological dominant finite morphism
         $X^{\prime}\stackrel{\jmath}{\rightarrow} X$ of schemes,

       \item[-]
        ${\cal O}_{X^{\prime}}^{\nc}$
         is the ${\cal O}_{X^{\prime}}$-algebra on $X^{\prime}$
         associated to ${\cal A}^{\nc}$ as an ${\cal A}$-algebra.
      \end{itemize}
   \end{itemize}
 A strict morphism is automatically a general morphism.
 A general morphism will also be called simply a {\it morphism}.
 Define $\Mor(X,Y)$ to be the set of morphisms from $X$ to $Y$.
 To simplify the notation, we will also denote
  $(X,{\cal O}_X,{\cal O}_X^{\nc})$ collectively by $X$ and
  both a strict morphism $(f,f^{\sharp},f^{\sharp\nc})$
   and a general morphism
   $(({\cal A},{\cal A}^{\nc}),(f,f^{\sharp},f^{\sharp\nc}))$
   collectively by $f:X\rightarrow Y$.
} 

\bigskip

\noindent
{\bf Definition/Example 1.1.2 [tautological morphism/surrogate].} {\rm
 With notations from Definition 1.1.1,
  the (strict) identity morphism
   $(X^{\prime},{\cal O}_X,{\cal O}_X^{\nc})
     \rightarrow (X^{\prime},{\cal O}_X,{\cal O}_X^{\nc})$
  defines a (general) morphism
   $X=(X,{\cal O}_X,{\cal O}_X^{\nc}) \rightarrow
                 X^{\prime}=(X^{\prime},{\cal O}_X,{\cal O}_X^{\nc})$.
 Given $X$, we will call an $X\rightarrow X^{\prime}$ arising this way
  a {\it tautological morphism} from $X$  and
  $X^{\prime}$ an {\it surrogate} of $X$.
} 

\bigskip

\noindent
{\bf Example 1.1.3 [noncommutative point].}
 Let $k$ be an algebraically closed field and $M_n(k)$
  be the $k$-algebra of $n\times n$-matrices with entries in $k$.
 Then, $X=(\Spec k, k, M_n(k)) =: \Space M_n(k)$
  defines an {\it Azumaya-type noncommutative point}.
 See Sec.~3.1 for more details.

\bigskip

\noindent
{\bf Example 1.1.4 [morphism of commutative schemes].}
 An Azumaya-type noncommutative space
  $X=(X,{\cal O}_X,{\cal O}_X^{\nc})$ is a commutative scheme
  if and only if ${\cal O}_X={\cal O}^{\nc}$.
 In this case, $X$ has no surrogates except $X$ itself and
  any morphism from $X$ to $Y=(Y, {\cal O}_Y, {\cal O}_Y^{\nc})$
  is a strict morphism from $X$ to $Y$.
 In particular, the natural inclusion
  $\;\Scheme\; \hookrightarrow\; \AzumayaSpace\;$
  of the category of commutative schemes into the category of
  Azumaya-type noncommutative spaces is fully faithful.

\bigskip

The foundation of Definition 1.1.1
 (i.e.\ of the sheaf ${\cal O}_X^{\nc}$)
 is on central localizations of (noncommutative) rings.
This will be discussed in Sec.~1.2.
The following lemma follows immediately from the definition:

\bigskip

\noindent
{\bf Lemma 1.1.5 [exhaustion].} {\it
  Let $X$ and $Y$ be Azumaya-type noncommutative spaces and
   $X^{\prime}$ be a surrogate of $X$.
  Then there is a canonical embedding
   $\Mor(X^{\prime},Y) \hookrightarrow \Mor(X,Y)$.
} 

\bigskip

\noindent
{\it Remark 1.1.6 $[$noncommutative geometry$]$.}
 Noncommutative algebraic geometry was developed with vigor
  by several schools of mathematicians immediately after Grothendieck's
  re-writing of commutative algebraic geometry in the 1960s.
 There are several classes of noncommutative spaces in existence;
  each is described in its own appropriate language.
 While many demanding fundamental issues have prevented it from reaching
  at the moment the same glory and a unified language as its commutative
  counterpart from Grothendieck's school, it is a constant growing subject.
 Readers are referred to, e.g.\ (in rough historical order)
  [Go], [vO-V], [A-Z], [J-V-V], [Ro1], [Ro2], [K-R1], [K-R2]
    from the algebraic aspect;
  [Co] from the analytic aspect; and
  [Man2], [Man3], [Kapr], [Lau], [leB1] from other aspects
  for details and more references.

\bigskip

\noindent
{\it Remark 1.1.7 $[$Azumaya-type noncommutative space$]$.}
 The class of noncommutative spaces we define here, namely
  $(X,{\cal O}_X,{\cal O}_X^{\nc})$, are chosen with
  D-branes in mind.
 While they may be thought of as noncommutative ``clouds"
   (i.e.\ ${\cal O}_X^{\nc}$) over (commutative) schemes
   (i.e.\ $(X,{\cal O}_X)$),
  the way we define a morphism from $X$ to $Y$ says that
  the main object of focus in the triple
  $(X,{\cal O}_X,{\cal O}_X^{\nc})$ is ${\cal O}_X^{\nc}$,
  rather than $(X,{\cal O}_X)$.
 This particular point is important in the realization of
  a D-brane of B-type as an Azumaya-type noncommutative space.
 We suggest readers to think of
  $$
  \begin{array}{ccl}
   (X,{\cal O}_X,{\cal O}_X^{\nc})
    & \hspace{.8em}\mbox{as}\hspace{.8em}
    & \begin{array}{l}
       {\cal O}_X^{\nc}\,,\hspace{1em}
       \mbox{{\it together with}
             the system $L_{{\cal O}_X^{\tinync}}$
             of sub-${\cal O}_X$-algebra pairs:}          \\[1.6ex]
       L_{{\cal O}_X^{\tinync}}\;
         =\;\left\{\,
              ({\cal A},{\cal A}^{\nc})\;
              \left|\;
              \begin{array}{l}
               {\cal O}_X \subset {\cal A} \subset {\cal A}^{\nc}
                          \subset {\cal O}_X^{\nc}\,;     \\[.6ex]
               \mbox{${\cal A}$, ${\cal A}^{\nc}\,$:
                                 sub-${\cal O}_X$-algebras}\,;\,
               {\cal A} \subset {\cal Z}({\cal A}^{\nc})\,
              \end{array}
              \right.
            \right\}.
      \end{array}
   \end{array}
  $$
 I.e.\ $(X,{\cal O}_X,{\cal O}_X^{\nc})$
  together with the system $\{X \rightarrow X^{\prime}\}_{X^{\prime}}$
  of surrogates in $\AzumayaSpace$.

\bigskip

\noindent
{\bf Example 1.1.8 [noncommutative point revisitd].}
 (Continuing Example 1.1.3.)
 A surrogate of the Azumaya-type noncommutative point
  $\Space M_n({k})$ over $k$
  is given by a sub-$k$-algebra pair $k\subset C\subset R\subset M_n(k)$
  with $C\subset Z(R)$.
 In particular,
  while $\Space M_n(k)$ consists geometrically of only one point
   (i.e.\ $\Spec k$),
  its surrogate $X^{\prime}=(\Spec C,C,R)$ can have more than one
   geometric points in $\Spec C$.
 All these $X^{\prime}$'s should be thought of as part of
  the ``geometry" of noncommutative point $\Space M_n(k)$.

\bigskip

\noindent
{\bf Definition 1.1.9 [left/right quasi-coherent/coherent sheaf].} {\rm
 A {\it left quasi-coherent sheaf}
  on $(X,{\cal O}_X, {\cal O}_X^{\nc})$
  is a sheaf of left ${\cal O}_X^{\nc}$-modules that is quasi-coherent
  on $(X,{\cal O}_X)$.
 Similarly for the definitin of
  a {\it right quasi-coherent sheaf},
  a {\it left coherent sheaf}, and
  a {\it right coherent sheaf} on $(X,{\cal O}_X,{\cal O}_X^{\nc})$.
 A {\it ${\cal O}_X^{\nc}$-module} is by convention
  a left ${\cal O}_X^{\nc}$-module.
} 

\bigskip

Hidden in the notion of ${\cal O}_X^{\nc}$
  in the tuple $(X,{\cal O}_X, {\cal O}_X^{\nc})$
 is the notion of central localizations,
 which we will discuss more thoroughly in Sec.~1.2.

\bigskip

\subsection{A noncommutative space as a gluing system of rings.}

The purely ring\footnote{For non-algebraic-geometers:
                                 A {\it ring} $R$ here is meant to be
                                  the {\it ring of functions} on
                                   a ``space" $X_R$ these functions
                                   are supposed to take
                                   as their defining domain, and
                                 a {\it ring-homomorphism}
                                  $R\rightarrow S$ is meant to be
                                  the {\it pulling-back} of functions
                                   on the underlying spaces when
                                   there is a map/morphism
                                   $X_S\rightarrow X_R$
                                   between the spaces.
                                 Algebraic geometers have turn the picture
                                  of ``space first, function-ring second"
                                  around to make the function-ring first
                                  and space --
                                  if functorially constructible at all
                                  -- second.
                                 Indeed, physicists have already adopted
                                  such ``function-ring first" philosophy
                                  (without knowing the ``space")
                                  when studying supersymmetry and
                                  superfields on a superspace.}-theoretic
 construction in this subsection enables us to talk about
 a ``noncommutative scheme" without having to construct one\footnote{The
                                      general functorial construction of
                                       noncommutative schemes that
                                       generalizes Grothendieck's school
                                       on commutative geometry
                                       is a subtle issue.
                                      See, e.g., [J-V-V: introduction]
                                       and Remark 1.1.6.}.
The ring system to be defined is meant to carry the same information
 as the noncommutative scheme associated to ${\cal O}_X^{\nc}$
 in $X=(X,{\cal O}_X,{\cal O}_X^{\nc})$ would.
Such a description will later be used to study $\Mor(X,Y)$.
Behind the messy notations is the notion of Grothendieck-descent-data
 description of spaces/stacks/sheaves and morphisms between them.

\bigskip

\begin{flushleft}
{\bf Noncommutative localizations.}
\end{flushleft}
%
The notion of noncommutative localizations can be traced back to
 Ore in [Or1] and [Or2] in 1930s.
Here we recall only definitions that will be needed later. See
e.g.\ [Ga], [Goldm], [Jat], [St]
 for more details and thorough discussions.

A {\it Gabriel filter} on a ring $R$ is a collection ${\frak F}$
 of ideals in $R$ that satisfies\footnote{Property (1)
                                     and Property (2) together define
                                      the notion of a {\it filter} of
                                      ideals in $R$;
                                     Property (3) and Property (4)
                                      together actually imply
                                      Property (1) and Property (2).}:
 \begin{itemize}
  \item[(1)]
   if $I\in{\frak F}$ and $J$ is an ideal that contains $I$,
   then $J\in {\frak F}$;

  \item[(2)]
   if $I$, $J\in {\frak F}$, then $I\cap J\in {\frak F}$;

  \item[(3)]
   if $I\in {\frak F}$, then $(I:r)\in {\frak F}$ for $r\in R$;

  \item[(4)]
   if $I\in {\frak F}$ and $J$ is an ideal such that
    $(J:r)\in{\frak F}$ for all $r\in I$,
   then $J\in{\frak F}$.
 \end{itemize}
Each Gabriel filter ${\frak F}$ on $R$ determines
 the subcategory ${\cal T}_{\frak F}$
  of ${\frak F}$-{\it torsion} objects  and
 the subcategory ${\cal F}_{\frak F}$
  of ${\frak F}$-{\it torsion-free} objects
 in the category $R$-$\Mod$ of (left) $R$-modules.
An object $M$ in ${\cal T}_{\frak F}$ is characterized by that
 each element $m$ of $M$ has its annihilator $\Ann(m)\in {\frak F}$;
  and
an object $N$ in ${\cal F}_{\frak F}$ is characterized by that
 $N$ contains no submodule in ${\cal T}_{\frak F}$ except
 the zero submodule $0$.
Each object in $M\in R$-$\Mod$ fits into an exact sequence
 $0 \rightarrow M^{\prime} \rightarrow M
    \rightarrow M^{\prime\prime} \rightarrow 0$,
 where
   $t_{\frak F}(M):= M^{\prime}\in {\cal T}_{\frak F}$ and
   $M^{\prime\prime}\in {\cal F}_{\frak F}$.
In particular, $M\in {\cal T}_{\frak F}$ (resp.\ ${\cal F}_{\frak F}$)
 if and only if $t_{\frak F}(M)=M$ (resp.\ $t_{\frak F}(M)=0$).
The {\it localization} $M_{\frak F}$ of $M\in R$-$\Mod$
 with respect to ${\frak F}$ is defined to be the
 ${\frak F}$-{\it injective envelop} $E_{\frak F}(M/t(M))$ of the
 ${\frak F}$-torsion-free quotient module $M/t_{\frak F}(M)$ of $M$.
When ${\frak F}$ is clear or omitted from the text,
 $E_{\frak F}$, ${\cal F}_{\frak F}$, $t_{\frak F}$,
 ${\cal T}_{\frak F}$, ``${\frak F}$-torsion", and
 ``${\frak F}$-torsion-free"
 will be denoted/called simply $E$, ${\cal F}$, $t$, ${\cal T}$,
 ``{\it torsion}", and ``{\it torsion-free}" respectively.

%
%
%
%
%
%

The following kind of localizations is closest to the localizations
 in the case of commutative rings.
It is the one used in Definition 1.1.1 for ${\cal O}_X^{\nc}$:

\bigskip

\noindent
{\bf Definition 1.2.1 [central localization]\footnote{Central
                                            localizations are particularly
                                             akin to Azumaya-type
                                             noncommutative spaces.
                                            It should be noted
                                             that most of the definitions,
                                             statements, and constructions
                                             we give based on
                                             central localizations
                                             {\it cannot} be taken
                                             directly for
                                             general localizations
                                             without additional works
                                             or modifications.}.}
{\rm
 Given a ring $R$, a {\it central localization} of $R$ is
  the localization $R_{{\frak F}_S}$ of $R$ with respect to
  the Gabriel filter ${\frak F}_S$ associated to
  a multiplicatively closed subset $S$ in the center $Z(R)$ of $R$.
} 

\bigskip

\noindent
Explicitly,
 the Gabriel filter in the above definition is given by
  ${\frak F}_S=\{I\,:\,\mbox{ideal of $R$}\;,\; I\cap S\ne \emptyset\}$
  and
 the central localization is given by
  $R_{{\frak F}_S} =
  [S^{-1}]R = R[S^{-1}] := (R\times S)/\sim$,
  where $(r_1,s_1)\sim (r_2,s_2)$ if and only if $s(r_1s_2-r_2s_1)=0$
   for some $s\in S$.

\bigskip

\noindent
{\bf Definition 1.2.2 [push-out, admissibility, and descent].} {\rm
 (1)
 Let
  $\varphi: R\rightarrow R^{\prime}$ be a ring-homomorphism,
  $S\subset Z(R)$ be a multiplicatively closed subset in $R$
  such that
   $\varphi(S)\subset Z(R^{\prime})$,  and
   $\psi:R\rightarrow R_{{\frak F}_{S}}$ be the central localization
    of $R$ with respect to $S$.
 Then
  the central localization
  $\psi^{\prime}:
    R^{\prime} \rightarrow R^{\prime}_{{\frak F}_{\varphi(S)}}$
  of $R^{\prime}$ is called the {\it push-out} of $\psi$ to $R^{\prime}$
  via $\varphi$.
 (2)
 Given central localizations
  $\psi:R\rightarrow R_{{\frak F}_{S}}$  and
   $\psi^{\prime}:R^{\prime}
            \rightarrow R^{\prime}_{{\frak F}_{S^{\prime}}}$,
 a ring-homomorphism $\varphi:R\rightarrow R^{\prime}$
   is called {\it admissible} to $(S,S^{\prime})$
  if $\varphi(S)\subset S^{\prime}$.
 For such $\varphi$, there is a canonical/unique ring-homomorphism
  $\varphi_{(S,S^{\prime})}:R_{{\frak F}_S}
     \rightarrow R^{\prime}_{{\frak F}_{S^{\prime}}}$
  that makes the following diagram commute:
  $$
   \begin{array}{cccl}
    R  & \stackrel{\varphi}{\longrightarrow}
       & R^{\prime}\hspace{.8em}                         \\
    \mbox{\scriptsize $\psi$}\downarrow\hspace{.8em}
       & & \downarrow\mbox{\scriptsize $\psi^{\prime}$}  \\[-.6ex]
    \hspace{.8em}R_{{\frak F}_S}
       & \stackrel{\varphi_{(S,S^{\prime})}}{\longrightarrow}
       & R^{\prime}_{{\frak F}_{S^{\prime}}}             & .
   \end{array}
  $$
 $\varphi_{(S,S^{\prime})}$ is called the {\it descent} of $\varphi$
  under the central localizations.
} 

\bigskip

\noindent
{\bf Example 1.2.3 [2-step consecutive central localization].}
 Given central localizations
  $\psi_1: R\rightarrow R_1$ and $\psi_2:R\rightarrow R_2$ of $R$,
 one has the push-outs $\psi_{12}: R_1\rightarrow R_{12}$ and
   $\psi_{21}: R_2\rightarrow R_{21}$
   of $\psi_2$ via $\psi_1$ and of $\psi_1$ via $\psi_2$ respectively.
 Then there is a canonical isomorphism $R_{12}\simeq R_{21}$ such that
  the following diagram
   $$
    \begin{array}{lcr}
      R    & \stackrel{\psi_2}{\longrightarrow}  & R_2\hspace{1ex}  \\
      \hspace{-.9em}\mbox{\scriptsize $\psi_1$}\downarrow
      & & \downarrow\mbox{\scriptsize $\psi_{21}$}\hspace{-.7em} \\[-.6ex]
      R_1  & \stackrel{\psi_{12}}{\longrightarrow}
           & R_{12}\; \simeq\; R_{21}
    \end{array}
   $$
   commutes.
 Both the compositions
  $\psi_{12}\circ\psi_1:R\rightarrow R_{12}$ and
  $\psi_{21}\circ\psi_2:R\rightarrow R_{21}$
  give central localizations of $R$.
 Such 2-step consecutive central localizations will appear in
  stating the cocycle conditions for the gluing of rings
  along their central localizations.

\bigskip

\begin{flushleft}
{\bf A ring-theoretic description of noncommutative spaces and
     their morphisms.}
\end{flushleft}
We give a description of a class of noncommutative spaces and
  their morphisms solely in terms of rings, ring-homomorphisms, and
  central localizations,
 {\it without} employing the notion of ``points" and ``topology"
  of a ``space".
This class contains the class of Azumaya-type noncommutative spaces
 introduced in Sec.~1.1 as a subclass.

\bigskip

\noindent
{\bf Definition 1.2.4 [finite central cover of a ring].} {\rm
 Let
  $A$ be a finite set and
  ${\cal U}:= \{\varphi_{\alpha}:R\rightarrow R_{\alpha}\}_{\alpha\in A}$
   be a finite collection of central localizations of $R$ with respect
   to Gabriel filters ${\frak F}_{\alpha}$, $\alpha\in A$, on $R$.
 We say that ${\cal U}$ is a {\it finite central cover} of $R$
  if $\sum_{\alpha\in A} I_{\alpha}=R$ for any tuple
     $(I_{\alpha})_{\alpha}\in \prod_{\alpha\in A}{\frak F}_{\alpha}$.
} 

\bigskip

\noindent
{\bf Definition 1.2.5 [gluing system of rings]\footnote{This
                                   is a Grothendieck's
                                   descent-data-of-objects description.}.}
{\rm
 A (finite) {\it gluing system} of rings
  $$
   {\cal R}\; =\; \left(
     \{R_{\alpha}\}_{\alpha\in A}
      \doublearrow \{R_{\alpha_1\alpha_2}\}_{\alpha_1,\alpha_2\in A}
                  \right)
  $$
  from central localizations consists of the following data:
  \begin{itemize}
   \item[(1)]
    [{\it local ring-charts}]\\[.6ex] 
    a finite collection $\{R_{\alpha}\}_{\alpha\in A}$ of rings;
    ($A$: the {\it index set} of ${\cal R}$)

   \item[(2)]
    [{\it transition ring-homomorphisms}]\\[.6ex] 
    a finite central cover
     $\{R_{\alpha_1}\rightarrow R_{\alpha_1\alpha_2}\}_{\alpha_2\in A}$
     for each $R_{\alpha_1}$ and
    a choice of ring-isomorphisms
     $\varphi_{\alpha_1\alpha_2}: R_{\alpha_1\alpha_2}
       \stackrel{\sim}{\rightarrow} R_{\alpha_2\alpha_1}$
     for each $(\alpha_1,\alpha_2)\in A\times A$
     such that
      $R_{\alpha\alpha}=R_{\alpha}$,
      $\varphi_{\alpha_1\alpha_2}=\varphi_{\alpha_2\alpha_1}^{\;\;-1}$,
       and
      $\varphi_{\alpha\alpha}=\Id_{R_{\alpha}}$;

   \item[$\cdot$] [{\it cocycle conditions}]\\[.6ex] 
    the ring-homomorphism $R_{\alpha_1}\rightarrow R_{\alpha_1\alpha_2}$
     pushes out the finite central cover
     $\{R_{\alpha_1}\rightarrow R_{\alpha_1\alpha_3}\}_{\alpha_3}$
      of $R_{\alpha_1}$
     to a finite central cover
      $\{R_{\alpha_1\alpha_2}\rightarrow
                           R_{\alpha_1\alpha_2\alpha_3}\}_{\alpha_3}$
      of $R_{\alpha_1\alpha_2}$
    and one has
     the canonical isomorphisms
      $R_{\alpha_1\alpha_2\alpha_3}\simeq R_{\alpha_1\alpha_3\alpha_2}$
      from the push-out diagrams;
    it is then required that
     the gluing ring-isomorphisms
      $R_{\alpha_1\alpha_2}\rightleftharpoons R_{\alpha_2\alpha_1}$
      descend to ring-isomorphisms
      $R_{\alpha_1\alpha_2\alpha_2}
        \rightleftharpoons R_{\alpha_2\alpha_1\alpha_3}$
      that make the following diagrams
      $$
       \begin{array}{ccc}
        R_{\alpha_1\alpha_2}
          & \rightleftharpoons    & R_{\alpha_2\alpha_1}  \\[.6ex]
        \downarrow    & & \downarrow          \\
        R_{\alpha_1\alpha_2\alpha_3}
          & \rightleftharpoons    & R_{\alpha_2\alpha_1\alpha_3}
       \end{array}
       \hspace{1em}\mbox{and}\hspace{1em}
       \begin{array}{c}
        R_{\alpha_1\alpha_2\alpha_3}
          \simeq R_{\alpha_1\alpha_3\alpha_2}           \\
        \swarrow\hspace{-1.2ex}\nearrow
         \hspace{10em} \nwarrow\hspace{-1.2ex}\searrow  \\[.6ex]
        R_{\alpha_2\alpha_1\alpha_3} \simeq R_{\alpha_2\alpha_3\alpha_1}
          \hspace{.8em}\rightleftharpoons\hspace{.8em}
         R_{\alpha_3\alpha_2\alpha_1}\simeq R_{\alpha_3\alpha_1\alpha_2}
       \end{array}
      $$
      commute.
    Note that under the requirement of the first diagram above
     the isomorphisms
     $R_{\alpha_1\alpha_2\alpha_3}
       \rightleftharpoons    R_{\alpha_2\alpha_1\alpha_3}$,
     when exists, are unique.
  \end{itemize}
 We will write $R_{\alpha}\in {\cal R}$ to indicate that
  $R_{\alpha}$ is a ring-chart in the system ${\cal R}$.
 A ({\it finite central}) {\it refinement} of ${\cal R}$ is a gluing system
  ${\cal R}^{\prime}
   =( \{R^{\prime}_{\alpha^{\prime}}\}_{\alpha^{\prime}\in A^{\prime}}
      \doublearrow
      \{R^{\prime}_{\alpha_1^{\prime},\alpha_2^{\prime}}\}
        _{\alpha_1^{\prime},\alpha_2^{\prime}\in A^{\prime}} )$
  of rings together with the following data:
  \begin{itemize}
   \item[$\cdot$]
    a surjective map $\tau:A^{\prime}\rightarrow A\,$;

   \item[$\cdot$]
    a central localization ring-homomorphism
     $R_{\alpha}\rightarrow R^{\prime}_{\alpha^{\prime}}$
     for each $\alpha\in A$ and $\alpha^{\prime}\in \tau^{-1}(\alpha)$
     such that
     \begin{itemize}
      \item[-]
       for each $\alpha\in A$,
       $\{R_{\alpha}\rightarrow R^{\prime}_{\alpha^{\prime}}\}
                  _{\alpha^{\prime}\in \tau^{-1}(\alpha)}$
        is a finite central cover of $R_{\alpha}$;

      \item[-]
       for all $(\alpha_1,\alpha_2)\in A\times A$ and
        $(\alpha^{\prime}_1,\alpha^{\prime}_2)
          \in \tau^{-1}(\alpha_1)\times\tau^{-1}(\alpha_2)$,
       $R_{\alpha_1}\rightarrow R^{\prime}_{\alpha^{\prime}_1}$
        descends to
        $R_{\alpha_1\alpha_2}\rightarrow
               R^{\prime}_{\alpha^{\prime}_1\alpha^{\prime}_2}$
        and
       all the diagrams
        $$
         \begin{array}{ccc}
          R_{\alpha_1\alpha_2}
           & \stackrel{\varphi_{\alpha_1\alpha_2}}{\longrightarrow}
           & R_{\alpha_2\alpha_1}                            \\[.6ex]
          \downarrow\hspace{1ex} & & \downarrow\hspace{1ex}  \\[-1.6ex]
          R^{\prime}_{\alpha^{\prime}_1\alpha^{\prime}_2}
           & \stackrel{\varphi^{\prime}
             _{\alpha^{\prime}_1\alpha^{\prime}_2}}{\longrightarrow}
           & R^{\prime}_{\alpha^{\prime}_2\alpha^{\prime}_1}
         \end{array}
        $$
        commute.
     \end{itemize}
  \end{itemize}
 We will denote ${\cal R}^{\prime}$, together with
   this data of arrows from ${\cal R}$ to ${\cal R}^{\prime}$,
   by ${\cal R}^{\prime}\precleftarrow {\cal R}$.
 Two gluing systems of rings ${\cal R}_1$ and ${\cal R}_2$ are said
   to be {\it equivalent}, in notation ${\cal R}_1\sim{\cal R}_2$,
  if there exists a gluing system ${\cal R}_3$ such that
   both ${\cal R}_3\precleftarrow{\cal R}_1$ and
   ${\cal R}_3\precleftarrow{\cal R}_2$ exist/hold.
 The {\it equivalence class} of ${\cal R}$ under refinements
  is denoted by $[{\cal R}]$.
} 

\bigskip

\noindent
{\bf Definition 1.2.6 [gluing system of ring-homomorphisms]\footnote{This
                                is a Grothendieck's
                                descent-data-of-morphisms description.}.}
{\rm
 A {\it gluing system of ring-homomorphisms}
  from a gluing system
  ${\cal R}=( \{R_{\alpha}\}_{\alpha\in A} \doublearrow
              \{R_{\alpha_1\alpha_2}\}_{\alpha_1,\alpha_2\in A} )$
  to another such system
  ${\cal S}=( \{S_{\beta}\}_{\beta\in B} \doublearrow
             \{S_{\beta_1\beta_2}\}_{\beta_1,\beta_2\in B} )$
  consists of the following data:
  \begin{itemize}
   \item[$\cdot$]
    a map $\tau: B\rightarrow A$ on the index sets;

   \item[$\cdot$]
    [{\it ring-homomorphisms on ring-charts}]\\[.6ex]  
    a collection
     $\{ \varphi_{\beta}: R_{\tau(\beta)}\rightarrow S_{\beta} \}
                                                       _{\beta\in B}$
     of ring-homomorphisms
    such that
    %
    %
    \begin{itemize}
     \item[-]
      [{\it compatibility with localizations}]\\[.6ex]  
      for all $\beta_1$, $\beta_2\in B$,
      $\varphi_{\beta_1}:R_{\tau(\beta_1)}\rightarrow S_{\beta_1}$
       is admissible
       and, hence, descends to a unique
       $\varphi_{\beta_1}|_{\beta_2}:R_{\tau(\beta_1)\tau(\beta_2)}
         \rightarrow S_{\beta_1\beta_2}$
       that makes the diagram
       $$
        \begin{array}{ccc}
          R_{\tau(\beta_1)}
           & \stackrel{\varphi_{\beta_1}}{\longrightarrow}
           & S_{\beta_1}                      \\[.6ex]
         \downarrow        && \downarrow      \\[-.6ex]
         R_{\tau(\beta_1)\tau(\beta_2)}
          & \stackrel{\varphi_{\beta_1}|_{\beta_2}}{\longrightarrow}
          & S_{\beta_1\beta_2}
        \end{array}
       $$
       commute, cf.\ Definition 1.2.2;

     \medskip
     \item[-]
      [{\it gluing conditions}]\\[.6ex]  
      the diagrams
      $$
       \begin{array}{ccc}
        S_{\beta_1\beta_2}  & \rightleftharpoons
                            & S_{\beta_2\beta_1}                \\[.6ex]
        \mbox{\scriptsize
              $\varphi_{\beta_1}|_{\beta_2}$}\uparrow\hspace{3em}
         && \hspace{2em}\uparrow
            \mbox{\scriptsize $\varphi_{\beta_2}|_{\beta_1}$}   \\
        R_{\tau(\beta_1)\tau(\beta_2)}
         & \rightleftharpoons  & R_{\tau(\beta_2)\tau(\beta_1)}
       \end{array}
      $$
      commute for all $(\beta_1,\beta_2)\in B\times B$.
    \end{itemize}
  \end{itemize}
 We will call the system $\Phi := (\tau,\{\varphi_{\beta}\}_{\beta})$
  also a {\it morphism} from ${\cal R}$ to ${\cal S}$.
} 

\bigskip

\noindent
{\bf Example 1.2.7 [refinement as a morphism].}
 A refinement ${\cal R}^{\prime} \precleftarrow {\cal R}$
  contains a system $\Phi:{\cal R}\rightarrow {\cal R}^{\prime}$
  of ring-homomorphisms in its data.
 In particular, a central cover $\{R\rightarrow R_{\alpha}\}_{\alpha}$
  of $R$ gives rise to a morphism
  $\{R\}\rightarrow \{R_{\alpha}\}_{\alpha}$.

\bigskip

Ring-homomorphisms have the following affine-gluing property:

\bigskip

\noindent
{\bf Lemma 1.2.8 [morphism: affine-gluing].}
{\it
 Given finitely generated rings $R$ and $S$,
 let
  $(\{S_{\alpha}\}_{\alpha\in A}
     \doublearrow \{S_{\alpha_1\alpha_2}\}_{\alpha_1,\alpha_2\in A})$
   be a gluing system of rings associated to a finite central cover
   $\{S\rightarrow S_{\alpha}\}_{\alpha\in A}$ of $S$  and
  $\Phi=\{ \varphi_{\alpha}: R \rightarrow S_{\alpha}\}_{\alpha\in A}$
   be a gluing system of ring-homomorphisms from $R$.
 Then, there exists a unique ring-homomorphism
  $\varphi:R\rightarrow S$ such that $\varphi$ descends to $\Phi$.
} 

\bigskip

\noindent
We will call $\varphi$ in the above lemma the {\it gluing}
 of the system $\Phi$.
A reverse of this lemma gives rise to the following definition:

\bigskip

\noindent
{\bf Definition 1.2.9 [refinement of morphism].} {\rm
 Given a morphism
   $\Phi=(\tau,\{\varphi_{\beta}\}_{\beta}):
                                   {\cal R}\rightarrow {\cal S}$  and
   a pair
    $( {\cal R}^{\prime}\precleftarrow {\cal R}\,,\,
       {\cal S}^{\prime}\precleftarrow {\cal S})$  of refinements,
   denote the index set of
    ${\cal R}$, ${\cal R}^{\prime}$, ${\cal S}$, ${\cal S}^{\prime}$
    by $A$, $A^{\prime}$, $B$, $B^{\prime}$ respectively.
  Let $\tau:B\rightarrow A$ and
      $( \tau_{A^{\prime},A}:A^{\prime}\rightarrow A\,,\,
         \tau_{B^{\prime},B}:B^{\prime}\rightarrow B )$
   be the maps on the index sets corresponding to $\Phi$ and
   the pair of refinements respectively.
 Then
  $( {\cal R}^{\prime}\precleftarrow {\cal R}\,,\,
     {\cal S}^{\prime}\precleftarrow {\cal S})$
   is said to be {\it $\Phi$-admissible}
  if, for all $\beta\in B$,
   $\varphi_{\beta}$ is admissible with respect to the
    localizations maps in the system pair
    $( {\cal R}^{\prime}\precleftarrow {\cal R}\,,\,
         {\cal S}^{\prime}\precleftarrow {\cal S})$;
   cf.\ Definition 1.2.2.
 When this is the case,
  fix a $\tau^{\prime}:B^{\prime}\rightarrow A^{\prime}$ so that
   the diagram
    $$
     \begin{array}{ccc}
      A & \stackrel{\tau}{\longleftarrow} & B  \\
      \hspace{-2em}\mbox{\scriptsize $\tau_{A^{\prime},A}$}\uparrow  &
       & \uparrow\mbox{\scriptsize $\tau_{B^{\prime},B}$}\hspace{-2.2em} \\
      A^{\prime}
       & \stackrel{\tau^{\prime}}{\longleftarrow} & B^{\prime}
     \end{array}
    $$
   commute.
 Then $\Phi$ descends to a unique morphism
  $\Phi^{\prime}
    =( \tau^{\prime},
       \{\varphi^{\prime}_{\beta^{\prime}}\}_{\beta^{\prime}} )
    : {\cal R}^{\prime}\rightarrow {\cal S}^{\prime}$,
  called a {\it refinement} of $\Phi$ with respect to
  $( {\cal R}^{\prime}\precleftarrow {\cal R}\,,\,
     {\cal S}^{\prime}\precleftarrow {\cal S})$.
} 

\bigskip

\noindent
{\bf Definition 1.2.10 [equivalence of morphisms].} {\rm
 Given equivalent ring-systems
   ${\cal R}_1\sim {\cal R}_2$ and ${\cal S}_1\sim {\cal S}_2$  and
  morphisms $\Phi_1:{\cal R}_1\rightarrow {\cal S}_1$ and
   $\Phi_2:{\cal R}_2\rightarrow {\cal S}_2$,
 we say that $\Phi_1$ and $\Phi_2$ are {\it equivalent},
   in notation $\Phi_1\sim \Phi_2$, if
 there exist common refinements
   ${\cal R}_1 \succrightarrow {\cal R}^{\prime}
                             \precleftarrow {\cal R}_2$  and
   ${\cal S}_1 \succrightarrow {\cal S}^{\prime}
                             \precleftarrow {\cal S}_2$
  such that
   (1)
    $({\cal R}^{\prime} \precleftarrow {\cal R}_1,
      {\cal S}^{\prime} \precleftarrow {\cal S}_1)$ and
    $({\cal R}^{\prime} \precleftarrow {\cal R}_2,
      {\cal S}^{\prime} \precleftarrow {\cal S}_2)$
    are $\Phi_1$- and $\Phi_2$-admissible respectively and
   (2)
    $\Phi_1$ and $\Phi_2$ can be descended to identical morphisms
    $\Phi_1^{\prime}=\Phi_2^{\prime}:
                   {\cal R}^{\prime} \rightarrow {\cal S}^{\prime}$.
 The equivalence class of $\Phi$ will be denoted by $[\Phi]$.
 An element in $[\Phi]$ will be called a {\it representative}
  of $[\Phi]$.
} 

\bigskip

\noindent
{\bf Definition 1.2.11 [strict morphism on equivalence classes].} {\rm
 By a {\it strict morphism} from $[{\cal R}_0]$ to $[{\cal S}_0]$,
 we mean an equivalence class $[\Phi:{\cal R}\rightarrow {\cal S}]$,
  where ${\cal R}\in [{\cal R}_0]$ and ${\cal S}\in [{\cal S}_0]$.
} 

\bigskip

By descending to a refinement ${\cal R}$ of ${\cal R}_0$ and
  taking the pre-composition with the localizations maps in
  ${\cal R}\precleftarrow {\cal R}_0$,
 one has the following lemma:

\bigskip

\noindent
{\bf Lemma 1.2.12 [one-side refinement enough].} {\it
 A strict morphism from $[{\cal R}_0]$ to $[{\cal S}_0]$
  can be represented by a $\Phi:{\cal R}_0\rightarrow {\cal S}$,
  for some ${\cal S}\in [{\cal S}_0]$.
} 

\bigskip

\noindent
Thus, in the discussion below, only the refinements on the
 $[{\cal S}_0]$-side are required.

%
%

\bigskip

\noindent
{\bf Definition 1.2.13 [injective strict morphism].} {\rm
 A {\it injective strict morphism}
  $[\Phi_0]:[{\cal R}]\rightarrow [{\cal S}_0]$
  is a strict morphism that can be represented by a
   $\Phi=(\tau,\{\varphi_{\beta}\}_{\beta})
                          :{\cal R}\rightarrow {\cal S}$,
   ${\cal S}\in [{\cal S}_0]$,
  such that
   (1) $\tau$ is surjective and
   (2) for each $R_{\alpha}\in {\cal R}$,
       there exists a $\beta\in \tau^{-1}(\alpha)$
        such that
        $\varphi_{\beta}:R_{\alpha}\rightarrow S_{\beta}\in {\cal S}$
        is a ring-monomorphism.
} 

\bigskip

\noindent
{\bf Definition 1.2.14 [(general) morphism].} {\rm
 A {\it general morphism} from $[{\cal R}]$ to $[{\cal S}_0]$
  consists of the following data:
  \begin{itemize}
   \item[$\cdot$]
     an injective strict morphism
      $[\Phi_0]: [{\cal S}_0^{\prime}]\rightarrow [{\cal S}_0]$,

   \item[$\cdot$]
    a strict morphism
     $[\Phi_0^{\prime}]:[{\cal R}]\rightarrow [{\cal S}_0^{\prime}]$.
  \end{itemize}
 We will denote
  the tuple $([{\cal S}_0^{\prime}], [\Phi_0], [\Phi_0^{\prime}])$
   collectively by $[\Phi_0^{\prime}]$  and
  a general morphism also by
  $[\Phi_0^{\prime}]:[{\cal R}]\rightarrow [{\cal S}_0]$.
 A {\it representative} of
  $[\Phi_0^{\prime}]:[{\cal R}]\rightarrow [{\cal S}_0]$
  is given by a {\it $3$-step ring-system-morphism diagram}
  $$
   {\cal R}\; \stackrel{\Phi^{\prime}}{\longrightarrow}\;
     {\cal S}^{\prime\prime}\;
     \precleftarrow\; {\cal S}^{\prime}\;
     \stackrel{\Phi}{\longrightarrow}\; {\cal S}
  $$
  with ${\cal S}\in [{\cal S}_0]$;
       ${\cal S}^{\prime}$,
        ${\cal S}^{\prime\prime}\in [{\cal S}_0^{\prime}]$;
       $\Phi\in [\Phi_0]$, and $\Phi^{\prime}\in [\Phi_0^{\prime}]$.
 A strict morphism is automatically a general morphism.
 A general morphism will also be called simply
  a {\it morphism}\footnote{For non-algebraic-geometers:
                            A few words follow on
                            why the {\it morphisms} in Sec.~1.1 and
                             here are defined as they are.
                            In the case of systems of
                             {\it commutative rings},
                             ``general morphism" is a redundant notion
                             as the 3-step ring-system-morphism diagram
                             ${\cal R}
                               \stackrel{\Phi^{\prime}}{\longrightarrow}
                                {\cal S}^{\prime\prime}
                               \precleftarrow {\cal S}^{\prime}
                               \stackrel{\Phi}{\longrightarrow} {\cal S}$
                             can always be reduced to a $2$-step diagram
                             ${\cal R}\;
                                \stackrel{\Phi}{\longrightarrow}\;
                                 {\cal S}^{\prime\prime\prime}\;
                                \precleftarrow {\cal S}$,
                             which represents a strict morphism.
                            In this case, $[{\cal R}]$ and $[{\cal S}]$
                             (resp.\ ${\cal R}$ and ${\cal S}$) are
                             contravariantly associated to schemes
                             (resp.\ atlases of affine charts on schemes).
                            This reducibility from a $3$-step diagram
                              to a $2$-step diagram no longer holds
                              in general
                              in the case of {\it noncommutative rings},
                             as the ring-homomorphisms $\varphi_{\beta}$
                              on ring-charts are required to be admissible
                              to the central localizations
                              in the construction
                              in order that gluings make sense and work.
                            On the other hand, when we shrink the rings
                             $S_{\beta}$
                             and take only a system of their subrings
                              $S^{\prime}_{\beta^{\prime}}$,
                             the center can increase:
                              $Z(S^{\prime}_{\beta^{\prime}})
                                               \supset Z(S_{\beta})$.
                            Thus, a ring-homomorphism that is not
                             admissible as a map to $S_{\beta}$ but with
                             the image contained in
                             $S^{\prime}_{\beta^{\prime}}$
                             may become admissible as a map to
                             $S^{\prime}_{\beta^{\prime}}$.
                            In other words, the notion of
                             general morphism partially takes care of
                             the more subtle issue of a functorial
                             construction of general localizations,
                             allowing us to stay in the much more
                             tractable central localizations.
                            This is not the whole story.
                            In the correspondence of the category of
                             {\it commutative rings} with the category
                             of (commutative) {\it affine schemes},
                             one has the canonical identification:
                             $\Mor(R,S)
                              =\Mor(\footnotesizeSpec S,
                                    \footnotesizeSpec R)$
                             by construction.
                            In the {\it noncommutative} case,
                             the functorial construction of the operation
                             ``\footnotesizeSpec" that associates
                              to a ring a ``space" is subtle.
                            Indeed, what Grothendieck's school accomplished
                             in the decade 1960s for
                             commutative algebraic geometry
                             is only partially realized through the work
                             of several independent schools on
                             noncommutative algebraic geometry
                             in the four decades after then.
                            There are several nonequivalent
                             constructions/realization
                             of the notion of ``\footnotesizeSpec",
                             with each maintaining part of the
                             equivalent characterizing properties of
                             $\footnotesizeSpec$
                             in the commutative case,
                            cf.\ sample references in Remark 1.1.6.
                            In the current work, we take rings and
                             ring-homomorphisms as more fundamental
                             for ``geometry" than the notion of
                             ``points" and ``topologies".
                            An injective strict morphism
                             $[\Phi_0]: [{\cal S}_0^{\prime}]
                                             \rightarrow [{\cal S}_0]$,
                             is then meant to give a dominant morphism
                             $\phi_0:  \footnotesizeSpace[{\cal S}_0]
                               \rightarrow
                               \footnotesizeSpace[{\cal S}_0^{\prime}]$,
                             should the latter ``spaces" be constructed
                              functorially.
                            Geometrically, a general morphism is then
                             simply an ordinary morphism precomposed
                             with a pinching and, hence, must be still
                             an allowable morphism if the setting
                             is natural.
                            From these hidden words to the main text,
                             one sees that we do want to include
                             general morphisms to
                             $\footnotesizeMor([{\cal R}],[{\cal S}_0])$
                             in any natural setting/definition.
                            Surprisingly, these independent
                             purely mathematical reasonings that
                             attempt to extend Grothendieck's language
                             of (commutative) algebraic geometry to
                             the noncommutative case give rise to
                             $\Mor([{\cal R}], [{\cal S}])$
                             that is also required for modeling
                             D-branes in string theory correctly!}.
 Define $\Mor([{\cal R}],[{\cal S}_0])$
  to be the set of morphisms from $[{\cal R}]$ to $[{\cal S}_0]$.
} 

\bigskip

\noindent
{\bf Example 1.2.15 [non-strict morphism].}
 Let $S$ be a subring of $S_0$ such that $Z(S)\supsetneqq Z(S_0)$ and
  $\Sigma=\{s_{\beta}\}_{\beta}$ be a finite subset in $Z(S)-Z(S_0)$
   such that $S=\sum_{s_{\beta}\in \Sigma}s_{\beta}\cdot S$.
 Let $S\rightarrow S_{\beta}$
  (resp.\ $S\rightarrow S_{\beta_1\beta_2}$)
  be the central localization with respect to $s_{\beta}$
  (resp.\ $s_{\beta_1}$ and then $s_{\beta_2}$) ,
 then $\{S\rightarrow S_{\beta}\}_{\beta}$ is a cover of $S$.
 Then the $3$-step diagram
  \begin{eqnarray*}
   \lefteqn{
    \left( \{S_{\beta}\}_{\beta}
     \doublearrow \{S_{\beta_1\beta_2}\}_{\beta_1, \beta_2} \right)
            }\\[.6ex]
   &&
   \stackrel{\scriptsizeId}{\longrightarrow}\;
   \left( \{S_{\beta}\}_{\beta}
     \doublearrow \{S_{\beta_1\beta_2}\}_{\beta_1, \beta_2} \right)\;
   \precleftarrow\;  \left( \{S\}\doublearrow \{S\} \right)\;
   \longrightarrow\; \left( \{S_0\}\doublearrow\{S_0\} \right)
  \end{eqnarray*}
  represents a morphism
   $[\Id]:( \{S_{\beta}\}_{\beta}
            \doublearrow\{S_{\beta_1\beta_2}\}_{\beta_1,\beta_2} )
                           \rightarrow (\{S_0\}\doublearrow \{S_0\})$
   that is not strict.
 See Sec.~4.2 for such examples with $S_0=M_n({\Bbb C})$.

\bigskip

\noindent
{\bf Grothendieck Ansatz [ring vs.\ space].} {\it
 We shall hiddenly think of
  an equivalence class $[{\cal R}]$ of ring-systems
   as a ``space" $\Space[\cal R]$
    with an equivalence class of atlases
     $\{\Space R_{\alpha}\}_{\alpha}$ (with the gluing data from
     the arrows
     $\{\Space R_{\alpha_1\alpha_2}\}_{\alpha_1,\alpha_2}
      \doublearrow \{\Space R_{\alpha}\}_{\alpha}$), and
  a morphism $[{\cal R}]\rightarrow [{\cal S}]$
   as a morphism $\Space[{\cal S}]\rightarrow \Space[{\cal R}]$.
  Cf.\ footnote 11.
} 

\bigskip

\noindent
{\it Remark 1.2.16 $[$morphism vs.\ map$]$.}
 In defining a {\it morphism} in a category of noncommutative spaces,
  we mean to keep both the {\it domain} and the {\it target} of
  the morphism {\it fixed}.
 In terms of the ring-system language,  this is reflected in
  the fact that
   a refinement of a ring-system ${\cal R}$
   is another ring system ${\cal R}^{\prime}$ {\it together with}
   a localization morphism ${\cal R}\rightarrow {\cal R}^{\prime}$ and
  the fact that the trivial localization is the identity map
   (not just a ring-isomorphism).
 In contrast, later (Sec.~4) when we discuss the space of {\it maps}
   or of {\it D0-brane probes},
  we remain to keep the target-space fixed but the {\it domain}-space
   will be taken as {\it not fixed}.
 The issue of automorphisms of the domain will then enter.

\bigskip

\noindent
{\it Remark 1.2.17 $[$$k$-algebra$]$.}
 When all the rings $R_{\alpha}\in {\cal R}$ involved are $k$-algebras
  for a fixed ground field $k$,
 we will take as a convention that all the ring-homomorphisms
  involved are then required to be $k$-algebra-homomorphisms
  unless otherwise noted.

\bigskip

\noindent
{\it Remark 1.2.18 $[$Azumaya-type noncommutative space$]$.}
 For an Azumaya-type noncommutative space
   $X=(X,{\cal O}_X,{\cal O}_X^{\nc})$,
  an affine cover $\{U_{\alpha}\}_{\alpha}$ of $(X,{\cal O}_X)$
  gives rise to a ring-system representation ${\cal R}_X$ of $X$
  defined by
  $$
   {\cal R}_X\;
   =\; (\{R_{\alpha}\}_{\alpha} \doublearrow
        \{R_{\alpha_1\alpha_2}\}_{\alpha_1,\alpha_2})\;
   :=\; (\{{\cal O}_X^{\nc}(U_{\alpha})\}_{\alpha} \doublearrow
         \{{\cal O}_X^{\nc}(U_{\alpha_1}\cap U_{\alpha_2})\}
                                          _{\alpha_1,\alpha_2})\,.
  $$
  Morphisms $X\rightarrow Y$
   between Azumaya-type noncommutative spaces can be expressed
   contravariantly as morphisms ${\cal R}_Y\rightarrow {\cal R}_X$
   of associated ring-systems.
  In particular,
   the notion of surrogates $X\rightarrow X^{\prime}$ of $X$
   corresponds to the notion of injective strict morphisms
   $[{\cal R}^{\prime}]\rightarrow [{\cal R}]$ into $[{\cal R}]$.

\bigskip

Before leaving this theme, we note that
 Lemma 1.2.8        
  and Lemma 1.2.12  
 together imply that:

\bigskip

\noindent
{\bf Lemma 1.2.19 [local description of morphisms].} {\it
 Let $R$ and $S$ be rings. Then
 $$
  \Mor(\Space[\{S\}], \Space[\{R\}])\;
  \stackrel{\mbox{\tiny Grothendieck Ansatz}}{:=}\;
  \Mor([\{R\}], [\{S\}])\; \simeq\; \Mor(R,S)
 $$
 canonically,
 where $\Mor(R,S)$ is the set of ring-homomorphisms from $R$ to $S$.
} 

\bigskip

\section{D-branes from the viewpoint of Grothendieck.}

\subsection{The notion of a space(-time): functor of points vs.\ probes.}

\begin{flushleft}
{\bf Space from a functor of points in algebraic geometry:
     a space without a space.}
\end{flushleft}
In the commutative case\footnote{Unfamiliar
                               readers are referred to [L-L-Y: Sec.~1]
                               for a brief introduction of and
                               literature guide for the notions of
                               {\it Grothendieck topology}, {\it site},
                               and {\it stack}.
                              All that is said here is standard
                               from algebraic deformation theory.},
let $\Scheme/S$ be the category of schemes over a base scheme $S$
 with a Grothendieck topology.
A {\it functor of points} on $\Scheme/S$ is a presheaf ${\cal F}$
 of sets on $\Scheme/S$.
For example, take $S=\Spec{\Bbb C}$, then a scheme $Y/{\Bbb C}$
 determines an ${\cal F}_Y$ on $\Scheme/{\Bbb C}$ with
 ${\cal F}_Y(Z) := \Mor_{{\scriptsizeBbb C}-\scriptsizescheme}(Z,Y)$
  for $Z\in \Scheme/{\Bbb C}$.
In this case, $Y$ can be recovered from ${\cal F}_Y$,
 cf.\ Yoneda lemma.

One can think of a functor of points ${\cal F}$
  as a generalized space ${\frak Y}_{\cal F}$ and
 ${\cal F}(Z)$ as the set $\Mor(Z,{\frak Y}_{\cal F})$
  of $Z$-valued points on ${\frak Y}_{\cal F}$.
The construction of the moduli space that satisfies
 the functorial/universal property for a moduli problem
 leads one in general to such a generalized space.
Encoded in the functor of points ${\cal F}$ on $\Scheme/S$ is
 the data of extension property of morphisms into ${\frak Y}_{\cal F}$.
In particular, ${\cal F}$ contains the information of
 tangent-obstruction structure of ${\frak Y}_{\cal F}$
 as well as of local properties like smoothness
 at a point (i.e.\ an element in, e.g.,
  ${\cal F}(\Spec{\Bbb C})=:\Mor(\Spec{\Bbb C},{\frak Y}_{\cal F})$)
 of ${\frak F}_{\cal F}$.
It is in this way that ${\cal F}$ describes the geometry of
 a ``space" {\it without} giving the space beforehand,
 for example, as a point-set with a topology and other structures.
Schemes, Deligne-Mumford stack (i.e.\ orbifolds), Artin stacks,
 and many moduli functors are all examples of functors of points.

There are diverse ways/versions to generalize the above to
 the noncommutative case.
The particular one that is selected from Sec.~1.2 is to consider
 the category $\RingSystemCategory$ of gluing systems of rings
 with a Grothendieck topology defined by central covers,
 \'{e}tale central covers, or fppf central covers.
(The \'{e}tale or fppf condition of a morphism can be defined
  purely ring-theoretically.)
Note that, as we are dealing directly with rings, all the arrows
  in the commutative case above are reversed here.
(However, if one wishes, one may write a ring system ${\cal R}$
 by a formal symbol $\Space {\cal R}$, meaning the associated
  space/geometry to ${\cal R}$, to preserve all the arrow directions.)
A {\it functor of points} on $\RingSystemCategory$ is then
 a presheaf ${\cal F}$ of sets on $\RingSystemCategory$.
Again, one can directly think of ${\cal F}$ as a generalized
 noncommutative space ${\frak Y}_{\cal F}$.
The data of extension properties of morphisms to ${\frak Y}_{\cal F}$
 is encoded in ${\cal F}$. Through this, ${\cal F}$ describes
 the geometry of a generalized noncommutative space ${\frak Y}_{\cal F}$
 without ${\frak Y}_{\cal F}$ being given beforehand.

\bigskip

\begin{flushleft}
{\bf Space(-time) from probes in QFT/string theory:
     space(-time)s emerge from QFT.}
\end{flushleft}
%
%
There are two particular classes of quantum field theories (QFT's)
 that are directly relevant to the notion of target space(-time):
\begin{itemize}
 \item[$\cdot$]
  {\it Nonlinear sigma models} are, by definition, quantum field
  theories whose field contents contain, among other fields,
  bosonic fields corresponding to maps from a domain
  (cf. world-volume of branes) to a target space(-time).

 \item[$\cdot$]
  In string theory, D0-brane physics is described by matrix theory.
  As the moduli space of a single D0-brane moving in a space(-time)
   is the space(-time) itself, the {\it moduli space of a single
   D0-brane} can be identified as the target space(-time).
 \end{itemize}

These two concretely target-space(-time)-related situations can
 be hidden implicitly in a general quantum field theory
 that is seemingly irrelevant to a target space-time.
Furthermore,
 depending on where we look at the theory in the related
  Wilson's theory-space\footnote{See
                         [L-Y1: appendix A.1] for highlights and
                         a literature guide for mathematicians on
                         this very important notion from
                         quantum field theory.
                        In particular, a {\it Wilson's theory-space}
                         goes with {\it universal objects} over it
                         that encode QFT contents, and
                        a {\it duality} is a local isomorphism
                         on Wilson's theory-space with
                         these structures.},
 there can be more than one target space(-time)s hidden in
  one combinatorial class of quantum field theories.
Even more,
 such target spaces can be taken
  either at the classical level -- which usually involve only
   algebraic manipulations of the Lagrangian of the theory --
  or at the quantum level -- which has to bring in the core
   techniques (and some arts as well) from quantum field theory.
A quantum-corrected target space(-time) can be different from
 its associate classical target space(-time).
The following three examples have been around for a while
 in the string-theory community:

\bigskip

\noindent
{\bf Example 2.1.1 [gauged linear sigma model].}
{\it Geometric phases} of a gauged linear sigma model are realized
 effectively by nonlinear sigma models.
Birationally equivalent target spaces emerge.
See [Wi1] and [M-P].

\bigskip

\noindent
{\bf Example 2.1.2 [D0-brane probe of space(-time) and singularities].}
A D0-brane moving in a singular space(-time) recognizes various
 (partial) resolutions of the singular space(-time) as
 the {\it moduli space of D0-branes} at different phases
 in the Wilson's theory-space of the $0+1$ dimensional matrix theory
 involved.
Birationally equivalent smooth or partially resolved target spaces
 emerge from a single singular target space.
See [D-G-M], [Do-M], and [G-L-R].

\bigskip

\noindent
{\bf Example 2.1.3 [conformal field theory with boundary].}
D0-branes are realized in a conformal field theory with boundary
 as a special class of boundary states.
The {\it moduli space of} such {\it boundary states} gives rise to
 a target space(-time).
See [M-M-S-S] and [S-S].

\bigskip

These examples suggest that quantum field theories,
  as probes to a target space(-time),
 can be more fundamental than the space(-time) itself.
The latter may even lose its absolute meaning under dualities of
 quantum field theories, like what happens in mirror symmetry.

\bigskip

\begin{flushleft}
{\bf Functor of points vs.\ probes.}
\end{flushleft}
A comparison of these two notions is given below:

\noindent\hspace{1em}
\begin{tabular}{lcl}
  & \hspace{1em} & \\
 \hspace{1.6em}
 {\it functor of points}
  &&  \hspace{7.6em} {\it QFT's as probes}  \\[.6ex]
 \hline  \\[-1.6ex]
 $\cdot$
  $\Scheme/S$
  && a category $\Brane$ of branes\\[1ex]
 $\cdot$
 a functor of point ${\cal F}$
  && a compatible system
     $\{\mbox{QFT}_{\Sigma}\}_{\Sigma\in\scriptsizeBrane}$
     of effective QFT\\
 \hspace{1.5ex}on $\Scheme/S$
  && on branes that have isomorphic target space(-time)s\\[1ex]
 $\cdot$
 ${\cal F}(T)\,$, $\,T\in \Scheme/S$
  && bosonic fields on a brane that correspond to maps\\
  && from the brane to a target space(-time) \\
  && \\
\end{tabular}

\noindent
Here, a `brane' means the defining domain of a quantum field theory.
For example, it can be the world-volume of a string, a D-brane,
 or an NS-brane.
Note also that
 a functor of points ${\cal F}$ encodes the data of a space
while an effective QFT from a QFT as a probe encodes more than
 just the information of the target space(-time).

\bigskip

\subsection{D-branes as Azumaya-type noncommutative spaces.}

\begin{flushleft}
{\bf Question: What is a D-brane intrinsically?}
\end{flushleft}
A {\it D-brane} (in full name: {\it Dirichlet brane} or
 {\it Dirichlet membrane})
 in string theory is by definition (i.e.\ by the very word `Dirichlet')
 a boundary condition for the end-points of open strings
 moving in a space-time.
In the geometric/target-space-time aspect\footnote{It
                                 should be noted that there are also
                                 algebraic properties of D-branes
                                 realized as states or operators
                                 in a $2$-dimensional conformal
                                 field theory with boundary.
                                These algebraic properties from
                                 the open-string world-sheet
                                 perspective reflect the geometric
                                 properties of D-branes in the
                                 target space-time of strings.
                                Our focus in this work is on
                                 the geometric aspect as given
                                 in [Pol3] and [Pol4].},
one may start by thinking of the world-volume
 (cf.\ Remark/Definition 2.2.4) of a D-brane
 as an embedded submanifold $f: Z\hookrightarrow M$
 in an open-string target space-time $M$
 such that:
 \begin{itemize}
  \item[$\cdot$]
   {\bf [defining property of D-brane: D = Dirichlet]}\\[.6ex]
   {\it The boundary of open-string world-sheets are mapped to $f(Z)$
        in $M$.}
 \end{itemize}
Via this defining property, open strings induce then additional
 structures on $Z$, including
  a gauge field
   (from the vibrations of open-strings with end-points on $f(Z)$) and
  a Chan-Paton bundle
   (from the Chan-Paton index on the end-points of such an open string)
  on $Z$.
Basic properties of D-branes under such a setting are given in
 [Pol3] and [Pol4].

To bring the relevant part of the work of Polchinski into
 the discussions and to enable a direct comparison/referral,
let us introduce notations the-same-as/as-close-as-possible-to
 those in [Pol4: vol.~I, Sec.~8.7]:
 let $\xi:=(\xi^a)_a$ be local coordinates on $Z$ and
  $X:=(X^a;X^{\mu})_{a,\mu}$ be local coordinates on $M$
  such that the embedding $f:Z\hookrightarrow M$ is locally
  expressed as
  $$
   X\; =\; X(\xi)\; =\; (X^a(\xi); X^{\mu}(\xi))_{a,\mu}\;
   =\; (\xi^a, X^{\mu}(\xi))_{a,\mu}\,;
  $$
 i.e., $X^a$'s (resp.\ $X^{\mu}$'s) are local coordinates along
       (resp.\ transverse to) $f(Z)$ in $M$.
This choice of local coordinates removes redundant degrees of freedom
 of the map $f$, and
$X^{\mu}=X^{\mu}(\xi)$ can be regarded as (scalar) fields on $Z$
 that collectively describes the postions/shapes/fluctuations
 of $Z$ in $M$ locally.
Here, both $\xi^a$'s, $X^a$'s, and $X^{\mu}$'s are ${\Bbb R}$-valued.
The gauge field on $Z$ is locally given by the connection
 $1$-form $A=\sum_a A_a(\xi)d\xi^a$ of a $U(1)$-bundle on $Z$.

When $n$-many such D-branes $Z$ are coincident, from the associated
 massless spectrum of (oriented) open strings with both end-points
 on $f(Z)$ one can draw the conclusion that
 \begin{itemize}
  \item[(1)]
   The gauge field $A=\sum_a A_a(\xi)d\xi^a$ on $Z$ is enhanced to
    $u(n)$-valued.

  \item[(2)]
   Each scalar field $X^{\mu}(\xi)$ on $Z$ is also enhanced
    to matrix-valued, cf.\ footnote~17.
 \end{itemize}
Property (1) says that there is now a $U(n)$-bundle on $Z$.
But
 \begin{itemize}
  \item[$\cdot$]
   {\bf Q.}\ {\it What is the meaning of Property (2)?}
 \end{itemize}
For this, Polchinski remarks that:
\begin{itemize}
 \item[$\cdot$]
 [{\sl quote from} [Pol4: vol.~I, Sec.~8.7, p.272]]\hspace{1em}
 ``{\it
  For the collective coordinate $X^{\mu}$, however, the meaning
   is mysterious: the collective coordinates for the embedding of
   $n$ D-branes in space-time are now enlarged to $n\times n$ matrices.
  This `noncommutative geometry' has proven to play a key role in
   the dynamics of D-branes, and there are conjectures that
   it is an important hint about the nature of space-time.}"
\end{itemize}
Particularly from the mathematical/geometric perspective,
Property (2) of D-branes when they are coincident,
  the above question, and Polchinski's remark
 are more appropriately incorporated into the following
 guiding question:
 \begin{itemize}
  \item[$\cdot$]
  {\bf Q.\ [D-brane]}$\;$
  {\it What is a D-brane intrinsically?}
 \end{itemize}
In other words, what is the {\it intrinsic} definition of D-branes
 so that {\it by itself} it can produce the properties of D-branes
 (e.g.\ Property (1) and Property (2) above)
 that are consistent with, governed by, or originally produced by
 open strings as well?\footnote{Since
                     the work of Ramond and of Neveu and Schwarz in 1971
                     that initiated string theory, there are by now
                     at least three ways to enter superstring theory:
                      Gate (1) the string-world-sheet/CFT way ($d=1+1$
                       or $d=2$ theory),
                      Gate (2) the target-space-time/supergravity/soliton
                       way ($d=9+1$ or $d=10+1$ theory),  and
                      Gate (3) the matrix-theory way ($d=0+1$ theory).
                    In Gate (1), after Wick-rotation, one can have
                     Riemann surfaces, conformal field theories,
                     moduli space of Riemann surfaces, ..., etc.\
                     before asking how strings move in a space-time.
                    D-branes entered string theory in the second half
                     of 1980s and took a central role after 1995 mainly
                     from the development of Gate (2) during 1990 - 1995.
                    In asking this question, we mean also to repeat
                     Gate (1) but for D-branes instead of for strings.
                    In other words, we are taking a ``D-brane"
                     as a fundamental object and asking,
                     ``What is (the definition of) a D-brane?",
                     before addressing how they ``move" in
                     -- i.e.\ are mapped into -- a space-time.}

\bigskip

\begin{flushleft}
{\bf The noncommutativity ansatz: from Polchinski to Grothendieck.}
\end{flushleft}
To understand Property (2) of D-branes, one has two aspects that
 are dual to each other:
 \begin{itemize}
  \item[(A1)]
   [{\it coordinate tuple as point}]\hspace{1em}
   A tuple $(\xi^a)_a$ (resp.\ $(X^a; X^{\mu})_{a,\mu}$)
    represents a point on the world-volume $Z$ of the D-brane
    (resp.\ on the target space-time $M$).

  \item[(A2)]
   [{\it local coordinates as generating set of
         local functions}]\hspace{1em}
   Each local coordinate $\xi^a$ of $Z$ (resp.\ $X^a$, $X^{\mu}$ of $M$)
    is a local function on $Z$ (resp.\ on $M$)  and
   the local coordinates $\xi^a$'s
    (resp.\ $X^a$'s and $X^{\mu}$'s) together
    form a generating set of local functions on the world-volume $Z$
    of the D-brane (resp.\ on the target space-time $M$).
 \end{itemize}
While Aspect (A1) leads one to the anticipation of a noncommutative
  space from a noncommutatization of the target space-time $M$
  when probed by coincident D-branes,
 Aspect (A2) of Grothendieck leads one to a different/dual\footnote{In
                                       what precise sense
                                        the noncommutativity of target
                                         space-time  and
                                        the noncommutativity of
                                         world-volume of branes
                                       are dual to each other
                                       deserves more thoughts.}
   conclusion:
  a noncommutative space from a noncommutatization of
  the world-volume $Z$ of coincident D-branes,
 as follows.

Denote by ${\Bbb R}\langle \xi^a\rangle_{a}$
  (resp.\ ${\Bbb R}\langle X^a; X^{\mu}\rangle_{a, \mu}$)
 the local function ring on the associated local coordinate chart
 on $Z$ (resp.\ on $M$).
Then the embedding $f:Z\rightarrow M$,
  locally expressed as
  $X=X(\xi)=(X^a(\xi); X^{\mu}(\xi))_{a,\mu}=(\xi^a; X^{\mu}(\xi))$,
 is locally contravariantly equivalent to a ring-homomorphism
 $$
  f^{\sharp}\;:\;
   {\Bbb R}\langle X^a; X^{\mu}\rangle_{a, \mu}\;
   \longrightarrow\; {\Bbb R}\langle \xi^a\rangle_{a}\,,
  \hspace{1em}\mbox{generated by}\hspace{1em}
  X^a\;\longmapsto\; \xi^a\,,\;
  X^{\mu}\;\longmapsto\;X^{\mu}(\xi)\,.
 $$
When $n$-many such D-branes are coincident, $X^{\mu}(\xi)$'s become
 $M_n({\Bbb C})$-valued.\footnote{Strictly
                            as induced by open-strings,
                            $X^{\mu}(\xi)'s$ are $u(n)$-valued
                             for oriented open strings and
                              either {\it so}$(n)$-
                               or {\it sp}$(n/2)$-valued
                             for unoriented open strings.
                           Instead of any of these Lie algebras,
                            here we directly think of $X^{\mu}(\xi)$
                            as $M_n(\footnotesizeBbb C)$-valued,
                             where $M_n(\footnotesizeBbb C)$
                              is regarded not as a Lie algebra with
                              a bracket (i.e.\ Lie product)
                              but rather as an associative algebra
                              (from the matrix multiplication)
                              with an identity ${\mathbf 1}$,
                            for two reasons:
                            \begin{itemize}
                             \item[(1)]
                              These Lie algebras are not associative
                               nor with an identity with respect to
                               the Lie product.
                              This makes the notion of localizations
                               and covers, which are crucial in algebraic
                               geometry for the local-to-global setup,
                               difficult to implement.
                              In view of noncommutative algebraic
                               geometry over ${\footnotesizeBbb C}$,
                               it is more natural to think of
                               $X^{\mu}(\xi)$'s as in a special class
                               of $M_n(\footnotesizeBbb C)$-valued
                               functions with $M_n({\footnotesizeBbb C})$
                               as an associative algebra with an identity.
                              Any associative algebra defines also
                               a tautological Lie algebra,
                               with the Lie product
                               $[x,y] := x\cdot y - y\cdot x$.
                              One can use this to translate back
                               to Lie algebras whenever needed.

                             \item[(2)]
                              In seeking the intrinsic
                               definition/structure of
                               a D-brane (or D-brane world-volume),
                              it is more natural to select the
                               structures thereon
                               as encompassing/universal as possible
                               so that they contain all what different
                               types of open strings can detect/see.
                              Each specific sector of structures on
                               D-brane world-volume seen by
                               a particular type of open strings
                               is then realized by a reduction from
                               the universal structures on D-brane
                               world-volume, as in the reductions of
                               the structure group of principal
                               {\it GL}$_n$ fiber bundles.
                           \end{itemize}
                           Cf.\ footnote~21.}
Thus, $f^{\sharp}$ is promoted to a new local ring-homomorphism:
 $$
  \hat{f}^{\sharp}\;:\;
   {\Bbb R}\langle X^a; X^{\mu}\rangle_{a, \mu}\;
   \longrightarrow\; M_n({\Bbb C}\langle \xi^a\rangle_{a})\,,
  \hspace{1em}\mbox{generated by}\hspace{1em}
  X^a\;\longmapsto\; \xi^a\cdot{\mathbf 1}\,,\;
  X^{\mu}\;\longmapsto\;X^{\mu}(\xi)\,.
 $$
Under Grothendieck's contravariant local equivalence of function rings
 and spaces, $\hat{f}^{\sharp}$ is equivalent to saying that we have
 now a map $\hat{f}: Z_{\scriptsizenoncommutative}\rightarrow M$.
Thus, the result of Polchinski re-read from the viewpoint of
 Grothendieck implies the following ansatz:

\bigskip

\noindent
{\bf Polchinski-Grothendieck Ansatz [D-brane: noncommutativity].}
{\it
 The world-volume of a D-brane carries a noncommutative structure
 locally associated to a function ring of the form $M_n(R)$
 for some $n\in {\Bbb Z}_{\ge 1}$ and ring $R$.\footnote{On
                                   purely mathematical ground,
                                  the $M_n(R)$ in the ansatz can
                                   be generalized in some cases.
                                  For example,
                                   in the case that $R$ is a Noetherian
                                    (commutative) integral domain,
                                   $M_n(R)$ can be replaced by
                                    the more general notion of
                                    an {\it $R$-order in a central
                                            simple $Q_R$-algebra},
                                    where $Q_R$ is the quotient field
                                    of $R$;
                                   cf.\ [Re].}
} 

\bigskip

\noindent
This ansatz is further enforced if one recalls that scalar fields on
 the world-volume of a brane are supposed to come from elements in
 the function ring of that world-volume  and the comparison of
 a functor of points vs.\ probes in Sec.~2.1.\footnote{From C.-H.L:
                    Several teachers and colleagues influenced
                     my painfully slow realization/appreciation
                     of this ansatz and its importance
                     through the personal journey of string theory:
                    {\it Orlando Alvarez} brought me to
                     the beauty of string theory and T-duality
                     at the dawn of its second revolution.
                    {\it Rafael Nepomechie} shared with me his
                     experience in the early days of higher-dimensional
                      extended objects before they became dominating
                     in ``string theory".
                    {\it Pei-Ming Ho} communicated the work
                     [Ho-W] to me.
                    The group meetings of the school of
                     {\it Philip Candelas} and
                     the insightful debates between
                      {\it Jacques Distler} and {\it Vadim Kaplunovsky}
                      promoted my understandings and kept me
                     aware of subtleties as well.
                    Teaching the late Professor {\it Raoul Bott}
                     mirror symmetry, fall 2000, assigned by
                     {\it Shing-Tung Yau} gave me a rare chance to
                      slow down and to map out what I had still been
                      ignorant of in the big picture.
                    The heat and enthusiasm {\it Shiraz Minwalla}
                     brought in to his various topic courses
                     from field theory to strings,
                     from phase structures in QFT to
                     supersymmetry, $\cdots\,$ over the years
                     helped me to access the mind of physicists
                     at the frontier.
                    {\it Shinobu Hosono} explained [H-S-T] to me in
                     March 2002, in which the subtle issue of
                      the multiplicity/wrapping of D-branes
                      in the torsion-sheaf picture
                      was brought out among other things.
                    Discussions with {\it Mihnea Popa}, spring 2002,
                     and his joint Seminar on Derived Category
                     with {\it Mircea Mustata}, fall 2002,
                     influenced my mathematical understanding
                     of D-branes of B-type.
                    The semester-long communications with
                     {\it Barton Zwiebach} on the draft of [Zw],
                     spring 2003, improved my understanding
                     of the physical fundamentals of string theory.
                    {\it Paul Aspinwall} emphasized many subtleties
                     of D-branes in his lectures at TASI 2003.
                    The topic courses and talks of {\it Kentaro Hori},
                     {\it Andrew Strominger}, and {\it Cumrun Vafa}
                     on string theory over the years printed
                     in my mind various pictures of how, physicists
                     think, D-branes should function.
                    The daily summary of work to each other with
                     {\it Ling-Miao Chou} over the years helped
                     to clarify my thoughts.
                    The vanishing lemma derived in [L-Y3] and its
                     comparison with [D-F] led me to a train of
                     discussions with {\it Duiliu-Emanuel Diaconescu},
                     December 2006, on the meaning of open-string
                     world-sheet instantons in the open/closed string
                     duality. 
                    These discussions propelled me to come back to
                     re-think about D-brane theory as a companion
                     theory to   
                     topological open strings and their instantons,
                     particularly the virtual ones.
                    Finally, it should be noted that,
                     even with this ansatz,
                     there are still other things missing
                     mathematically to understand D-branes fully
                     in a larger scope, cf.\ footnote~20.

                    Incidentally, while this work is under writing,
                     {\it William Thurston} came to give a talk,
                     May 2007, on the future of $3$-dimensional geometry
                     and topology after the justification of
                     the geometrization conjecture of $3$-manifolds.
                    Hyperbolic geometry has now applications
                     to cosmology and 
                     AdS/CFT correspondence.
                    It is surprising how a change of course of life
                     of a teacher can lead to a completely unexpected
                     journey of his student.
                    This detour  
                     is very demanding, yet only particularly lucky
                     one is given a chance to it.}

\bigskip

\noindent
{\it Remark 2.2.1 $[$D-brane and noncommutative geometry$]$.} {\rm
 The observation that D-brane should be related to noncommutative
  geometry was made soon after the second-revolution year 1995
  of string theory;
 see [Dou4] and [Dou5] for a survey and,
  e.g., [Ho-W] for an earlier study and
        [Laz] for a more recent study
  in the differential/symplectic geometry category.
 Noncommutative structures on a D-brane itself and on a space-time
  are two related but separate issues, e.g.\ [Dou2], [C-H1], and [C-H2].
 It is worth pointing out that,
  {\it from the viewpoint of Grothendieck,
   it is the noncommutative structure on the world-volume of
    a D-brane that comes first}.
 It is exactly because of such a structure on D-branes that
  a space-time may reveal its noncommutative nature
  when probed by a D-brane.
 Said algebro-geometrically in terms of function rings,
  since a ring-homomorphism from a noncommutative ring $R$
   to a commutative ring $S$ must factor through a ring-homomorphism
   $R/[R,R]\rightarrow S$ from the commutatization $R/[R,R]$ of $R$,
  D-branes without a noncommutative structure thereon cannot
   probe/sense any noncommutativity, if any, of a space-time at all.
} 

\bigskip

\noindent
{\it Remark 2.2.2 $[$B-field and noncommutativity on D-brane$]$.} {\rm
 It is known that when the target space(-time) $M$ has the B-field
  $B$ turned on, the gauge theory on a D-brane world-volume $Z$
  can be expressed as a noncommutative gauge theory;
  (see [Ch-K] and [S-W2]
   for details and more references on this subject.)
 From the underlying formulation, this implies in particular that,
  in this case, the commutative product of a local function ring $R$
  on $Z$ is deformed to a noncommutative $\ast$-product
  depending on $B$.
 When $n$-many D-branes $Z$ coincide,
  these string-induced property on D-branes
   compared with our discussion above says that:
  \begin{itemize}
   \item[$\cdot$]
    If $B=0$, then a local function ring on the world-volume of
     the coincident D-branes is of the form $M_n(R)$,
     where $R$ is commutative.

   \item[$\cdot$]
    If $B\ne 0$, then a local function ring on the world-volume of
     the coincident D-branes can become $M_n(R_B)$,
      where $R_B$ is a noncommutatization of $R$
      depending on/induced by $B$.
  \end{itemize}
In this work, we ignore the effect of B-field.
} 

\bigskip

{\it The Polchinski-Grothendieck Ansatz for D-branes applies
 to both nonsupersymmetric and supersymmetric D-branes,  and
 to both D-branes of A-type and D-branes of B-type}
  (cf.\ [B-B-St], [H-I-V], and [O-O-Y]) {\it in the latter case.}
Due to the different languages
 used in differential geometry and in algebraic geometry
 for noncommutative geometry
 (though the philosophy to equate a space and a function ring
  in each category is common),
we will focus entirely on supersymmetric D-branes of B-type,
 for which algebro-geometric language is appropriate.
The ansatz leads thus to a prototype\footnote{For
                         non-string-theorists:
                         There are two reasons we call this
                          a ``prototype" definition.
                         The first one is mild: we focus only on
                          the most essential fields on the brane and
                          ignore the others.
                         The second one is the true reason:
                          the definition we give here reflects only
                          what one should think mathematically about
                          a D-brane {\it in a special region of the
                          relevant Wilson's theory-space of string theory}
                          (cf.\ [L-Y1: appendix A.1]) and,
                         furthermore, we ignore also here the variation
                          to the definition required to incorporate
                          all forms of {\it D-brane bound states}.
                         Once we move away from this region, what one
                          should think of D-branes can become more
                          complicated or even not that clear when
                          trying to incorporate both mathematics and
                          physics involved.
                         However, since the mathematical definition
                          given here naturally reproduces the key
                          features of D-branes in its beginning years
                          after Polchinski [Pol2], it is our strong belief
                          that those more involved and languagewise
                          more demanding features/descriptions of D-branes
                          by string theorists in its growing years
                          can finally be reached, beginning with the
                          current prototype definition.
                         While the detail of this advanced step remains
                          challenging, there is definitely a related
                          Floer-Gromov-Witten-type theory involved
                          so that the coupling of D-branes and strings
                          is always incorporated, cf.\ footnote 1.}
 intrinsic definition of D-branes of B-type as follows:

\bigskip

\noindent
{\bf Definition 2.2.3 [D-brane of B-type and Chan-Paton sheaf].}
{\rm
 (1)
 A {\it D-brane of B-type} is an Azumaya-type noncommutative space
   $(X,{\cal O}_X, {\cal O}_X^{\nc})$ over ${\Bbb C}$,
  together with a fundamental ${\cal O}_X^{\nc}$-module ${\cal E}_X$.
 ${\cal E}_X$ is called the {\it Chan-Paton sheaf} on the D-brane $X$.
 We say that ${\cal E}_X$ has {\it rank} $r$ if it has rank $r$
  as an ${\cal O}_X$-module.
 Note that
  ${\cal E}_X|_{\eta} \simeq
     \kappa_{\eta}^{n_1}\oplus\,\cdots\,\oplus\kappa_{\eta}^{n_s}$
  at a generic point
   $\eta$ of $(X,{\cal O}_X)$ with residue field $\kappa_{\eta}$
  if ${\cal O}_X^{\nc}|_{\eta}/J({\cal O}_X^{\nc}|_{\eta})
       \simeq  M_{n_1}(\kappa_{\eta})
                  \times\,\cdots\,\times M_{n_s}(\kappa_{\eta})$.
  Here
   ${\cal O}_{X}^{\nc}|_{\eta}$ is the fiber of ${\cal O}_X^{\nc}$
    at $\eta$  and
   $J(\,\cdot\,)$ is the Jacobson radical of $(\,\cdot\,)$.
 \hspace{1ex}(2)
 A {\it D-brane} (of B-type) {\it in a target space $Y$}
  is a morphism $\Phi:X\rightarrow Y$.
 Here, $Y$ can be a (commutative) scheme,
  an Azumaya-type noncommutative space,
  a noncommutative space represented by a ring-system, or
  whatever noncommutative space
   to which the notion of morphisms from $X$ can be defined.
 The image Azumaya-type noncommutative space $\Phi(X)$ is called
  the {\it image D-brane} of $X$ in $Y$.
 \hspace{1ex}(3)
 The {\it Chan-Paton sheaf of a D-brane $\Phi:X\rightarrow Y$ on $Y$}
  is the push-forward $\Phi_{\ast}{\cal E}_X$ of ${\cal E}_X$,
  a coherent sheaf supported on $\Phi(X)$ in $Y$.
} 

\bigskip

\noindent
{\bf Remark/Definition 2.2.4 [D-brane vs.\ D-brane world-volume].}
{\rm
 The world-volume of a D-brane is what a D-brane sweeps out
  in a space-time and, hence, has the extra time-dimension than
  the D-brane has.
 It has a Lorentzian structure by definition.
 The world-volume after Wick rotation is called a Euclidean
  D-brane world-volume, which has now a Riemannian structure.
 We will define a {\it Eulcidean D-brane world-volume of B-type}
  the same as in Definition 2.2.3
  with `D-brane' replaced by `Euclidean D-brane world-volume'.
 Similarly, for a {\it Euclidean D-brane world-volume} (of B-type)
  {\it in a target space} $Y$ and the {\it Chan-Paton sheaf}
  and {\it its push-forward} on $Y$.
 In general, we keep the word `Euclidean' implicit and call it simply
  {\it D-brane world-volume} (of B-type)
   (resp.\ {\it D-brane world-volume} (of B-type) {\it in $Y$}).
 Readers should compare these simplified terminologies with the term
  `world-sheet' in the commonly used statement by physicists:
  ``The world-sheet of a string is a Riemann surface.",
  which takes the same interpretation implicitly.
} 

\bigskip

How these two definitions fit in string theory and, by themselves,
 reproduce three key open-string-induced properties of D-branes
 can be summarized/highlighted as follows:
 \begin{itemize}
  \item[(1)] [{\it interaction with open strings}]
   \begin{itemize}
    \item[$\cdot$]
    The Chan-Paton sheaf ${\cal E}_X$ should be identified
     with a singular coherent analytic sheaf on $X$ with
      a (singular) connection $A$
     via a Kobayashi-Hitchin correspondence.
    An end-point of an open string in $Y$ can then be coupled to
     the D-brane $X$ via a morphism $\Phi:X\rightarrow Y$ and
     the connection $A$, regarded as on ${\cal E}_X$.
   \end{itemize}
 \end{itemize}
 \begin{itemize}
  \item[(2)] [{\it source of Ramond-Ramond fields}]
   \begin{itemize}
    \item[$\cdot$]
     (Subject to that $X$ here has to be interpreted
      as a Euclidean D-brane world-volume.)
     Identify $(X, {\cal O}_X)$ canonically with an analytic space
      $X_{\an}$ (with the structure sheaf ${\cal O}_{X_{\an}}$
      of analytic functions).
     A Ramond-Ramond field (i.e.\ a differential form) on $Y$
      can be pulled back and integrate over $X_{\an}$
      via $\Phi:X\rightarrow Y$.
   \end{itemize}
 \end{itemize}
 \begin{itemize}
  \item[(3)] [{\it Higgsing/un-Higgsing associated to
                   un-stacking/stacking of D-brane}]
  \item[]
   The Azumaya-type noncommutative structure ${\cal O}_X^{\nc}$
    on $X$ makes the deformations of $\Phi:X\rightarrow Y$
    locally matrix-valued, as in [Pol4].
   It realizes the {\it Higgsing/un-Higgsing} behavior of
    the gauge theory on D-branes on $Y$ via (a continuous family of)
    deformations of a morphism $\Phi:X\rightarrow Y$,
    as explained below:

  \item[]
  \begin{itemize}
   \item[(3.1)]
    Associated to the (associative, unital) ${\cal O}_X$-algebra
     ${\cal O}_X^{\nc}$ is the (non-associative, non-unital)
     Lie ${\cal O}_X$-algebra
     ${\cal O}_X^{\nc,\Lie}:= ({\cal O}_X^{\nc}, [\,\cdot\,,\,\cdot\,])$
     with the commutator product
     $[s_1, s_2]:= s_1\cdot s_2 - s_2\cdot s_1$ for local sections
     of ${\cal O}_X^{\nc}$.
    A gauge theory on the D-brane $X$ corresponds to a choice
     of a gauge sheaf ${\cal G}_X$ embedded in ${\cal O}_X^{\nc,\Lie}$.
    Here, a {\it gauge sheaf} is a sheaf of ${\cal O}_X$-Lie-algebras
     that generalizes the notion of the Lie-algebra bundle associated
     to the adjoint representation of the gauge group of
     a principal bundle.\footnote{Readers
                         may wonder why we do not take
                         ${\cal O}_X^{\tinync,\tinyLie}$ or ${\cal G}_X$
                         directly to define the noncommutative structure
                         on $X$.
                        There are two reasons:
                         (1) The ``geometry" (in the sense
                             of ``points" and ``topology")
                             associated to a non-associative,
                             non-unital ring is less clear than that
                             for an associative unital ring
                             at the moment.
                         (2) Since the function ring of local charts
                             of the target space is associative and
                             unital, if we use
                             ${\cal O}_X^{\tinync,\tinyLie}$ for $X$,
                             we will have to consider ring-homomorphisms
                             from an associative unital ring
                             to a Lie ring.
                             The only such ring-homomorphism
                             is the zero-homomorphism.
                             This renders such setting containing
                             no contents as long as
                             ``probing a space(-time) via morphisms
                               into it" is concerned.
                         Cf.\ footnote~17.}
    This renders ${\cal E}_X$ a ${\cal G}_X$-module.
    Thus, it is enough to consider ${\cal O}_X^{\nc}$ and
     ${\cal E}_X$ as an ${\cal O}_X^{\nc}$-module.   
    %

   \medskip

   \item[(3.2)]
    A D-brane $\Phi: X\rightarrow Y$ on $Y$ determines a sheaf
     ${\cal O}_X \subset {\cal A}^{\nc} \subset {\cal O}_X^{\nc}$
     of subalgebras of ${\cal O}_X^{\nc}$,
     namely the image of the ring-system homomorphism
      ${\cal R}_Y\rightarrow {\cal R}_X$ that defines $\Phi$.
    The associated {\it gauge symmetry on the D-brane on} $Y$
     is given by the sheaf
     $\Centralizersheaf_{{\cal O}_X^{\tinync}}({\cal A}^{\nc})$
     of {\it centralizer subalgebras} of ${\cal A}^{\nc}$
     in ${\cal O}_X^{\nc}$.
    A continuous family $\Phi_t:X_t\rightarrow Y$ of deformations
     of the morphism $\Phi:X\rightarrow Y$ gives rise to
     a (not-necessarily flat) family
     $\Centralizersheaf_{{\cal O}_{X_t}^{\tinync}}({\cal A}^{\nc}_t)$
     of sheaves of algebras.
    This realizes the Higgsing/un-Higgsing behavior of
     the gauge symmetry on D-branes on $Y$ under deformations of
     D-branes on $Y$.\footnote{(1)
                                A priori,
                                one has a choice of whether or not
                                 the Higgsing/un-Higgsing of D-branes
                                 should be described as nearby points
                                 in the to-be-constructed moduli space
                                 of D-branes.
                                For a fixed string target-space $Y$,
                                 the Wilson theory-space of ``D-branes"
                                 in the region where they are still
                                 branes resembles the Wilson theory-space
                                 of a gauge system.
                                With the type of the gauge system fixed,
                                 we have a continuum for the latter
                                 theory-space.
                                The gauge group and hence the gauge bundle
                                 under Higgsing/un-Higgsing jump
                                 discontinuously but the situation is
                                 like that on the theory-space in
                                 Seiberg-Witten theory:
                                there is a continuum as the theory-space.
                                Another similar situation occurs in the
                                 geometric engineering of gauge theories,
                                 in which the compactification of
                                 a superstring theory
                                 on a degeneration family ${\cal X}$
                                 of Calabi-Yau $3$-spaces
                                 over a base $B$ gives rise to a family
                                 $\{\footnotesizeQFT_b\}_{b\in B}$
                                 of $d=4$ effective field theories,
                                 parameterized by $B$,
                                 whose gauge symmetry is enhanced at
                                 special locus of $B$ that corresponds
                                 to singular fibers of ${\cal X}/B$.
                                Mathematicians may also recall
                                 the moduli space ${\cal M}$ of coherent
                                 sheaves of a fixed Hilbert polynomial
                                 on a projective variety.
                                Even when ${\cal M}$ is connected,
                                 the function on ${\cal M}$ that assigns
                                  to an $[{\cal F}]\in {\cal M}$
                                  its sheaf-cohomology dimensions
                                  or Betti numbers is in general
                                  discontinuous.
                                The upper-semicontinuity of
                                 such a function, in particular $h^0$
                                 from the global section functor,
                                 on ${\cal M}$ can be taken as
                                 a resemblance of the phenomenon of
                                 enhancement of gauge symmetry due
                                 to additional zero/massless modes.

                                (2)
                                 It can happen that the ``good part" of
                                 the (coarse) moduli space of objects
                                 of different nature admit canonical
                                 identifications.
                                For example,
                                  the moduli space of maps,
                                  the moduli space of subschemes, and
                                  the moduli space of cycles
                                 canonically coincide when the maps are
                                 embeddings of reduced schemes
                                 with the trivial automorphism group.
                                Ignoring the issue of automorphisms,
                                it is the behavior under degenerations
                                 (i.e.\ moving away from such
                                  ``good part" of the moduli space)
                                 that the nature of the objects we intend
                                  to parameter reveals itself.
                                It is only when the degeneration feature
                                 distinct for each moduli problem is
                                 captured in the setting may one now hope
                                 to have a correct description of
                                 the objects and hence their moduli space.
                                Definition 2.2.3 is made
                                 with both (1) and (2) in mind.}
    \end{itemize}
 \end{itemize}

These highlights explain why we take Definition 2.2.3
 as a prototype intrinsic definition for D-branes
 (or D-brane world-volumes) in the region of
 the theory-space where ``branes are still branes".
Details of the case of D0-branes are given in Sec.~3 and Sec.~4.
The general higher-dimensional brane case can be thought of as
 sheafifying/smearing the discussion for D0-branes along
 a higher-dimensional cycle, chain, or more generally
 current in the sense of [G-H] or [Fe]; cf.\ [L-Y4].

\bigskip

\noindent
{\it Remark 2.2.5 $[$other intrinsic definitions$]$.} {\rm
 There have been other working mathematical intrinsic definitions for
  D-branes by other authors aiming also to understanding D-branes (in
  the region of Wilson's theory-space where ``branes are still branes").
 For example, there were
  the interpretation of D-branes as {\it stable torsion sheaves},
   given, e.g., in [H-S-T] in the algebro-geometric category
   from the viewpoint of BPS states and Gopakumar-Vafa invariants, and
  the notion of `{\it flat D-branes}',
   given in [B-M-R-S] in the smooth differential-geometric category
   from the viewpoint of K-theory.
 Each of these definitions singles out important key
  properties/features of D-branes in stringy literatures.
 In contrast, our prototype intrinsic definition of D-branes follows
  from the Grothendieck's viewpoint of Polchinski's work,
  phrased as the Polchinski-Grothendieck Ansatz for D-branes.
 This starting point is lower than these other existing
  intrinsic definitions and can reach up/be linked, for example,
  to [H-S-T] by considering D-brane images
   with the push-forward Chan-Paton sheaf on the target space  and
  to [B-M-R-S] by considering formal linear combinations of
   D-branes with Chan-Paton sheaves and their equivalence classes
   in the K-group of the D-brane.
} 

\bigskip

\section{$\Mor(\Space M_n({\LargeBbb C}), Y)$ as a coarse moduli space.}

We realize in this section the space
 $$
  \Mor(\Space M_n({\Bbb C}), Y)\;
  :=\; \Mor([{\cal R}], [\{M_n({\Bbb C})\}])\;
   =\; \Mor({\cal R}, [\{M_n({\Bbb C})\}])
 $$
 of morphisms from $\Space M_n({\Bbb C})$ to
  $Y = \Space[{\cal R}]
     = \Space( [\{R_{\gamma}\}_{\gamma\in C} \doublearrow
                \{R_{\gamma_1\gamma_2}\}_{\gamma_1,\gamma_2\in C}] )$
 as a constructible set in a topological space from an adhesion of
 affine varieties/${\Bbb C}$.

\bigskip

\subsection{Central localizations of Artinian rings and their modules.}

Recall first the {\it Structure Theorem of Artinian Rings}:

\bigskip

\noindent
{\bf Theorem 3.1.1 [Artinian ring].} ([A-M], [A-N-T], and [Jat].) {\it
 \begin{itemize}
  \item[$(1)$]
   Let $R$ be an Artinian ring.
   The center $Z(R)$ of $R$ is a commutative Artinian ring
    and hence has finitely-many maximal ideals.
   Let $t$ be the number of maximal ideals in $Z(R)$.

  \item[$(2)$]
   There exist a unique collection $\{e_1,\,\cdots\,,\,e_t\}$
    of orthogonal primitive idempotents in $Z(R)$ such that
    $1=e_1+\,\cdots\,+e_t$ and that
   $R$ is the direct sum of the two-sided ideals
    $R=Re_1+\,\cdots\,+Re_t$.
   Up to permutations, the collection
    $\{Re_1\,,\, \cdots\,,\,Re_t\}$
    is unique with respect to the following property:
    \begin{itemize}
     \item[$\cdot$]
      $R=I_1+\,\cdots\,+I_{t^{\prime}}$, where
       $I_i$ are two-sided ideals of $R$,
       $I_i\cdot I_j=0$ for $i\ne j$, and
       each $I_i$ is indecomposable in the sense that
        $I_i$ cannot be decomposed as a direct sum
         $I_{i^{\prime}}+I_{i^{\prime\prime}}$
         with $I_{i^{\prime}}$ and $I_{i^{\prime\prime}}$
         non-zero two-sided ideals.
    \end{itemize}
   Under such decomposition of $R$,
    each $R_i:=Re_i$ is itself an Artinian ring with identity $e_i$
     and
    the decomposition $R=Re_1+\,\cdots\,+Re_t$ can be written
    as the product of rings $R=R_1\times\,\cdots\,\times R_t$.
   This decomposition restricts to a decomposition
    $Z(R)=Z(R_1)\times\,\,\cdots\,\times Z(R_t)$
    with each $Z(R_i)$, $i=1,\,\ldots\,, t$, an Artinian local ring.

  \item[$(3)$]
   Let $J(R)$ be the Jacobson radical of $R$.
   Then there is an orthogonal idempotent decomposition
     $$
       1\;=\; \sum_{j_1=1}^{l_1}e_{1j_1}\,
             +\,\cdots\,+\sum_{j_s=1}^{l_t}e_{tj_t}
     $$
     in $R$ that refines the decomposition
     $1=e_1+\,\cdots\,+e_t$ in $Z(R)$,
     with $e_i=\sum_{j_i=1}^{l_i}e_{ij_i}$,
    such that the image $\bar{e}_{ij_i}$ of $e_{ij_i}$
      in $R/J(R)$ lies in $Z(R/J(R))$ and
     that
      $$
       \bar{1}\;=\;  \sum_{j_1=1}^{l_1}\bar{e}_{1j_1}\,
                     +\,\cdots\, +\, \sum_{j_t=1}^{l_t}\bar{e}_{tj_t}
      $$
      is an orthogonal primitive idempotent decomposition
      in $Z(R/J(R))$.
   Let
    $$
     {\frak m}_{ij_i}\; :=\;  R\, (1-e_{ij_i})\,R\,,
     \hspace{2em}\mbox{for $1\le i\le t$ and $1\le j_i\le l_i$}\,,
    $$
    and $\Spec R$ be the set of all prime ideals in $R$.
   Then all prime ideals in $R$ are maximal ideals and
    $$
     \Spec R\; =\;
      \{\, {\frak m}_{ij_i}\;:\; 1\le i\le t\,,\, 1\le j_i\le l_i\, \}\,.
    $$

  \item[$(4)$]
   Consider the directed graph $\Gamma_R$
    with the set of vertices $\Spec R$ and
     a directed edge
      ${\frak m}_{i_1j_{i_1}}\rightarrow {\frak m}_{i_2j_{i_2}}$
      for each pair $(e_{i_1j_{i_1}}, e_{i_2j_{i_2}})$
          with $e_{i_1j_{i_1}} J(R) e_{i_2j_{i_2}}\ne 0$.
   Then $\Gamma_R$ has exactly $t$-many connected components
    $\Gamma_R^{(i)}$, $i=1,\,\cdots\,,t$,
    with the set of vertices of $\Gamma_R^{(i)}$ being
     $\{{\frak m}_{ij_i}\,:\, 1\le j_i\le l_i\}\,$.
   The two graphs $\Gamma_R^{(i)}$ and $\Gamma_{R_i}$
    are canonically isomorphic.
   In particular, each $\Gamma_{R_i}$ is connected.

  \item[$(5)$]
   By definition,
    $J(R)=\cap_{i=1}^t\cap_{j=1}^{l_i}\, {\frak m}_{ij_i}$.
   The quotient ${\frak m}_{ij_i}/J(R)$,
    with the induced addition and multiplication from those of $R$,
    is a simple ring and hence is isomorphic to a matrix ring
    $M_{n_{ij_i}}(k_{ij_i})$ for some skew-field $k_{ij_i}$.
   The decomposition $R=Re_1+\,\cdots\,+Re_t$ restricts to
    a decomposition $J(R)=J(R)_1+\,\cdots\,+J(R)_t$, which can be
     written canonically as $J(R)=J(R_1)\times\,\cdots\,\times J(R_t)$.
   With respect to this, one has isomorphisms
    $$
     R/J(R)\;\simeq\; \prod_{i=1}^t\,R_i/J(R_i)\;
      \simeq\;
       \prod_{i=1}^t\prod_{j_i=1}^{l_i}\, M_{n_{ij_i}}(k_{ij_i})\,.
    $$
 \end{itemize}
} 

\bigskip

\noindent
{\it Remark 3.1.2 $[$quiver$]$.} {\rm
 The graph $\Gamma_R$ associated to an Artinian ring $R$
  (as an $R$-module) in Theorem is an example of
  (various) quivers associated to an $R$-module.
 See Sec.~4.1 and footnote 36 for a theme in which
  we bring this in again.
} 

\bigskip

The theorem gives a visualization of an Artinian algebra $R/{\Bbb C}$
 (e.g.\ $M_n({\Bbb C})$ and its subalgebras)
 as a noncommutative space of the form:
 \begin{quote}
  \hspace{-1ex}``a
   finite collection of commutative points (i.e.\ $\Spec Z(R)$),
    with each point dominated/shadowed by a noncommutative cloud
     (i.e.\ $Z(R_i)\subset R_i$, where $R_i:=Re_i$);
   associated to each noncommutative cloud (i.e.\ $R_i$) over
    a commutative point (i.e.\ $\Spec Z(R_i)$) are
    a refined collection of commutative points
    (i.e.\ $\Spec(\sum_{j=1}^{l_i}{\Bbb C}\cdot e_{ij})$)
    split off from and stacked over that point
     (more precisely,
      $\Spec Z(R_i)_{\redscriptsize}$) and
    are dominated/shadowed by that cloud
    (i.e.\  $\sum_{j=1}^{l_i}e_{ij}=e_i$ and
            ${\Bbb C}\cdot e_i
             \subset \sum_{j=1}^{l_i}{\Bbb C}\cdot e_{ij}
             \subset R_i$) and
    bound by directed bonds
    (i.e.\ $e_{ij_1} J(R_i) e_{ij_2}$
      with the direction from $e_{ij_1}$ to $e_{ij_2}$)
    created through that cloud (i.e.\ $R_i$)".
 \end{quote}

\bigskip

The following are immediate consequences of the theorem.

\bigskip

\noindent
{\bf Lemma 3.1.3 [central non-zero-divisor invertible].} {\it
 Let $R$ be an Artinian ring and $r\in Z(R)$ be
 a non-zero-divisor in $R$. Then $r$ is invertible in $R$.
} 

\bigskip

\noindent
{\bf Lemma 3.1.4 [direct-sum decomposition of module].}
(Cf.\ {\sl Peirce decomposition}.)
{\it
  Let
   $R$ be an Artinian ring and
   $R=Re_1+\,\cdots\,+Re_t =: R_1+\,\cdots\,+R_t$
    be a decomposition of $R$ as in Theorem 3.1.1 (2).
  Let $M$ be an $R$-module.
  Then,
   $M=e_1M+\,\cdots\,+e_tM =: M_1+\,\cdots\,+M_t$
     is a direct-sum decomposition of $M$
    such that $R_iM_i=M_i$ and $R_jM_i=0$ for $j\ne i$.
  In particular, $M_i$ is a $R_i$-module for $i=1,\,\ldots\,, t$.
} 

\bigskip

\noindent
{\bf Corollary 3.1.5 [localization = quotient].} {\it
 {\rm (With notations from above.)}
 $R_i$ is canonically isomorphic to
  both the quotient $R/(e_j:j\ne i)=R/(\sum_{j\ne i}e_j)$ of $R$
  and  the localization $R[S_i^{-1}]$ of $R$,  where $S_i$ is
             the multiplicatively closed subset $\{1, e_i\}$.
 Similarly,
 $M_i$ is canonically isomorphic to
  both the quotient $M/(\sum_{j\ne i} M_j)$ of $M$
  and the localization $M[S_i^{-1}]$ of $M$.
} 

\bigskip

\noindent
{\bf Corollary 3.1.6 [localization: standard form].} {\it
 $(1)$
 Any nonzero central localization $R\rightarrow R^{\prime}$
  of an Artinian ring $R$ is realized by inverting
  a finite multiplicatively closed subset $S \subset Z(R)$
  that consists only of idempotents.
 I.e.\ $R^{\prime}=R[S^{-1}]$ and
  $R^{\prime}\rightarrow R$ is $R\rightarrow R[S^{-1}]$
  for an afore-mentioned $S$.
 $(2)$
 Any central localization $f:R\rightarrow R^{\prime}$
  of an Artinian ring $R$ is a quotient of $R$ that admits
  a ring-set-homomorphism\footnote{See Definition 3.2.2.}
  $g:R^{\prime}\rightarrow R$ such that $f\circ g=\Id_{R^{\prime}}$.
 $(3)$
 Fix a direct-sum decomposition $R=R_1+\,\cdots\,+R_t$
  from Theorem 3.1.1 (2).   
 Then the localization $f: R\rightarrow R^{\prime}$ in $(2)$
   is simply the projection of $R$ onto the sum
   $R_{i_1}+\,\cdots\,+R_{i_{t^{\prime}}}$ of some direct summands and
  $g:R^{\prime}\rightarrow R$ in $(2)$ can be taken to be the inclusion
   of $R_{i_1}+\,\cdots\,+R_{i_{t^{\prime}}}$ into $R$.
} 

\bigskip

\noindent
{\it Proof.}
 Let $R=R_1+\,\cdots\,+R_t$ be a direct-sum decomposition
  of $R$ from Theorem 3.1.1 (2).   
 Then $S=S_1+\,\cdots\,+ S_t$, where $S_i:=e_iS \subset Z(R_i)$,
  is a direct-sum decomposition of $S$ and
 $R[S^{-1}]=R_1[S_1^{-1}]\times\,\cdots\,\times R_t[S_t^{-1}]$
  canonically.
 This reduces the proof to
  the case that $t=1$ in the decomposition of $R$
  (i.e.\ the case $Z(R)$ is an Artinian local ring).

 When $Z(R)$ is an Artinian local ring,
  $R[S^{-1}]=0$ if $S$ contains an element in the maximal ideal
   of $Z(R)$, as such an element is nilpotent.
 Otherwise, all elements of $S$ are not in the maximal ideal of $Z(R)$;
 then they are all invertible and, hence, $R[S^{-1}]=R$.
 In the former (resp.\ latter) case, we may replace $S$ by
  $\{1,0\}$ (resp.\ $\{1\}$).
 The corollary now follows.

\noindent\hspace{15cm}$\Box$

\bigskip

\noindent
{\bf Lemma 3.1.7 [localization in terms of generators of $S$].} {\it
 Let $R$ be an Artinian ring and $S$ be a multiplicatively closed
  subset in $Z(R)$, generated by\footnote{I.e.\
                                         an element of $S$ is either
                                         the identity $1$ or a monomial
                                         of $s_1,\,\cdots\,, s_l$.}
  $\{s_1,\,\cdots\,, s_l\}$.
 Let $n_0$ be a positive integer such that
  every nilpotent element $r$ of $R$ satisfies $r^{n_0}=0$.
 Then $R[S^{-1}]=R/\sum_{i=1}^l(s_i^{n_0})^{\perp}$,
  where $(\,\bullet\,)^{\perp}:=\{r\in R: (\,\bullet\,)\cdot r =0\}$.
} 

\bigskip

\noindent
{\it Proof.}
 This follows immediately from Corollary 3.1.6.

\noindent\hspace{15cm}$\Box$

\bigskip

\subsection{$\Mor(\Space M_n({\largeBbb C}),Y)$
            as a coarse moduli space.}

\noindent
{\bf Definition 3.2.1 [ring-subset].} {\rm
 Let $R=(R, 0,1,+,\,\cdot\,)$ be a ring, with the identity $1$.
 An additive subgroup $R^{\prime}\subset R$
   is called a {\it ring-subset} of $R$
  if, in addition,
  (1) $R^{\prime}$ is closed under the multiplication $\cdot$ in $R$,
   and
  (2) there is an element $e\in R^{\prime}$ such that
      $(R^{\prime},0,e,+,\,\cdot\,)$ is a ring with the identity $e$.
} 

\bigskip

\noindent
{\bf Definition 3.2.2 [ring-set-homomorphism].} {\rm
 Let $R$ and $S$ be rings with the identity $1_R$ and $1_S$ respectively.
 A map $\varphi: R\rightarrow S$ is called a {\it ring-set-homomorphism}
  if $\varphi$ satisfies all the requirement for a ring-homomorphism
   {\it except} that it is not required that $\varphi(1_R)=1_S$.
} 

\bigskip

Note that $e$ in Definition 3.2.1
 is unique and satisfies $e^2=e$.

\bigskip

\noindent
{\bf Example 3.2.3 [ring-subset].}
 The image $\varphi(R)$ in Definition 3.2.2
  is a ring-subset of $S$ with the identity $\varphi(1_R)$.
 In particular, $\{0\}\subset R$ is the minimal ring-subset of $R$.

\bigskip

We will retain these terminologies for algebras and
 algebra-homomorphisms over a fixed ground field as well.

\bigskip

\begin{flushleft}
{\bf Surrogates of the Azumaya-type noncommutative point
     $\Space M_n({\Bbb C})$.}
\end{flushleft}
$M_n({\Bbb C})$ is a simple ring in the sense that
 it is semi-simple as a left $M_n({\Bbb C})$-module and
  has the only two-sided ideals the zero-ideal $({\mathbf 0})$ and
   itself $M_n({\Bbb C})$.
In particular,
 the only prime ideal of $M_n({\Bbb C})$ is $({\mathbf 0})$  and
 the center $Z(M_n({\Bbb C}))$ of $M_n({\Bbb C})$ is given
  by ${\Bbb C}\cdot{\mathbf 1}$.
There are only two Gabriel filters on $M_n({\Bbb C})$:
 ${\frak F}_0$ that is generated by $({\mathbf 0})$ and is given
  by the set of all left ideals of $M_n({\Bbb C})$ and
 ${\frak F}_1:= \{M_n({\Bbb C})\}$.
The localization of $M_n({\Bbb C})$ with respect to ${\frak F}_0$
  (resp.\ ${\frak F}_1$) is the zero-ring $0$
  (resp.\ $M_n({\Bbb C})$ itself).
The former (resp.\ latter) covers the notion of the localization
 of $M_n({\Bbb C})$ with respect to a non-invertible
 (resp.\ invertible) element.
Thus, directly on $M_n({\Bbb C})$,
 we see only a seemingly barren geometry.
Things change when we bring in the notion of surrogates introduced
 in Sec.~1.1.

A surrogate of the Azumaya-type noncommutative point
 $\Space M_n({\Bbb C})=(\Spec{\Bbb C},{\Bbb C}, M_n({\Bbb C}))$
 is given ring-theoretically by a subalgebra pair
 ${\Bbb C}\subset C \subset R \subset M_n({\Bbb C})$
 with $C\subset Z(R)$.
It follows from Corollary 3.1.6 that
a finite central cover of the sub-${\Bbb C}$-algebra $R$
 of $M_n({\Bbb C})$ can be described by a finite collection
 $\{(R_{\alpha}, e_{\alpha})\}_{\alpha\in A}$
 of ring-subsets of $R$ (and hence of $M_n({\Bbb C})$)
 that satisfies the following conditions:
 \begin{itemize}
  \item[$(0)$]
   $R=\sum_{\alpha\in A} R_{\alpha}$.

  \item[$(1)$]
   $e_{\alpha_1}$ commutes with elements of $R_{\alpha_2}$
    for all $\alpha_1$, $\alpha_2\in A$.

  \item[$(2)$]
   $e_{\alpha_1}R_{\alpha_2}=e_{\alpha_2}R_{\alpha_1}$
    for all $\alpha_1$, $\alpha_2\in A$.

  \item[$(3)$]
   $e_{\alpha_1}e_{\alpha_2}\in R_{\alpha_1}$
   for all $\alpha_1$, $\alpha_2\in A$.

  \item[$(4)$]
   Fix a well-ordering of the index set $A$; then
   \begin{eqnarray*}
    \lefteqn{
     1\;=\; \sum_{\alpha} e_{\alpha}\;
           -\, \sum_{\alpha_1<\alpha_2} e_{\alpha_1}e_{\alpha_2}\;
           +\, \sum_{\alpha_1<\alpha_2<\alpha_3}
                  e_{\alpha_1}e_{\alpha_2}e_{\alpha_3} }\\[.6ex]
     && \hspace{6em}
         \pm\; \cdots\;
           +\, (-1)^{|A|+1}\,
                \sum_{\alpha_1<\,\cdots\,<\alpha_{|A|}}
                    e_{\alpha_1}\,\cdots\,e_{\alpha_{|A|}}\,.
   \end{eqnarray*}
 \end{itemize}
Conditions (1), (2), and (3) imply that
 $e_{\alpha_2}R_{\alpha_1}
  =R_{\alpha_1}\cap R_{\alpha_2} = e_{\alpha_1}R_{\alpha_2}$,
 which is itself a ring with the identity $e_{\alpha_1}e_{\alpha_2}$.
In particular,
 $(e_{\alpha_2}R_{\alpha_1}, e_{\alpha_2}e_{\alpha_1})
   =((e_{\alpha_1}e_{\alpha_2})R_{\alpha_1},
      e_{\alpha_1}e_{\alpha_2})$
 is a ring-subset of both rings $(R_{\alpha_1},e_{\alpha_1})$
 and $(R_{\alpha_2},e_{\alpha_2})$.
Condition (4) simplifies to $1=\sum_{\alpha}e_{\alpha}$
 when $R=\sum_{\alpha}R_{\alpha}$ is a direct sum.

Conversely, one has the following proposition:

\bigskip

\noindent
{\bf Proposition 3.2.4
     [subring in terms of a collection of ring-subsets].}
{\it
 $(1)$
 Let $\{(R_{\alpha},e_{\alpha})_{\alpha\in A}\}$
  be a finite collection of ring-subsets of $M_n({\Bbb C})$
  that satisfies Conditions $(1)$, $(2)$, $(3)$, and $(4)$ above.
 Then $R:=\sum_{\alpha\in A} R_{\alpha}$
  contains the identity ${\mathbf 1}$ of $M_n({\Bbb C})$ and
  is a sub-${\Bbb C}$-algebra of $M_n({\Bbb C})$.
 $(2)$
 There are tautological ring-homomorphisms
  $R\rightarrow R_{\alpha}$, $\alpha\in A$, that render
  the collection $\{R\rightarrow R_{\alpha}\}_{\alpha\in A}$
  a finite central cover of $R$.
} 

\bigskip

\noindent
{\it Proof.}
 Observe that for a ring-subset $(P,e_P)$ and
   an idempotent $e^{\prime}$ of a ring $Q$ that
   commutes with the elements in $P$,
  $(e^{\prime}P, e^{\prime}e_P)$ is another ring-subset of $Q$.
 In particular, elements in $e^{\prime}P$ are closed under
  the multiplication in $Q$.
 Moreover, if, in addition, $e^{\prime}e_p\in P$,
 then $P=(e_P-e^{\prime}e_P)P+(e^{\prime}e_P)P$
  is an orthogonal direct-sum decomposition for $P$
  (when neither summand is zero).
 Using these observations, one can show that
  Properties (1), (2), and (3) imply that
  $R_{\alpha_1}R_{\alpha_2}\subset R_{\alpha_1}+R_{\alpha_2}$
  for all $\alpha_1$, $\alpha_2\in A$.
 This proves that $R:=\sum_{\alpha}R_{\alpha}$ is closed under
  the multiplication in $M_n({\Bbb C})$ as
  $R\cdot R
    \subset \sum _{\alpha_1,\alpha_2\in A} (R_{\alpha_1}+R_{\alpha_2})
    = R$.

 Now let
   \begin{eqnarray*}
    \lefteqn{
     e\; :=\; \sum_{\alpha} e_{\alpha}\;
           -\, \sum_{\alpha_1<\alpha_2} e_{\alpha_1}e_{\alpha_2}\;
           +\, \sum_{\alpha_1<\alpha_2<\alpha_3}
                  e_{\alpha_1}e_{\alpha_2}e_{\alpha_3} }\\[.6ex]
     && \hspace{6em}
         \pm\; \cdots\;
           +\, (-1)^{|A|+1}\,
                \sum_{\alpha_1<\,\cdots\,<\alpha_{|A|}}
                    e_{\alpha_1}\,\cdots\,e_{\alpha_{|A|}}\,.
   \end{eqnarray*}
 Then it follows from the above that $e\in R$.
 For $r\in R_{\alpha_i}$, $\alpha_i\in A$, one can check directly
  that $re=r$, using the property that $re_{\alpha_i}=r$ and
   the above defining expression of $e$.
 This implies that $er=r$ for every $r\in R$.
 It follows that
  $(R, e)$ is a ring-subset of $M_n({\Bbb C})$.
 The additional Condition (4), $e=1$, implies then
  that $R$ is a subalgebra of $M_n({\Bbb C})$.
 This proves Statement (1).

 Condition (1) implies that $\{e_{\alpha}\}_{\alpha\in A}\subset Z(R)$.
 For each $\alpha\in A$, the commutativity, idempotent property, and
  that both $e_{\alpha}R_{\alpha}=R_{\alpha}$
       and  $e_{\alpha}(e_{\alpha}R)=e_{\alpha}R$ hold
   imply that
  the orthogonal direct-sum decomposition
   $R=e_{\alpha}R+(1-e_{\alpha})R$ of $R$ coincides with
   the decomposition $R=R_{\alpha}+e_{\alpha}^{\perp}$,
   where $e_{\alpha}^{\perp}:= \{r\in R: e_{\alpha}r=0\}$.
 This shows that the projection map $R\rightarrow R_{\alpha}$
  from the above decomposition is identical with
  the central localization of $R$ with respect to
  the multiplicatively closed subset $\{1,e_{\alpha}\}$.
 Furthermore, $\sum_{\alpha\in A}e_{\alpha}$ is invertible in $Z(R)$.
 Thus, $\{(R\rightarrow R_{\alpha})\}_{\alpha\in A}$
  is a central finite cover of $R$.
 This proves Statement (2).

\noindent\hspace{15cm}$\Box$

\bigskip

\begin{flushleft}
{\bf The space $\Mor^{\ringsetscriptsize}(R,M_n(\Bbb C))$
         of ring-set-homomorphisms from $R$ to $M_n({\Bbb C})$.}
\end{flushleft}
Let
 $R$ be a finitely-presentable algebra over ${\Bbb C}$ and
 $\Mor^{\ringsetscriptsize}(R,M_n(\Bbb C))$
  be the set of ring-set-homomorphisms from $R$ to $M_n({\Bbb C})$.
We will construct a topology on
 $\Mor^{\ringsetscriptsize}(R,M_n(\Bbb C))$ in this theme.

Let
\begin{itemize}
 \item[$\cdot$]
  $R\, =\,
    \langle g_0, g_1,\,\cdots\,, g_l\rangle /(r_1,\,\cdots\,, r_m)\;$
   be a presentation of $R$ as a quotient
    of the free unital associative ${\Bbb C}$-algebra
     $\langle g_0, g_1,\,\cdots\,, g_l\rangle$
     generated by $g_0, g_1,\,\cdots\,, g_l$
    by the two-sided ideal $(r_1,\,\cdots\,,r_m)$ generated by
     $\{r_i=r_i(g_0,\,\cdots\,,g_l): i=1,\,\ldots\,,m\}$.
   Here, for later use, we have
    the redundant generator $g_0 =$ the identity $1$ and
    the redundant relators $g_0g_i=g_ig_0=g_i$, $i=0,1,\,\ldots\,,l$,
     contained in the relator set $\{r_1,\,\cdots\,, r_m\}$.

 \item[$\cdot$]
  $\Gr^{(2)}(n; d,n-d)\,\simeq\,
   \GL_n({\Bbb C})/(\GL_d({\Bbb C})\times\GL_{n-d}({\Bbb C}))\;$
  be the Grassmannian manifold of ordered pairs $(\Pi_1, \Pi_2)$
  of ${\Bbb C}$-linear subspaces of ${\Bbb C}^n$ with
  $\dimm \Pi_1=d$, $\dimm \Pi_2=n-d$, and $\Pi_1+\Pi_2={\Bbb C}^n$;

 \item[$\cdot$]
  ${\mathbf 1}_d$, $d=0,\,\ldots\,,n$,
   be the diagonal matrix
   $\Diag(1,\,\cdots\,,1,0,\,\cdots\,,0)$ in $M_n({\Bbb C})$
    whose first $d$ diagonal entries are $1$ and the rest $0$,
   (here, ${\mathbf 1}_0=$ the zero-matrix $0$ and
          ${\mathbf 1}_n={\mathbf 1}$ by convention);  and

 \item[$\cdot$]
  `$m_1\sim m_2$'$\;$
   means that $m_1$ and $m_2$ are in the same adjoint
   $\GL_n({\Bbb C})$-orbit in $M_n({\Bbb C})$.
\end{itemize}
Let $\Rep^{\ringsetscriptsize}(R,M_n(\Bbb C))$
 be the subvariety of the affine space
 ${\Bbb A}^{n^2}_{(0)}\times{\Bbb A}^{n^2}_{(1)}\times\,\cdots\,
   \times {\Bbb A}^{n^2}_{(l)}$
 (here ${\Bbb A}^{n^2}_{(i)}$ has the polynomial coordinate ring
  ${\Bbb C}[m_{i, jk}: 1\le j, k\le n]$, $i=0,\,\ldots\,,l$)
 determined\footnote{I.e.\ taking the reduced scheme associated to
                                    the possibly nonreduced subscheme
                                    described by the ideal generated by
                                    these equations.}
 by the system of equations
 $$
  r_1(M_0, M_1,\,\cdots\,, M_l)\;=\; \cdots\;
  =\; r_m(M_0, M_1,\,\cdots\,, M_l)\;
  =\; \mbox{the zero-matrix}\; 0\,\in\, M_n({\Bbb C})\,,
 $$
 where $M_i=(m_{i,jk})_{jk}$.
Note that this is like the ordinary representation variety of
 $R$ in $M_n({\Bbb C})$ {\it except} that it is not required
 that $M_0=$ the identity ${\mathbf 1}\in M_n({\Bbb C})$.
For convenience, we will call the reduced affine scheme
 $\Rep^{\ringsetscriptsize}(R,M_n(\Bbb C))$
 the {\it representation variety} in our discussion.
By construction, we have
 the Zariski topology on $\Rep^{\ringsetscriptsize}(R,M_n(\Bbb C))$
  and
 the analytic topology on the set
   $\Rep^{\ringsetscriptsize}(R,M_n(\Bbb C))_{\scriptsizeBbb C}$
   of ${\Bbb C}$-points of
    $\Rep^{\ringsetscriptsize}(R,M_n(\Bbb C))$.
Regard ${\Bbb C}^n$ as the unique non-zero irreducible
 $M_n({\Bbb C})$-module.
Then, the correspondence
 $e\mapsto (e\cdot{\Bbb C}^n, e^{\perp})$,
  where $e^{\perp}$ here $=\{v\in {\Bbb C}^n: e\cdot v=0\}$,
 gives rise to a (continuous) map
  from the set of idempotents $\sim {\mathbf 1}_d$ in $M_n({\Bbb C})$
  to $\Gr^{(2)}(n; d, n-d)$.
It follows that
the projection map
 $\pi_{(0)}:
   {\Bbb A}^{n^2}_{(0)}\times{\Bbb A}^{n^2}_{(1)}\times\,\cdots\,
     \times {\Bbb A}^{n^2}_{(l)}  \rightarrow  {\Bbb A}^{n^2}_{(0)}$
 restricts to a map
 $$
  \pi_{(0)}\; :\;
   \Rep^{\ringsetscriptsize}(R,M_n(\Bbb C))\;
    \longrightarrow\; \amalg_{d=0}^n \Gr^{(2)}(n; d, n-d)\,.
 $$
Let
 $$
  \Rep^{\ringsetscriptsize}(R,M_n(\Bbb C))_{(d)}\;
   :=\; \pi_{(0)}^{-1}(\Gr^{(2)}(n; d, n-d))\,.
 $$

As a set,
 $\Mor^{\ringsetscriptsize}(R,M_n({\Bbb C}))
  =\Rep^{\ringsetscriptsize}(R,M_n({\Bbb C}))_{\scriptsizeBbb C}$.
This identification defines
 a preliminary analytic topology ${\cal T}_0$ on
 $\Mor^{\ringsetscriptsize}(R,M_n(\Bbb C))
   = \amalg_{d=0}^n \Mor^{\ringsetscriptsize}(R,M_n(\Bbb C))_{(d)}$
 by bringing over the analytic topology on
 $\Rep^{\ringsetscriptsize}(R,M_n({\Bbb C}))_{\scriptsizeBbb C}$.
We then modify this preliminary analytic topology,
  following an analytic format of a {\it valuative criterion},
 so that
  each $\Mor^{\ringsetscriptsize}(R,M_n(\Bbb C))_{(d^{\prime})}$,
  $d^{\prime}<d$, adheres to
  $\Mor^{\ringsetscriptsize}(R,M_n(\Bbb C))_{(d)}$ appropriately,
  for $d=1,\,\cdots\,, n$, in the new topology.
Let
 $T$ be a (commutative, Noetherian) integral domain over ${\Bbb C}$
  and
 $(D:=(\Spec T)_{\scriptsizeBbb C}\,,\,p)$
  be the associated analytic space
  together with a base ${\Bbb C}$-point $p$.
Note that the residue field $\kappa_p$ of $T$ at $p$ is canonically
 isomorphic to ${\Bbb C}$.
Let $T_p$ be the localization of $T$ at $p$ and
 $Q_T$ be the field of fractions of $T$.
Then, $T\subset T_p\subset Q_T$,
 $T_p$ is a valuation ring of $Q_T$
  (regarded now as the field of fractions of $T_p$), and
 $M_n(T)\subset M_n(T_p)\subset M_n(Q_T)$.

\bigskip

\noindent
{\bf Definition 3.2.5 [limit of family ring-set-homomorphisms].} {\rm
 Let $\phi:R\rightarrow M_n(Q_T)$ be a ring-set-homomorphism
  such that there exists a unique idempotent $e\in M_n(T_p)$
  such that
   (1) $e\in Z(\image\phi)$;
   (2) $e\cdot\phi$ is a ring-set-homomorphism from $R$ to $T_p$;
       in particular,
       $(e\cdot\phi)|_p: R\rightarrow M_n(\kappa_p)=M_n({\Bbb C})$
       makes sense;
  (3) $\image(e\cdot\phi)|_p$ is the unique maximum
      (with respect to inclusion) in the set of ring-subsets
      $\image(e^{\prime}\cdot\phi)$ of $M_n(\kappa_p)$,
      where $e^{\prime}$ satisfies Condition $(1)$ and Condition $(2)$
      above.
 For such a $\phi$, we call $(e\cdot\phi)|_p$ the {\it limit}
  of $\phi$ over $D$ at $p$.
} 

\bigskip

\noindent
Such a $\phi$ defines a rational map
 $\Phi_{\phi}: (D,p) \dottedrightarrow
               \Mor^{\ringsetscriptsize}(R,M_n(\Bbb C))^{{\cal T}_0}$
 that is assigned the value $(e\cdot\phi)|_p$ at $p$.

\bigskip

\noindent
{\bf Definition 3.2.6 [$\Mor^{\ringsetscriptsize}(R,M_n(\Bbb C))$
                      with analytic topology].} {\rm
 With the notations from above,
 let ${\cal T}$ be the weakest topology on
   $\Mor^{\ringsetscriptsize}(R,M_n(\Bbb C))_{\scriptsizeBbb C}$
  such that
   \begin{itemize}
    \item[(1)]
     the tautological inclusion
     $\Mor^{\ringsetscriptsize}(R,M_n(\Bbb C))^{{\cal T}_0}_{(d)}
      \hookrightarrow
      \Mor^{\ringsetscriptsize}(R,M_n(\Bbb C))^{\cal T}$
     of sets is an embedding of topological spaces,
     for $d=0,\,\cdots\,,n$, and that

    \item[(2)]
     $\Phi_{\phi}$ is continuous at $p$ for all $T$, $(D,p)$,
     and $\phi$ in Definition 3.2.5.\footnote{This
                      is a valuative criterion.
                     The meaning of this topology in terms of analytic
                      geometry is as follows.
                     Under deformations of a morphism from
                      $\footnotesizeSpace M_n({\footnotesizeBbb C})$
                      to $\Space R$, some connected components of
                      the image points of
                      $\footnotesizeSpace M_n({\footnotesizeBbb C})$
                      may move away toward the boundary at infinity
                      of $\footnotesizeSpace R$ and disappear in the end.
                     This corresponds to a drop from
                      $M_0\sim {\mathbf 1}_d$ to some
                      $M_0\sim {\mathbf 1}_{d^{\prime}}$
                      with $d^{\prime}<d$.
                     When we consider only
                      $\Mor(\Space M_n({\footnotesizeBbb C}), \Space R)$
                      by itself, $M_0\sim {\mathbf 1}$ must always hold.
                     However, when we consider
                      $\Mor(\Space M_n({\footnotesizeBbb C}), \Space R)$
                       that occurs as a subset in
                       $\Mor(\Space M_n({\footnotesizeBbb C}),
                                                        \Space{\cal R})$
                       for a gluing system ${\cal R}$ of rings
                       that contains $R$ as a member,
                     it can happen that
                      some of the connected components of the image
                      of a morphism
                      $\Space M_n({\footnotesizeBbb C})
                                         \rightarrow \Space{\cal R}$
                      is not contained in $\Space R$.
                     This explains geometrically why,
                       in the equivalent ring-theoretic language,
                      we enlarge here the class of maps from
                       ring-homomorphisms to ring-set-homomorphisms.
                     Furthermore, when the morphism deforms,
                      the number of connected components in $\Space R$
                      of the image of morphisms from
                      $\Space M_n({\footnotesizeBbb C})$ to
                      $\Space{\cal R}$ can change.
                     The topology on
                       $\Mor^{\ringsettiny}(R,M_n({\footnotesizeBbb C}))$
                       defined in Definition 3.2.6
                       ring-theoretically
                      takes all these issues into account.
                     Such treatment automatically comes up and
                      is rerquired in building a (general) morphism
                      from $[{\cal R}]$ to $[\{M_n({\Bbb C})\}]$,
                      following Definition 1.2.14.}
                     %
   \end{itemize}
 ${\cal T}$ is called the {\it analytic topology} on
  $\Mor^{\ringsetscriptsize}(R,M_n(\Bbb C))$.
} 

\bigskip

\noindent
{\bf Proposition 3.2.7 [independence of presentation].} {\it
 $(\Mor^{\ringsetscriptsize}(R,M_n(\Bbb C)),{\cal T})$ is independent
 of the choice of presentations of $R$ in the construction.
} 

\bigskip

\noindent
{\it Proof.}
 Associated to a new presentation
  $$
   R\; =\;
    \langle g^{\prime}_0, g^{\prime}_1,\,\cdots\,,
            g^{\prime}_{l^{\prime}}\rangle /
     (r^{\prime}_1,\,\cdots\,, r^{\prime}_{m^{\prime}})
  $$
  of $R$ is a canonical ring-isomorphism
  $$
    f^{\sharp}\;:\;
     \langle g^{\prime}_0, g^{\prime}_1,\,\cdots\,,
             g^{\prime}_{l^{\prime}}\rangle /
     (r^{\prime}_1,\,\cdots\,, r^{\prime}_{m^{\prime}})\;
   \stackrel{\sim}{\longrightarrow}\;
    \langle g_0, g_1,\,\cdots\,, g_l\rangle /(r_1,\,\cdots\,, r_m)\,,
  $$
  represented by a noncanonical ring-homomorphism
  $\tilde{f}^{\sharp}:
    \langle g^{\prime}_0, g^{\prime}_1,\,\cdots\,,
           g^{\prime}_{l^{\prime}}\rangle
    \rightarrow
    \langle g_0, g_1,\,\cdots\,, g_l\rangle$.
 $\tilde{f}^{\sharp}$ induces contravariantly a morphism
  $$
   \tilde{f}\;:\;
     {\Bbb A}^{n^2}_{(0)}\times{\Bbb A}^{n^2}_{(1)}
                   \times\,\cdots\,\times {\Bbb A}^{n^2}_{(l)}\;
   \longrightarrow\;
   {\Bbb A}^{n^2}_{(0)}\times{\Bbb A}^{n^2}_{(1)}\times\,\cdots\,
      \times {\Bbb A}^{n^2}_{(l^{\prime})}
  $$
  that restricts to a morphism
  $$
   f\;:\;
    \Rep^{\ringsetscriptsize}(R,M_n(\Bbb C))\;
     \longrightarrow\;
     \Rep^{\ringsetscriptsize}(R,M_n(\Bbb C))^{\prime}\,,
  $$
  where
   $\Rep^{\ringsetscriptsize}(R,M_n(\Bbb C))^{\prime}
     \subset
     {\Bbb A}^{n^2}_{(0)}\times{\Bbb A}^{n^2}_{(1)}\times\,\cdots\,
            \times {\Bbb A}^{n^2}_{(l^{\prime})}$
   is the representation variety associated to the new presentation
   of $R$.
 Reverse this argument, now from
   $\langle g_0, g_1,\,\cdots\,, g_l\rangle /(r_1,\,\cdots\,, r_m)$
  to
   $\langle g^{\prime}_0, g^{\prime}_1,\,\cdots\,,
              g^{\prime}_{l^{\prime}}\rangle /
       (r^{\prime}_1,\,\cdots\,, r^{\prime}_{m^{\prime}})$,
  implies that $f$ is indeed an isomorphism.

 Since a ring-isomorphism sends the identity to the identity,
  $f$ restricts to isomorphisms
  $f_{(d)}:\Rep^{\ringsetscriptsize}(R,M_n(\Bbb C))_{(d)}
           \stackrel{\sim}{\rightarrow}
           \Rep^{\ringsetscriptsize}(R,M_n(\Bbb C))^{\prime}_{(d)}$,
  for $d=0,\,\cdots\,,n$.
 In other words,
  $f_{(d)}:\Mor^{\ringsetscriptsize}(R,M_n(\Bbb C))^{{\cal T}_0}_{(d)}
           \stackrel{\sim}{\rightarrow}
           \Mor^{\ringsetscriptsize}(R,M_n(\Bbb C))
                                  ^{\prime\,{\cal T}^{\prime}_0}_{(d)}$,
  for $d=0,\,\cdots\,,n$.
 Furthermore, each valuative criterion setup
   $\Phi_{\phi}:
     (D,p)\dottedrightarrow \Mor^{\ringsetscriptsize}(R,M_n(\Bbb C))$
   gives a valuative criterion setup
  $\Phi_{\phi^{\prime}}= f\circ \Phi_{\phi}:
    (D,p)\dottedrightarrow
     \Mor^{\ringsetscriptsize}(R,M_n(\Bbb C))^{\prime}$
 and vice versa.
 As we choose the topology ${\cal T}$ on
  $\Mor^{\ringsetscriptsize}(R,M_n(\Bbb C))$
  (resp.\ ${\cal T}^{\prime}$ on
          $\Mor^{\ringsetscriptsize}(R,M_n(\Bbb C))^{\prime}$)
  to be the weakest topology that renders
   all inclusions
    $\Mor^{\ringsetscriptsize}(R,M_n(\Bbb C))^{{\cal T}_0}_{(d)}
     \hookrightarrow 
     \Mor^{\ringsetscriptsize}(R,M_n(\Bbb C))$
    (resp.\
    $\Mor^{\ringsetscriptsize}(R,M_n(\Bbb C))^{\prime\,{\cal T}_0}_{(d)}
     \hookrightarrow
     \Mor^{\ringsetscriptsize}(R,M_n(\Bbb C))^{\prime}$)
   embeddings of topological spaces and
  all $\Phi_{\phi}$'s (resp.\ $\Phi_{\phi^{\prime}}$'s) continuous,
 this implies that
  $$
   f\; :\;  (\Mor^{\ringsetscriptsize}(R,M_n(\Bbb C)),{\cal T})\;
      \longrightarrow\;
       (\Mor^{\ringsetscriptsize}(R,M_n(\Bbb C))^{\prime},
                                                {\cal T}^{\prime})
  $$
  is an isomorphism.
 This completes the proof.

\noindent\hspace{15cm}$\Box$

\bigskip

By construction, there is a canonical (continuous) bijective embedding
$$
 \tau_{R,n}\,:\;
  \Rep^{\ringsetscriptsize}(R,M_n(\Bbb C))_{\scriptsizeBbb C}\;
  \longrightarrow\; \Mor^{\ringsetscriptsize}(R,M_n(\Bbb C))\,.
$$

\bigskip

\noindent
{\it Remark 3.2.8 $[$moduli problem$]$.}
By construction, $\Mor^{\ringsetscriptsize}(R,M_n(\Bbb C))$
 is a coarse moduli space of ring-set-homomorphisms
 from $R$ to $M_n({\Bbb C})$.
Since a ring-set-homomorphism with a fixed domain and target
 does not have non-trivial automorphisms,
it is instructive to think of
 $\Mor^{\ringsetscriptsize}(R,M_n(\Bbb C))$
 as representing the functor
 $$
  \begin{array}{ccccc}
   {\cal F} & :
    & \left(\,\mbox{\parbox{27ex}{(commutative) varieties/${\Bbb C}$\\
                with analytic topology}}\,\right)^{\circ}
      & \longrightarrow  & (\,\mbox{sets}\,) \\[2.6ex]
   &&
    V & \longmapsto
      & \Mor_{{\cal O}_V\mbox{\scriptsize\it -Alg}}
         ( {\cal O}_V\otimes R\,,\,
           {\cal O}_V\otimes M_n({\Bbb C}) )
  \end{array}
 $$
 similar to a functor of points.
Here, $(\,\cdots\,)^{\circ}$ is the category $(\,\cdots\,)$
 with the arrows reversed.

\bigskip

\begin{flushleft}
{\bf $\Mor(\Space M_n({\largeBbb C}),Y)$ as a coarse moduli space.}
\end{flushleft}
Let $Y$ be a noncommutative space presented as a gluing system of
 finitely-presentable rings
 ${\cal R}=
   (\{R_{\alpha}\}_{\alpha\in A}
     \doublearrow \{R_{\alpha_1\alpha_2}\}_{\alpha_1,\alpha_2\in A})$.
We fix a well-ordering of the index set $A$ for convenience.
Denote the identity of $R_{\alpha}$ by $1_{R_{\alpha}}$.
Assume that each central localization
 $\varphi_{\alpha_1\alpha_2}:R_{\alpha_1}\rightarrow R_{\alpha_1\alpha_2}$
 is associated to a finitely-generated multiplicatively closed subset
 $S_{\alpha_1\alpha_2}$ in $Z(R_{\alpha_1})$.

\bigskip

\noindent
{\bf Definition 3.2.9 [admissible tuple].} {\rm
 A tuple
  $(\varphi_{\alpha}:R_{\alpha}\rightarrow M_n({\Bbb C}))_{\alpha\in A}$
  of ring-set-homomorphisms to $M_n({\Bbb C})$ is called
  {\it admissible}
 if it satisfies the following conditions:
  \begin{itemize}
   \item[$(1)$]
    $\varphi_{\alpha_1}(1_{R_{\alpha_1}})$ commutes with
     elements of $\varphi_{\alpha_2}(R_{\alpha_2})$
     for all $\alpha_1$, $\alpha_2\in A$.

   \item[$(2)$]
    $\varphi_{\alpha_1}(1_{R_{\alpha_1}})\,
         \varphi_{\alpha_2}(R_{\alpha_2})
     = \varphi_{\alpha_2}(1_{R_{\alpha_2}})\,
         \varphi_{\alpha_1} (R_{\alpha_1})$
     for all $\alpha_1$, $\alpha_2\in A$.

   \item[$(3)$]
    $\varphi_{\alpha_1}(1_{R_{\alpha_1}})\,
         \varphi_{\alpha_2}(1_{R_{\alpha_2}})
     \in \varphi_{\alpha_1}(R_{\alpha_1})$
    for all $\alpha_1$, $\alpha_2\in A$.

   \item[$(4)$]
    Let ${\mathbf 1}$ be the identity matrix in $M_n({\Bbb C})$.
    Then
    \begin{eqnarray*}
     \lefteqn{
      {\mathbf 1}\;=\; \sum_{\alpha} \varphi_{\alpha}(1_{R_{\alpha}})\;
           -\, \sum_{\alpha_1<\alpha_2}
                 \varphi_{\alpha_1} (1_{R_{\alpha_1}})\,
                 \varphi_{\alpha_2} (1_{R_{\alpha_2}})\;  } \\[.6ex]
     && \hspace{6.2em}
        +\, \sum_{\alpha_1<\alpha_2<\alpha_3}
             \varphi_{\alpha_1} (1_{R_{\alpha_1}})\,
             \varphi_{\alpha_2} (1_{R_{\alpha_2}})\,
             \varphi_{\alpha_3} (1_{R_{\alpha_3}})          \\[.6ex]
     && \hspace{6.2em}
        \pm\; \cdots\;
          +\, (-1)^{|A|+1}\,
               \sum_{\alpha_1<\,\cdots\,<\alpha_{|A|}}
                 \varphi_{\alpha_1} (1_{R_{\alpha_1}})\,\cdots\,
                 \varphi_{\alpha_{|A|}} (1_{R_{\alpha_{|A|}}})\,.
    \end{eqnarray*}

   \item[$(5)$]
    $\varphi_{\alpha_2}(1_{R_{\alpha_2}})
     \cdot \left( \varphi_{\alpha_1}(s)^{\perp}
                  \cap \varphi_{\alpha_1}(R_{\alpha_1}) \right) = 0$,
    where
     $\varphi_{\alpha_1}(s)^{\perp}
      := \{ m\in M_n({\Bbb C}): \varphi_{\alpha_1}(s)\cdot m =0 \}$,
     for all $\alpha_1$, $\alpha_2\in A$  and
             $s\in S_{\alpha_1\alpha_2}$.
    This condition is equivalent to the existence of push-out
     $\varphi_{\alpha_1}|_{\alpha_2}$ under localizations
     in the following commutative diagram:
     $$
      \begin{array}{lccl}
       R_{\alpha_1}  & \stackrel{\varphi_{\alpha_1}}{\longrightarrow}
                     & \varphi_{\alpha_1}(R_{\alpha_1})     \\
       \hspace{1ex}\downarrow    & & \downarrow             \\[-1ex]
       R_{\alpha_1\alpha_2}
         & \stackrel{\varphi_{\alpha_1}|_{\alpha_2}}{\longrightarrow}
         & \varphi_{\alpha_2}(1_{R_{\alpha_2}})
            \cdot \varphi_{\alpha_1}(R_{\alpha_1}) &.
      \end{array}
     $$

   \item[$(6)$]
    $\varphi_{\alpha_1}|_{\alpha_2} =
     \varphi_{\alpha_2}|_{\alpha_1} \circ \varphi_{\alpha_1\alpha_2}$
     %
     %
    for all $\alpha_1$, $\alpha_2\in A$.
  \end{itemize}
} 

\bigskip

\noindent
The meaning of these conditions is given below.
\begin{itemize}
 \item[$\cdot$]
  Conditions (1) - (4):
  The finite collection
   $\{ ( \varphi_{\alpha}(R_{\alpha}),
         e_{\alpha}:=\varphi_{\alpha}(1_{R_{\alpha}}) )
      \}_{\alpha\in A}$
   of ring-subsets of $M_n({\Bbb C})$
   glue to $\sum_{\alpha\in A}\varphi_{\alpha}(R_{\alpha})$
   that is a subalgebra of $M_n({\Bbb C})$.
   Cf.\ Proposition 3.2.4.

 \item[$\cdot$]
  Condition (5):
   Elements in $\varphi_{\alpha_1}(S_{\alpha_1\alpha_2})$
    become invertible after being mapped to
    $e_{\alpha_2}\cdot \varphi_{\alpha_1}(R_{\alpha_1})$  and,
   hence,
   $\varphi_{\alpha_1}$ can be pushed out to a ring-homomorphism
    $\varphi_{\alpha_1}|_{\alpha_2}$ from $R_{\alpha_1\alpha_2}$
    to the localization
    $e_{\alpha_2} \cdot \varphi_{\alpha_1}(R_{\alpha_1})$
    of $\varphi_{\alpha_1}(R_{\alpha_1})$.
   Cf.\ Lemma 3.1.3.

 \item[$\cdot$]
  Condition (6):
   The gluing conditions on the tuple
    $\{\varphi_{\alpha}:R_{\alpha}\rightarrow M_n({\Bbb C})\}
       _{\alpha\in A}$
    as a system of ring-homomorphisms from ${\cal R}$
    to $( \{\varphi_{\alpha}(R_{\alpha})\}_{\alpha\in A}
          \doublearrow
          \{ e_{\alpha_2} \cdot\varphi_{\alpha_1}(R_{\alpha_1}) \}
           _{\alpha_1,\alpha_2\in A} )$.
  Cf.\ Condition (2) above and Definition 1.2.6.
\end{itemize}

Thus, Conditions (1) - (6) are necessary conditions for the tuple
  $\{\varphi_{\alpha}\}_{\alpha\in A}$ to represent a morphism
  from $(\Spec{\Bbb C}, {\Bbb C}, M_n({\Bbb C}))$ to $Y$.
It follows from Definition 1.2.14
 that they are also sufficient and
 that such presentations are effective in the sense that
  different admissible tuples give different morphisms.
This proves the following lemma:

\bigskip

\noindent
{\bf Lemma 3.2.10 [admissible tuple = morphism].} {\it
 A tuple
  $\Phi =
   (\varphi_{\alpha}:R_{\alpha}\rightarrow M_n({\Bbb C}))_{\alpha\in A}$
  of ring-set-homomorphisms to $M_n({\Bbb C})$ corresponds to
  a morphism from $\Space M_n({\Bbb C})$ to $Y=\Space {\cal R}$
 if and only if $\Phi$ is admissible.
 As sets,
  $\Mor(\Space M_n(\Bbb C), Y) = \{\,\mbox{admissible tuples}\,\}$.
} 

\bigskip

Fix now the following data of presentations and representatives:
 \begin{itemize}
  \item[$\cdot$] [{\it ring chart}]\\  
   a finite presentation for each ring-chart $R_{\alpha}$ in ${\cal R}$
    $$
     R_{\alpha}\; =\;
      \langle g^{(\alpha)}_0, g^{(\alpha)}_1,\,\cdots\,,
             g^{(\alpha)}_{l^{(\alpha)}}\rangle
      /(r^{(\alpha)}_1,\,\cdots\,, r^{(\alpha)}_{m^{(\alpha)}})\,,
    $$
    with the redundant generator $g^{(\alpha)}_0 = 1_{R_{\alpha}}$ and
    the redundant relators
     $$
      g^{(\alpha)}_0g^{(\alpha)}_i\;
      =\; g^{(\alpha)}_ig^{(\alpha)}_0\; =\; g^{(\alpha)}_i\,,
      \hspace{1ex}
      i\,=\,0,\,1,\,\ldots\,,l^{(\alpha)}\,,
     $$
     contained in the relator set
     $\{r^{\alpha}_1,\,\cdots\,, r^{\alpha}_{m^{\alpha}}\}$,
    as before;

  \item[$\cdot$] [{\it localization}]\\  
   a lifting (as sets) $\tilde{S}_{\alpha_1\alpha_2}$
    of $S_{\alpha_1\alpha_2}$
    in $\langle g^{(\alpha_1)}_0, g^{(\alpha_1)}_1,\,
                \cdots\,, g^{(\alpha_1)}_{l^{(\alpha_1)}}\rangle$
    for each $(\alpha_1, \alpha_2)\in A\times A$;

  \item[$\cdot$] [{\it transition data}]\\  
   a representative
    in the induced presentation of $R_{\alpha_2\alpha_1}$
    for each $g^{(\alpha_1)}_i$,
         $i=0,\,\ldots\,, l^{(\alpha_1)}$,
        $\tilde{s}\in \tilde{S}_{\alpha_1\alpha_2}$,
        and $(\alpha_1, \alpha_2)\in A\times A$,
    $$
     (g^{(\alpha_1\alpha_2)}_i, s^{(\alpha_1\alpha_2)}_i)\;,\;\;
     (g^{(\alpha_1\alpha_2)}_{\tilde{s}},
                     s^{(\alpha_1\alpha_2)}_{\tilde{s}})\;
      \in\;  \langle g^{(\alpha_2)}_0, g^{(\alpha_2)}_1,\,\cdots\,,
                     g^{(\alpha_2)}_{l^{(\alpha_2)}}\rangle
             \times \tilde{S}_{\alpha_2\alpha_1}
    $$
   so that
    $\varphi_{\alpha_1\alpha_2}(g^{(\alpha_1)}_i, 1_{R_{\alpha_1}})
     =(g^{(\alpha_1\alpha_2)}_i,s^{(\alpha_1\alpha_2)}_i)$
      and
    $\varphi_{\alpha_1\alpha_2}(1_{R_{\alpha_1}},\tilde{s})
     =(g^{(\alpha_1\alpha_2)}_{\tilde{s}},
                       s^{(\alpha_1\alpha_2)}_{\tilde{s}})$.
   (Here, to simplify notations, we identify elements in
    a presentation of a ring with the corresponding elements
    in that ring.)
 \end{itemize}
Let
 \begin{itemize}
  \item[$\cdot$]
   ${\Bbb A}^{n^2}_{(\alpha,\,i)}$,
    $\alpha\in A$, $i=0,\,\ldots\,,l^{(\alpha)}$,
    be the affine space with the polynomial coordinate ring
    ${\Bbb C}[m^{(\alpha)}_{i, jk}: 1\le j, k\le n]$;

  \item[$\cdot$]
   ${\mathbf A}_{\alpha}$
    be the affine space
    ${\Bbb A}^{n^2}_{(\alpha, 0)}\times{\Bbb A}^{n^2}_{(\alpha, 1)}
     \times\,\cdots\, \times {\Bbb A}^{n^2}_{(\alpha,\,l^{(\alpha)})}$
     and
   ${\mathbf A}$ be the affine space
    $\prod_{\alpha\in A}{\mathbf A}_{\alpha}
     =\; {\Bbb A}^{\sum_{\alpha\in A}(1+l^{(\alpha)})\,n^2}$;

  \item[$\cdot$]
   $R({\mathbf A}_{\alpha})
    := {\cal O}_{{\mathbf A}_{\alpha}}({\mathbf A}_{\alpha})
     = \otimes_{i=0}^{l^{(\alpha)}}
        {\Bbb C}[m^{(\alpha)}_{i, jk}: 1\le j, k\le n]$
      and
   $R({\mathbf A}) :={\cal O}_{\mathbf A}({\mathbf A})
    = \otimes_{\alpha\in A} R({\mathbf A}_{\alpha})$;

  \item[$\cdot$]
   $\Psi_{\alpha}:
     R({\mathbf A})\otimes_{\scriptsizeBbb C}
        \langle g^{(\alpha)}_0, g^{(\alpha)}_1,\,\cdots\,,
                g^{\alpha}_{l^{(\alpha)}}\rangle
     \rightarrow
     R({\mathbf A})\otimes_{\scriptsizeBbb C} M_n({\Bbb C})
      = M_n(R({\mathbf A}))$
   be the tautological $R({\mathbf A})$-algebra-homomorphism
   defined/generated by\footnote{Recall
                                 that the multiplication $\cdot$ in
                                 the tensor product
                                 ${\footnotesizeBbb C}$-algebra
                                 $R\otimes_{\tinyBbb C}S$
                                 of two ${\footnotesizeBbb C}$-algebras
                                 $R$ and $S$ is
                                 ${\footnotesizeBbb C}$-linearly
                                 generated by defining
                                 $(r_1\otimes s_1)\cdot(r_2\otimes s_2)
                                  = (r_1r_2)\otimes (s_1s_2)$.}
   $$
    1\otimes g^{(\alpha)}_i\;
     \longmapsto\; 1\otimes \left(m^{(\alpha)}_{i,jk}\right)_{jk}
   $$
    and
   $\image\Psi_{\alpha}$ be the image $R({\mathbf A})$-submodule
    of $\Psi_{\alpha}$ in $M_n(R({\mathbf A}))$;

  \item[$\cdot$]
   $E_{\mathbf A}={\mathbf A}\times M_n({\Bbb C})$
    be the trivialized trivial vector bundle on ${\mathbf A}$
     with fiber the ${\Bbb C}$-algebra $M_n({\Bbb C})$;
   the associated sheaf of local sections of $E_{\mathbf A}$
    is ${\cal O}_{\mathbf A}\otimes M_n({\Bbb C})$;
   elements and sub-$R({\mathbf A})$-modules in
    ${\cal O}_{\mathbf A}\otimes M_n({\Bbb C})$
    are canonically identified respectively with
    global sections and constructible sets in $E_{\mathbf A}$.
 \end{itemize}
Define
 $$
  \Rep^{\ringsetscriptsize}({\cal R},M_n(\Bbb C))\;
   \subset\;
   \prod_{\alpha\in A}
     \Mor^{\ringsetscriptsize}(R_{\alpha},M_n(\Bbb C))
 $$
  to be the locus in the indicated product space
  determined by the following system of constraints
  from the defining conditions of admissible tuples,
 via the canonical bijective embedding
 $$
  \prod_{\alpha\in A}
    \Mor^{\ringsetscriptsize}(R_{\alpha},M_n(\Bbb C))\;
  \stackrel{\prod_{\alpha\in A}\tau_{R_{\alpha},n}}
           {\longleftarrow\hspace{-1ex}\mbox{------------}}\;
  \prod_{\alpha\in A}
    \Rep^{\ringsetscriptsize}(R_{\alpha},M_n(\Bbb C))
                                        _{\scriptsizeBbb C}\;
  \subset\; {\mathbf A}\,:
 $$

\medskip

\begin{itemize}
 \item[] \parbox[t]{3em}{(0.1)}
  \vspace{-2.4em}
  \begin{eqnarray*}
   \lefteqn{\hspace{3.6em}
   r^{(\alpha)}_1(M_{\alpha,0}, M_{\alpha,1},\,\cdots\,,
                M_{\alpha,\,l^{(\alpha)}})\;
    =\; \cdots\;
    =\; r^{(\alpha)}_{m^{(\alpha)}}
          (M_{\alpha,0}, M_{\alpha, 1},\,
                 \cdots\,, M_{\alpha,\,l^{(\alpha)}})\; } \\[,6ex]
   && \hspace{4em}
      =\; \mbox{the zero-matrix $0\,\in\, M_n({\Bbb C})$,
                where
                 $M_{\alpha,\,i}
                  =\left(m^{(\alpha)}_{i,jk}\right)_{jk}$.}
          \hspace{6em}
  \end{eqnarray*}

 \item[] \parbox[t]{3em}{(1.1)}
  $M_{\alpha_1,0}\, M_{\alpha_2, i}=M_{\alpha_2,i}\, M_{\alpha_1,0}\;$
   for all $\alpha_1$, $\alpha_2\in A$, $i=0,\,\ldots\,, l^{(\alpha_2)}$.

 \item[] \parbox[t]{3em}{(1.2)}
  $M_{\alpha_1,0}\, \image\Psi_{\alpha_2}
                 = M_{\alpha_2,0}\, \image\Psi_{\alpha_1}$
   for all $\alpha_1$, $\alpha_2\in A$.

 \item[] \parbox[t]{3em}{(1.3)}
  $M_{\alpha_1,0}\,
       M_{\alpha_2,0} \in \image\Psi_{\alpha_1}$
  for all $\alpha_1$, $\alpha_2\in A$.

 \item[] \parbox[t]{3em}{(1.4)}
  (${\mathbf 1}\in M_n({\Bbb C})$ is the identity)
  \begin{eqnarray*}
   \lefteqn{\hspace{2em}
    {\mathbf 1}\;=\; \sum_{\alpha} M_{\alpha,0}\;
      -\, \sum_{\alpha_1<\alpha_2}
            M_{\alpha_1,0}\, M_{\alpha_2,0}\,
      +\, \sum_{\alpha_1<\alpha_2<\alpha_3}
            M_{\alpha_1,0}\, M_{\alpha_2,0}\, M_{\alpha_3,0} }\\[.6ex]
   && \hspace{10em}
      \pm\; \cdots\;
        +\, (-1)^{|A|+1}\,
             \sum_{\alpha_1<\,\cdots\,<\alpha_{|A|}}
               M_{\alpha_1,0}\,\cdots\, M_{\alpha_{|A|,0}}\,.
  \end{eqnarray*}

 \item[] \parbox[t]{3em}{(1.5)}
  $M_{\alpha_2,0}
   \cdot \left( \Psi_{\alpha_1}(\tilde{s})^{\perp}_{E_{\mathbf A}}
                \cap \image\Psi_{\alpha_1} \right) = 0$
   for all $\alpha_1$, $\alpha_2\in A$  and
           $\tilde{s}\in \tilde{S}_{\alpha_1\alpha_2}$.
  Here\footnote{Caution that $\Psi_{\alpha_1}(\tilde{s})^{\perp}
                                                       _{E_{\mathbf A}}$
                                          here is defined to be
                                          the union of fiberwise $\perp$
                                          of $\Psi_{\alpha_1}(\tilde{s})$
                                          in $E_{\mathbf A}$.
                                         In general, it is not
                                          a sub-$R({\mathbf A})$-module
                                          of $M_n(R({\mathbf A}))$.}
   $$
    \Psi_{\alpha_1}(\tilde{s})^{\perp}_{E_{\mathbf A}}\;
     :=\;
     \{ m\in E_{\mathbf A}: \Psi_{\alpha_1}(\tilde{s})\cdot m =0 \}\,.
   $$

 \item[] \parbox[t]{3em}{(1.6)}
  \vspace{-1.4em}
  \begin{eqnarray*}
   \lefteqn{
     \left( M_{\alpha_2,0}\, M_{\alpha_1, i} \rule{0ex}{2ex}\right)\,
     \left( M_{\alpha_1,0}\;
       s^{(\alpha_1\alpha_2)}_i
        (M_{\alpha_2,0},\,\cdots\,, M_{\alpha_2,\,l^{(\alpha_2)}})
      \right)               }\\[.6ex]
     && =\; M_{\alpha_1,0}\:
            g^{(\alpha_1\alpha_2)}_i
            (M_{\alpha_2,0},\,\cdots\,, M_{\alpha_2,\,l^{(\alpha_2)}})
         \hspace{11.2em}
  \end{eqnarray*}
  $\mbox{\hspace{3.2em}}$ and
  \begin{eqnarray*}
   \lefteqn{
    \left( M_{\alpha_2,0}\,M_{\alpha_1,0} \rule{0ex}{2ex}\right)\,
    \left( M_{\alpha_1,0}\;
           s^{(\alpha_1\alpha_2)}_{\tilde{s}}
           (M_{\alpha_2,0},\,\cdots\,, M_{\alpha_2,\,l^{(\alpha_2)}})
      \right)                 }\\[.6ex]
    && =\;
    \left( M_{\alpha_1,0}\:
      g^{(\alpha_1\alpha_2)}_{\tilde{s}}
       (M_{\alpha_2,0},\,\cdots\,, M_{\alpha_2,\,l^{(\alpha_2)}}
     \right)\,
    \left( M_{\alpha_2,0}\;
      \tilde{s}
       (M_{\alpha_1,0},\,\cdots\,, M_{\alpha_1,\,l^{(\alpha_1)}})
      \right)\,
  \end{eqnarray*}
  $\mbox{\hspace{3.2em}}$
  for all $\alpha_1$, $\alpha_2\in A$,
      $i=0,\,\ldots\,,l^{(\alpha_1)}$, and
      $\tilde{s}\in \tilde{S}_{\alpha_1\alpha_2}$.
\end{itemize}

\bigskip

\noindent
{\bf Proposition 3.2.11 [$\Mor(\Space M_n({\largeBbb C}),Y)$].} {\it
 $\Mor(\Space M_n({\largeBbb C}),Y)$
  is given by a constructible set in the product space
 $\prod_{\alpha\in A}
       \Mor^{\ringsetscriptsize}(R_{\alpha},M_n(\Bbb C))$,
 independent of the data of presentation chosen in the construction.
} 

\bigskip

\noindent
{\it Proof.}
 Conditions (0.1), (1.1), and (1.4) are closed conditions.
 Condition (1.6) can be restricted to a finite generating
  set of $S_{\alpha_1\alpha_2}$ and, hence, gives
  also a closed condition.
 Conditions (1.2), (1.3), and (1.5) involve image
  $R({\mathbf A})$-submodules $\image\Psi_{\bullet}$
  in $M_n(R({\mathbf A}))$.
 Let $S_{\alpha_1\alpha_2}^0$ be a finite generating set
   of $S_{\alpha_1\alpha_2}$  and
  $\tilde{S}_{\alpha_1\alpha_2}^0
     \subset \tilde{S}_{\alpha_1\alpha_2}$
    its corresponding lifting in
    $\langle g^{(\alpha_1)}_0, g^{(\alpha_1)}_1,\,
               \cdots\,, g^{(\alpha_1)}_{l^{(\alpha_1)}}\rangle$.
 Then, it follows from Lemma 3.1.7
    and
   the fact that every nilpotent element $m$ of $M_n({\Bbb C})$
    satisfies $m^n=0$
 that the seemingly possibly-infinite system of constraints from
  Condition (1.5) can be replaced by the following finite system:
  \begin{itemize}
   \item[] \parbox[t]{3em}{(1.5)$^{\prime}$}
    $M_{\alpha_2,0}
     \cdot \left( \Psi_{\alpha_1}(\tilde{s}^n)^{\perp}_{E_{\mathbf A}}
                  \cap \image\Psi_{\alpha_1} \right) = 0$
     for all $\alpha_1$, $\alpha_2\in A$  and
            $\tilde{s}\in \tilde{S}_{\alpha_1\alpha_2}^0$.
 \end{itemize}
 Thus, the solution set to Conditions (1.2), (1.3), and (1.5)
  is described by a finite intersection of constructible sets on
  ${\mathbf A}$ described via determinantal varieties.

 This shows that the solution set to the system of constraints from
  Condition (0.1) and Conditions (1.1) - (1.5) is a constructible set
  in $\prod_{\alpha\in A}
       \Rep^{\ringsetscriptsize}(R_{\alpha},M_n(\Bbb C))$  and,
 hence, in
  $\prod_{\alpha\in A}
    \Mor^{\ringsetscriptsize}(R_{\alpha},M_n(\Bbb C))$.
 That different choices of data of presentations give isomorphic
  solution sets (with the subset topology) follows the same
  discussion as that in the proof of Proposition 3.2.7.
 Since
  $\Mor(\Space M_n({\Bbb C}), Y)
   =\Rep^{\ringsetscriptsize}({\cal R},M_n({\Bbb C}))$ as sets,
 this concludes the proof.

\noindent\hspace{15cm}$\Box$

\bigskip

We remark that from the proof above, the constructible set referred
 to in Proposition 3.2.11 is of algebraic kind.
It is the set of ${\Bbb C}$-points (with the analytic topology)
 of a finite union of constructible sets in varieties/${\Bbb C}$.

Finally, note that in discussing the space of morphisms
  from $\Space M_n({\Bbb C})$ to $Y$,
 both $\Space M_n({\Bbb C})$ and $Y$ are thought of as fixed.
The automorphism group of $M_n({\Bbb C})$
 as a ${\Bbb C}$-algebra is given by $\GL_n({\Bbb C})$
 via the adjoint $\GL_n({\Bbb C})$-action on $M_n({\Bbb C})$.
This induces a $\GL_n({\Bbb C})$-action on
 $\Mor(\Space M_n({\Bbb C}), Y)$.

\bigskip

\noindent
{\bf Definition 3.2.12 [isomorphism between morphisms].} {\rm
 Two morphisms from $\Space M_n({\Bbb C})$ to $Y$ are said to be
  {\it isomorphic}, in notation $\Phi_1\sim \Phi_2$,
 if they are in the same $\GL_n({\Bbb C})$-orbit in
  $\Mor(\Space M_n({\Bbb C}), Y)$.
 Define the {\it space $\Map(\Space M_n({\Bbb C}), Y)$ of maps}
  from $\Space M_n({\Bbb C})$ to $Y$ to be the quotient space
  $\Mor(\Space M_n({\Bbb C}), Y)/\!\sim\,$ (with the quotient topology).
 It parameterizes isomorphism classes of morphisms from
  $\Space M_n({\Bbb C})$ to $Y$.
} 

\bigskip

\section{D0-branes on a commutative quasi-projective variety.}

A D0-brane in the sense of Definition 2.2.3
 is simply an Azumaya-type noncommutative point
 $\Space M_n({\Bbb C})$ (cf.\ Example 1.1.3 and Example 1.1.8)
 together with the irreducible $M_n({\Bbb C})$-module ${\Bbb C}^n$
 as the Chan-Paton space/module.
A D0-brane on a target space $Y$ is given by an isomorphism class
 of morphisms from $\Space M_n({\Bbb C})$ to $Y$.
The moduli space of D0-branes on $Y$ in this sense is given then by
 $\Map((\Space M_n({\Bbb C});{\Bbb C}^n),Y)
  =\Map(\Space M_n({\Bbb C}),Y) = \Mor(\Space M_n({\Bbb C}),Y)/\!\sim$.
This moduli space for the case of $Y$ being a (commutative)
 complex quasi-projective smooth curve/surface, or a variety
 is given in this section to illustrate Sec.~1 - Sec.~3.
These examples already reveal simplified key features of D-branes
 that are fundamental for beyond.
Details involving only linear algebras in, e.g., [Ho-K] or
 straightforward manipulations are omitted.
\bigskip

\subsection{D0-branes on the complex affine line ${\Bbb A}^1$.}

Various themes concerning D0-branes on ${\Bbb A}^1$ are given
 in this subsection to illustrate the far-reaching/power of
 the Polchinski-Grothendieck Ansatz for D-branes,
 in particular the reproduction
 of D-brane properties in the work of Polchinski.
Same/Similar phenomena occur also for other targets
 in later subsections by same/similar reasons,
 which we then omit but focus mainly on the moduli problem.
The general discussions in Sec.~1 - Sec.~3 are intentionally
 made explicit in this example.
For that reason, some important algebro-geometric notions
 are slightly repeated in this subsection for concreteness.

\bigskip

\begin{flushleft}
{\bf The moduli space
     $\Map((\Space M_n({\Bbb C});{\Bbb C}^n)\,,\,{\Bbb A}^1)$
     of D0-branes on ${\Bbb A}^1$.}
\end{flushleft}
Let $Y={\Bbb A}^1=\Spec {\Bbb C}[y]$ be the affine line over ${\Bbb C}$.
Then the Grothendieck Satz or Lemma 1.2.19 says that
 $\Mor(\Space M_n({\Bbb C}), Y)=\Mor({\Bbb C}[y], M_n({\Bbb C}))$.
The corresponding ${\Bbb C}$-algebra representation variety
 $\Rep({\Bbb C}[y], M_n({\Bbb C}))$ is given by ${\Bbb A}^{n^2}$
 with a closed point represented by $m=(m_{ij})_{i,j}\in M_n({\Bbb C})$
  corresponding to the ${\Bbb C}$-algebra-homomorphism
  $$
   \varphi_m: {\Bbb C}[y]\rightarrow M_n({\Bbb C})\,,
   \hspace{1ex}
   \mbox{generated by $1\mapsto{\mathbf 1}$ and $y\mapsto m$}\,.
  $$
We will call the $\GL_n({\Bbb C})$-action on
 $\Rep({\Bbb C}[y], M_n({\Bbb C}))$ by post-compositions with
 the conjugations on $M_n({\Bbb C})$ still the {\it adjoint action}.
It follows that
 $$
  \Map((\Space M_n({\Bbb C});{\Bbb C}^n), {\Bbb A}^1)\;
  =\; \Map(\Space M_n({\Bbb C}), {\Bbb A}^1)\;
  =\; \Rep({\Bbb C}[y],M_n({\Bbb C}))/\!\sim\,,
 $$
 the orbit-space\footnote{We shall always think of such
                                     an orbit-space $M/\!\sim$
                                     as an Artin stack with atlas $M$.
                                    When $M$ is smooth,
                                     it is in this sense that we define
                                     a smooth map to $M/\!\sim$.}
 of the adjoint action with the quotient topology.
This space is a connected non-Hausdorff topological space,
 well-understood in other contents from algebraic geometry and
 Lie groups and Lie algebras as follows.

Each adjoint-orbit $O_{\varphi_m}$ is represented by
 a Jordan form $J_m$ of $m$, unique up to permutations of
 diagonal blocks of $J_m$ with distinct characteristic values.
An adjoint-orbit on $\Rep({\Bbb C}[y], M_n({\Bbb C}))$ is closed
 if and only of it is represented by
 $\varphi_m$ associated to a diagonal matrix $m$.
Given an orbit $O_{\varphi_m}$, let $\overline{O_{\varphi_m}}$
 be the closure of $O_{\varphi_m}$ in ${\Bbb A}^{n_2}$.
It has the property that
 $O_{\varphi_m}$ is an open dense subset in $\overline{O_{\varphi_m}}$
  and that
 $\overline{O_{\varphi_m}}$ is a union of $O_{\varphi_m}$ and
  finitely many lower-dimensional orbits, (e.g.\ [Stei]).
Note that any two orbits $O_{\varphi_{m_1}}$ and $O_{\varphi_{m_2}}$
 satisfy
  either $O_{\varphi_{m_1}}\cap \overline{O_{\varphi_{m_2}}}=\emptyset$
  or     $O_{\varphi_{m_1}}\subset \overline{O_{\varphi_{m_2}}}$.

\bigskip

\noindent
{\bf Definition 4.1.1
     [partial order on $\Rep({\Bbb C}[y], M_n({\Bbb C}))/\!\sim$].}
{\rm
 Define a {\it partial order} on the orbit-space
  $\Rep({\Bbb C}[y], M_n({\Bbb C}))/\!\sim$
  by setting $O_{\varphi_{m_1}}\prec O_{\varphi_{m_2}}$
  if $O_{\varphi_{m_1}}\subset\overline{O_{\varphi_{m_2}}}$.
} 

\bigskip

\noindent
This partial order can be described in terms of Jordan forms,
 as follows.

Let $J^{(\lambda)}_j\in M_j({\Bbb C})$ be the matrix
{\scriptsize
 $$
  \left[
   \begin{array}{cccc}
    \lambda   &          &        & 0       \\
    1         & \lambda  &        &         \\
              &  \ddots  & \ddots &         \\
    0         &          &   1    & \lambda
   \end{array}
  \right]_{j\times j}
 $$
{\normalsize A}}
Jordan form $J$ in $M_n({\Bbb C})$ is a matrix of the following form
$$
 \left[
  \begin{array}{ccc}
   A_1 &         & 0 \\
       & \ddots  &   \\
   0   &         & A_k
  \end{array}
 \right]
  \hspace{1em}
   \mbox{with each $A_i\in M_{n_i}({\Bbb C})$ of the form}\hspace{1em}
 \left[
  \begin{array}{ccc}
   J^{(\lambda_i)}_{d_{i1}} &         & \\
             & \ddots  &              \\
             &         &  J^{(\lambda_i)}_{d_{ik_i}}
  \end{array}
 \right]\,.
$$
Here, omitted entries are all zero,
 $n_1\ge \,\cdots\,\ge n_k> 0$, and $d_{i1}\ge\,\cdots\,\ge d_{ik_i}>0$.
We thus have a {\it double partition} of $n$
 by non-increasing positive integers:
 $$
  \pi(n)\;:\; n\;=\; n_1+\,\cdots\,+ n_k\,; \hspace{1em}
  \pi(n_i)\;:\; n_i\; =\;  d_{i1}+\,\cdots\,+ d_{ik_i}\,,\;
                i=1\,,\,\ldots\,, k\,.
 $$
We will call this double partition the {\it type},
 in notation $\type(J)$, of $J$.
Denote also the set of all such double partitions of $n$ by
 $\PP(n)$.
Then the admissible permutations of the blocks $A_i, \, \cdots\,, A_k$
 induces a finite group action on $\PP(n)$.
The quotient set is denoted by $\PP(n)/\!\sim$.
For a general $m\in M_n({\Bbb C})$, define its type by
 $\type(m)=\type(J_m)$,
 which is uniquely defined after passing to $\PP(n)/\!\sim$.

\bigskip

\noindent
{\bf Definition 4.1.2 [partial order between Jordan forms].} {\it
 Given two Jordan forms $J_1$ and $J_2$, we say that
  $J_1\prec J_2$ if the following two conditions are satisfied:
  \begin{itemize}
   \item[$(1)$]
    $J_1$, $J_2$ have the same characteristic values
      $\lambda_1,\,\cdots\,,\lambda_k$ of
     the same multiplicities $n_i$ for $\lambda_i$.

   \item[$(2)$]
    Let
     $A_{1i},\, A_{2i}\in M_{n_i}({\Bbb C})$ be the diagonal blocks
      of $J_1$ and $J_2$ respectively that are associated to $\lambda_i$
       and
     ${\mathbf 1}_{n_i}$ be the identity of $M_{n_i}({\Bbb C})$.
    Then $\rank((A_{1i}-\lambda_i{\mathbf 1}_{n_i})^j)
          \le \rank ((A_{2i}-\lambda_i{\mathbf 1}_{n_i})^j)$
     for all $j\in {\Bbb N}$.
  \end{itemize}
} 

\bigskip

\noindent
This defines a partial order $\prec$ on the set of Jordan matrices
 in $M_n({\Bbb C})$ that is invariant under admissible permutations
 of diagonal blocks of distinct characteristic values.

\bigskip

\noindent
{\bf Proposition 4.1.3 [partial order of orbits via Jordan forms].}
([M-T], [Ge], [Dj].) {\it
  $$\mbox{
   $O_{\varphi_{m_1}}\prec O_{\varphi_{m_2}}$
   if and only if $J_{m_1}\prec J_{m_2}$. }
  $$
} 

\smallskip

The following simplified/coarser partial order helps us to see things
 more directly.

\bigskip

\noindent
{\bf Definition 4.1.4 [isotopic decay].}\footnote{For topologists:
                                        Here the term ``{\it isotopic}"
                                         comes from the notion of
                                         ``isotope" in physics/chemitry,
                                         not topology.
                                        The reason why we choose this term
                                         is partially enlightened
                                         in footnote 35.}
{\rm
 The composition of a sequence of operations
  of the form
  $J^{(\lambda)}_j \rightarrow
   \Diag(J^{(\lambda)}_{j_1}\,,\, J^{(\lambda)}_{j_2})$
  with $j=j_1+j_2$, $j_1\ge j_2$,
 will be called an {\it isotopic decay}.
} 

\bigskip

\noindent
Given two Jordan forms $J_1$ and $J_2$,
 define $J_1\pprec J_2$ if $J_1$ is obtained from $J_2$ by
 a sequence of isotopic decays and an re-arrangement of
 the sub-blocks in each diagonal block associated to
 a characteristic value.

\bigskip

\noindent
{\bf Lemma 4.1.5 [coarser partial order].} {\it
 $(1)$
 $O_{m_1}\prec O_{m_2}$ if $J_{m_1}\pprec J_{m_2}$.
 $(2)$
 $\prec$ and $\pprec$ generate the same equivalence relation,
  in notation $\approx$, on the set of Jordan forms.
} 

\bigskip

\noindent
{\it Remark 4.1.6 $[$orbit dimension drop under $\pprec$$]$.}
(E.g.\ [We]; also [Ge] or [Bas].)
 Let
  {\scriptsize
  $$
   \begin{array}{c}
    \\[1ex]
   T^{(b_1,\,\cdots\,b_i)}_{i\times j}\\[1em]
   (i\le j)
   \end{array}\;
   =\;
   \left[
    \begin{array}{cccccc}
      b_1    &         &         &      &         \\[.6ex]
      b_2    & b_1     &         &      &         \\
             & b_2     & \ddots  &      &     & \hspace{1em}0    \\
      \vdots & \ddots  & \ddots  & b_1  &         \\[.6ex]
      b_i    & \cdots  &         & b_2  & b_1 &
    \end{array}
   \right]_{i\times j}\,,
   \hspace{1em}
   \begin{array}{c}
    \\[1ex]
   T^{(b_1,\,\cdots\,b_j)}_{i\times j}\\[1em]
   (i\ge j)
   \end{array}\;
   =\;
   \left[
    \begin{array}{ccccc}
             &         &    0                     \\[1em]
      b_1    &         &         &      &         \\[.6ex]
      b_2    & b_1     &         &      &         \\
             & b_2     & \ddots  &      &         \\
      \vdots & \ddots  & \ddots  & b_1  &         \\[.6ex]
      b_j    & \cdots  &         & b_2  & b_1
    \end{array}
   \right]_{i\times j}\,.
  $$
 {\normalsize Here}}, all omitted entries are zero.
 The centralizer of $J$, in the form given previously, consists
  of all matrices of the form
  $$
   \left[
    \begin{array}{ccc}
     B_1 &         & 0 \\
         & \ddots  &   \\
     0   &         & B_k
    \end{array}
   \right]
  $$
  with each $B_i\in M_{n_i}({\Bbb C})$ of the block form
   $\left[ B_{i,rs}\right]_{k_i\times k_i}$
   where
   $$
    B_{i,rs}\; =\; T^{(b_{i,rs; 1}, \,\cdots\,b_{i,rs; di_r})}
                                                _{di_r\times di_s}
     \hspace{1ex}\mbox{for $r\ge s$},\hspace{1em}
    B_{i,rs}\; =\; T^{(b_{i,rs; 1}, \,\cdots\,b_{i,rs; di_s})}
                                               _{di_r\times di_s}
     \hspace{1ex}\mbox{for $r<s$.}
   $$
 (Again, omitted entries are all zero.)
 The dimension of the stabilizer of $J$, as given, is
  thus
  \begin{eqnarray*}
   \lefteqn{
    n\; \le \;
     \dimm_{\scriptsizeBbb C}\Stab(J) }\\[.6ex]
    && =\;\sum_{i=1}^k\, \left(\rule{0em}{1em}
          (d_{i1}+\,\cdots\,+ d_{ik_i}) + 2(d_{i2}+\,\cdots\,+ d_{ik_i})
            + \, \cdots\; + 2(d_{ik_i})
                      \right)\; \le\; n^2\,.
  \end{eqnarray*}
 Thus, for each
  $J^{(\lambda)}_j \rightarrow
    \Diag(J^{(\lambda)}_{j_1}\,,\,J^{(\lambda)}_{j_2})$ with $j_1\ge j_2$
  the corresponding new adjoint-orbit drops the dimension by
  an integral amount $\ge j_2$.

\bigskip

Some properties of $\Map((\Space M_n({\Bbb C});{\Bbb C}^n), {\Bbb A}^1)$
 are listed below:
\begin{itemize}
 \item[(1)]
 The equivalence relation $\approx$ in Lemma 4.1.5
  descends to an equivalence relation, still denoted by $\approx$,
  on the topological space
  $\Map((\Space M_n({\Bbb C});{\Bbb C}^n), {\Bbb A}^1)$.
  The associated quotient space
  $\Map((\Space M_n({\Bbb C});{\Bbb C}^n), {\Bbb A}^1)/\!\!\approx$
  is the $n$-th symmetric product
  $S^n{{\Bbb A}^n}
   :=({\Bbb A}^1)^n/\Sym_n\simeq\footnote{However, caution that
                                       under this isomorphism that
                                        comes from the ring generated
                                        by elementary symmetric polynomials,
                                       the diagonal locus
                                        in $({\footnotesizeBbb A}^1)^n$
                                        becomes a complicated
                                        discriminant locus
                                        in ${\footnotesizeBbb A}^n$.}\;
    {\Bbb A}^n$ of ${\Bbb A}^1$,
  where $\Sym_n$ is the permutation group of $n$ letters.
 Each $\approx$-equivalence class of points
  on $\Map((\Space M_n({\Bbb C});{\Bbb C}^n), {\Bbb A}^1)$
  contains a unique maximal point and a unique minimal point
  with respect to $\prec$ on
  $\Map((\Space M_n({\Bbb C});{\Bbb C}^n), {\Bbb A}^1)$.
 Any other point in the same class is sandwiched between the two
  by $\prec$.

 \item[(2)]
 The types of Jordan forms give rise to a finite stratification
  $\{S_t\}_t$ of\\ $\Map((\Space M_n({\Bbb C});{\Bbb C}^n), {\Bbb A}^1)$.
 The stratum associated to the double partition
  $$
   \pi(n)\;:\; n\;=\; n_1+\,\cdots\,+ n_k\,; \hspace{1em}
   \pi(n_i)\;:\; n_i\; =\;  d_{i1}+\,\cdots\,+ d_{ik_i}\,,\;
                 i=1\,,\,\ldots\,, k\,,
  $$
  of $n$
  is homeomorphic to
  $({\Bbb C}^k-\mbox{(diagonal locus)})/\Sym_k$.
 Here, `diagonal locus' means the set of all points
  whose coordinates have some identical entries.
 The stratum $S_{(n=1+\,\cdots\,+1)}$ is open dense
  in $\Map((\Space M_n({\Bbb C});{\Bbb C}^n), {\Bbb A}^1)$.
\end{itemize}
%

\bigskip

\begin{flushleft}
{\bf The Chan-Paton space/module on D0-branes on ${\Bbb A}^1$.}
\end{flushleft}
Let
 $m\in M_n({\Bbb C})$ with the Jordan form as given above  and
 $\langle {\mathbf 1}, m\rangle$ be the sub-algebra of
  $M_n({\Bbb C})$ generated by ${\mathbf 1}$ and $m$.
$\langle {\mathbf 1}, m\rangle$ is commutative.
The {\it c}haracteristic polynomial and the {\it min}imal polynomial
 of $m$ are then respectively
 $$
  f_m^c(\lambda)\;
   =\; (\lambda-\lambda_1)^{n_1}\,\cdots\,(\lambda-\lambda_k)^{n_k}
  \hspace{1em}\mbox{and}\hspace{1em}
  f_m^{\scriptsizemin}(\lambda)\;
   =\; (\lambda-\lambda_1)^{d_{11}}\,\cdots\,(\lambda-\lambda_k)^{d_{k1}}\,.
 $$

\bigskip

\noindent
{\bf Lemma 4.1.7 [interpolation formula].} {\it
 Given
  $g(\lambda)
    := (\lambda-\lambda_1)^{d_1}\,\cdots\,(\lambda-\lambda_k)^{d_k}
    \in {\Bbb C}[\lambda]$
  with the $\lambda_i$'s distinct from each other,
 then
  the inverse of $g(\lambda)/(\lambda-\lambda_i)^{d_i}$
  in ${\Bbb C}[\lambda]/((\lambda-\lambda_i)^{d_i})$ exists,
  for $i=1,\,\ldots\,, k$.
 Denote this inverse by $(1/g_{(i)})(\lambda)$,
  which is a polynomials of degree $\le d_i-1$.
 Let $d=d_1+\,\cdots\,+d_k$ and
  $f(\lambda)$ be a polynomial of degree $< d$.
 Then
  there exist unique polynomials $f_i(\lambda)$
   with $\degree f_i(\lambda) < d_i$
  such that
   $$
    f(\lambda)\; =\;
    \sum_{i=1}^k\:
      f_i(\lambda) \cdot (1/g_{(i)})(\lambda)
      \cdot \frac{g(\lambda)}{(\lambda-\lambda_i)^{d_i}}\,.
   $$
  Indeed, $f_i(\lambda)$ is the Taylor expansion of
   $f(\lambda)$ in $(\lambda-\lambda_i)$
   up to (including) degree $d_i-1$.
} 

\bigskip

\noindent
It follows that, as a ${\Bbb C}$-algebra,
 \begin{eqnarray*}
  \lefteqn{
   \langle{\mathbf 1},m \rangle\;
    \simeq\; {\Bbb C}[\lambda]/(f_m^{\scriptsizemin}(\lambda)) }\\[.6ex]
  &&\hspace{1.4em}
  =\; \sum_{i=1}^k
        \left(
         (1/{f_{m\; (i)}^{\scriptsizemin}})(\lambda)
          \cdot \frac{ f_m^{\scriptsizemin}(\lambda)}
                                  {(\lambda-\lambda_i)^{d_{i1}} }
         \right)\;
  \simeq\;
    \prod_{i=1}^k
      \left( {\Bbb C}[\lambda]/(\lambda-\lambda_i)^{d_{i1}}  \right)\,.
 \end{eqnarray*}
The sum in the above expression is a direct sum of
 orthogonal indecomposable ideals in\\
 ${\Bbb C}[\lambda]/(f_m^{\scriptsizemin}(\lambda))$
 associated to the decomposition
 $$
  1 \;
    =\; \sum_{i=1}^k
          (1/{f_{m\;(i)}^{\scriptsizemin}})(\lambda)
           \cdot \frac{ f_m^{\scriptsizemin}(\lambda)}
                                  {(\lambda-\lambda_i)^{d_{i1}} }
 $$
 through the complete set of primitive orthogonal idempotents
 in ${\Bbb C}[\lambda]/(f_m^{\scriptsizemin}(\lambda))$.
The length $l_{\langle{\mathbf 1},m\rangle}$
 of $\langle{\mathbf 1},m\rangle$ is
 $\degree f_m^{\scriptsizemin}(\lambda)=d_{11}+\,\cdots\,+d_{k1}$.

Let ${\Bbb C}^n$ be the unique non-zero irreducible representation of
 $M_n({\Bbb C})$.
Up to the $\GL_n({\Bbb C})$ adjoint action, we may assume that $m$ is
 already a Jordan form $J=\Diag(A_1,\,\cdots\,,A_k)$ given earlier.
Let ${\mathbf 1}_{(i)}
      = \Diag(0,\,\cdots\,, 0, {\mathbf 1}_{n_i},0, \,\cdots\,,0)$,
 where
  ${\mathbf 1}_{n_i}$ in the $i$-th position is the identity matrix
   $\in M_{n_i}({\Bbb C})$  and
  the $0$ in the $j$-th position are the zero-matrix
   $\in M_{n_j}({\Bbb C})$
   for $j=1,\,\cdots\,, i-1, i+1, \,\cdots\,, k$.
Then
 $$
  \left.
  \left(
   (1/{f_{J\;(i)}^{\scriptsizemin}})(\lambda)
        \cdot \frac{ f_J^{\scriptsizemin}(\lambda)}
                             {(\lambda-\lambda_i)^{d_{i1}} }
   \right)\right|_{\,\lambda\,=\,J}\; =\; {\mathbf 1}_{(i)}\,.
 $$
This implies that
 ${\mathbf 1}_{(i)}\in \langle{\mathbf 1}, J\rangle$
  for $i=1,\,\ldots\,,k$  and that
 ${\mathbf 1}={\mathbf 1}_{(1)}+\,\cdots\,+{\mathbf 1}_{(k)}$
  is an orthogonal primitive idempotent decomposition
  in $\langle{\mathbf 1},J\rangle$.
The corresponding direct-sum decomposition,
 now as $\langle{\mathbf 1}, J\rangle$-modules,
 $$
  {\Bbb C}^n\; =\; {\mathbf 1}_{(1)}\cdot {\Bbb C}^n \,+\,\cdots\,
                   +\, {\mathbf 1}_{(k)}\cdot {\Bbb C}^n\;
  =\; {\Bbb C}^{n_1}\,+\,\cdots\,+\,{\Bbb C}^{n_k}\;
  =:\; V_1\,+\,\cdots\,+V_k
 $$
 is the same decomposition of ${\Bbb C}^n$ that renders $J$
 the given diagonal block form.
As a $\langle{\mathbf 1},J\rangle$-module,
 $V_i$ $(={\Bbb C}^{n_i})$ decomposes into a direct sum
  $V_i = {\Bbb C}^{d_{i1}}+\,\cdots\,+{\Bbb C}^{d_{ik_i}}
       =: V_{i1}+\,\cdots\,+V_{ik_i}$
 of indecomposable $\langle{\mathbf 1},J\rangle$-modules.
$\Spec\langle{\mathbf 1},J\rangle$ has $k$-many connected components,
 associated respectively to ideals $({\mathbf 1}-{\mathbf 1}_{(i)})$
 in $\langle{\mathbf 1},J\rangle$,
 $i=1,\,\ldots\,,k$.
One has that
 $$
  \langle{\mathbf 1}, J\rangle/({\mathbf 1}-{\mathbf 1}_{(i)})\;
   =\;  \langle{\mathbf 1}, J\rangle\cdot{\mathbf 1}_{(i)}\;
   \simeq\; \langle{\mathbf 1}_{n_i}, A_i\rangle\;
   \simeq\;  {\Bbb C}[\lambda]/((\lambda-\lambda_i)^{n_i})
 $$
 and that the annihilator $\Ann(V_i)$ of $V_i$ $(={\Bbb C}^{n_i})$
  as an $\langle{\mathbf 1}, J\rangle$-module
  is $({\mathbf 1}-{\mathbf 1}_{(i)})$.
In terms of
 $\langle{\mathbf 1}, J\rangle
  \simeq \prod_{i=1}^k\langle{\mathbf 1}_{n_i}, A_i\rangle$,
 the $\langle{\mathbf 1},J\rangle$-modules
  $V_i$, $V_{i1}, \,\cdots\,, V_{ik_i}$
  are also $\langle{\mathbf 1}_{n_i}, A_i\rangle$-modules
  automatically.

The above algebraic statements correspond to
 the following geometric picture of Chan-Paton modules
 on the associated D0-branes on ${\Bbb A}^1\,$:
\begin{itemize}
 \item[(1)]
 Under Grothendieck Ansatz or Lemma 1.2.19,
  $\varphi_J: {\Bbb C}[y]\rightarrow M_n({\Bbb C})$
  gives(/is equivalent to) a morphism
  $\hat{\varphi}_J:\Space M_n({\Bbb C})\rightarrow {\Bbb A}^1$
  with the image subscheme
  $\image\hat{\varphi}_J \simeq \Spec\langle{\mathbf 1},J\rangle$
  associated to the ideal
  $$
   \Ker(\varphi_J)\;
   =\; (f_m^{\scriptsizemin}(y))\;
   =\; \left(
        (y-\lambda_1)^{d_{i1}}\,\cdots\,(y-\lambda_k)^{d_{k1}}
       \right)
  $$
  in ${\Bbb C}[y]$.
 Thus, on ${\Bbb A}^1$ there are $k$-many (generally non-reduced)
  points located respectively at $y=\lambda_1,\,\cdots\,,\lambda_k$
  (in the underlying complex plane ${\Bbb C}$ of ${\Bbb A}^1$)
  where D0-branes in Polchinski's sense may sit upon.
 These are the {\it D0-branes on ${\Bbb A}^1$ associated to
  $\varphi_J$} in the sense of Definition 2.2.3.
 From the discussion, for a general $\varphi_m$, they depend only on
  the minimal polynomial $f_m^{\scriptsizemin}(\lambda)$ of $m$.

 \item[(2)]
 The push-forward\footnote{For non-algebro-geometers:
                             ${\footnotesizeBbb C}^n$
                             as a $\langle{\mathbf 1},J\rangle$-module
                             is now a ${\footnotesizeBbb C}[y]$-module
                             via $\varphi_J$, with annihilator
                             $\Ker(\varphi_J)$.
                            Thus, though $\hat{\varphi}_J$
                             is not directly defined,
                             $\hat{\varphi}_{J\ast}{\footnotesizeBbb C}^n$
                             is well-defined.
                            This is the Grothendieck Ansatz
                             on quasi-coherent modules versus
                             quasi-coherent sheaves,
                             similar to that on rings versus spaces.}
  $\hat{\varphi}_{J\ast}{\Bbb C}^n
   =\sum_{i=1}^k\hat{\varphi}_{J\ast}V_i
   =\sum_{i=1}^k\sum_{j=1}^{k_i}\hat{\varphi}_{J\ast}V_{ij}$
  is now an ${\cal O}_{\scriptsizeimage\hat{\varphi}_J}$-module
  of length $n$.
 Decompose $\image\hat{\varphi}_J$ into a disjoint union
  $\amalg_{i=1}^kZ_i$, where $Z_i$ is the subscheme of ${\Bbb A}^1$
  associated to the ideal $((y-\lambda_i)^{d_{i1}})$.
 Then $\hat{\varphi}_{J\ast}V_i$ is supported on $Z_i$ and, hence,
  is an ${\cal O}_{Z_i}$-module of length $n_i$.
 The decomposition
  $\hat{\varphi}_{J\ast}V_i
    =\sum_{j=1}^{k_i}\hat{\varphi}_{J\ast}V_{ij}$
  is automatically a direct-sum decomposition
  as ${\cal O}_{Z_i}$-modules as well.
 Let $Z_i^{(l)}$, $l\le n_i$, be the subscheme of $Z_i$
  associated to the ideal $((y-\lambda_i)^l)$ in ${\Bbb C}[y]$.
 Note that $Z_i^{(l)}$ has length $l$  and that
  $Z_i^{(1)}$ is the ${\Bbb C}$-point in $Z_i$ and $Z_i^{(n_i)}=Z_i$.
 Then,
  $\hat{\varphi}_{J\ast}V_{ij}$ is a rank-$1$ ${\cal O}_{Z_i}$-module
   of length $d_{ij}$ and is supported on $Z_i^{(d_{ij})}$.
 As ${\cal O}_{Z_i}$-modules,
  $$
   \hat{\varphi}_{J\ast}V_{i1}\; \simeq\; {\cal O}_{Z_i}
  $$
  and
  $$
   \mbox{$\hat{\varphi}_{J\ast}V_{ij}\;
      \simeq\;$
       the ideal $(y-\lambda_i)^{d_{i1}-d_{ij}}\cdot{\cal O}_{Z_i}$
        of ${\cal O}_{Z_i}\;
      \simeq\;$
        the quotient ${\cal O}_{Z_i^{(d_{ij})}}$ of ${\cal O}_{Z_i}$}\,.
  $$
 In our setting\footnote{See
                           footnote 35 for remarks on the original
                           setting in string theory.},
  we call $\hat{\varphi}_{J\ast}V_i$
  the {\it Chan-Paton module on the D0-brane supported on
  $Z_i\subset {\Bbb A}^1$ associated to $\varphi_J$}.
 From the discussion, for a general $\varphi_m$, their isomorphism class
  depends only on both $f_m^{\scriptsizemin}(\lambda)$ and the type of $m$.
\end{itemize}

\bigskip

\begin{flushleft}
{\bf Comparison with Hilbert schemes and Chow varieties.}
\end{flushleft}
The Hilbert scheme $\Hilb^n_{{\scriptsizeBbb A}^1}=: ({\Bbb A}^1)^{[n]}$
 of $n$ points on ${\Bbb A}^1$ parameterizes $0$-dimensional subschemes
 of length $n$ on ${\Bbb A}^1$.
Such a subscheme of ${\Bbb A}^1$ is given uniquely by an ideal
 $(f)\subset {\Bbb C}[y]$, where $f$ is a monic polynomial of degree $n$.
In terms of matrices, it is thus represented by an $m\in M_n({\Bbb C})$
 such that both the characteristic polynomial and the minimal polynomial
 of $m$ are $f$.\footnote{In other words,
                                   $m$ is a regular matrix in
                                   $M_n({\footnotesizeBbb C})$.}
Observe that the Jordan form of
 (omitted entries are zero; the multiplicity of $\lambda_i=n_i$)

\bigskip

{\tiny
 $$
  \mbox{\normalsize
   $J^{(\lambda_1,\,\cdots\,,\lambda_k)}_+\; :=\;$}
  \left[
  \begin{array}{ccccccccccccc}
   \lambda_1 & \\
    1   & \ddots \\
        & \ddots &\ddots \\
    &   &     1  & \lambda_1   \\[3ex]
    & & & {\mathbf 1}    & \lambda_2 \\
    &&& & 1   & \ddots \\
    &&&&& \ddots  & \ddots \\
    &&&&&& 1    &  \lambda_2 \\
    &&&&&&&  \ddots & \ddots \\[1ex]
    &&&&&&&&  {\mathbf 1}     & \lambda_k \\
    &&&&&&&&&  1  & \ddots  \\
    &&&&&&&&&& \ddots  & \ddots \\
    &&&&&&&&&&&  1     & \lambda_k
  \end{array}
  \right]\,,
 $$} 

 \bigskip

 \noindent
 where
 $(\lambda_1,\,\cdots\,,\,\lambda_k)$ is the $n$-tuple
  $(\lambda_1,\,\cdots\,,\,\lambda_1,\,
    \lambda_2,\,\cdots\,,\,\lambda_2,\,
    \cdots\,,\,
    \lambda_k,\,\cdots\,,\,\lambda_k)$
  with the specified multiplicity $n_i$ for $\lambda_i$,
 is $\Diag(J^{(\lambda_1)}_{n_1}\,,\,\cdots\,,\,J^{(\lambda_k)}_{n_k})$,
 up to a permutation of the blocks.
Its characteristic polynomial and minimal polynomial are identical:
 $(y-\lambda_1)^{n_1}\,\cdots\,(y-\lambda_k)^{n_k}$.
Let ${\Bbb C}^n$ parameterizes the ordered tuples of roots of
 monic polynomial of degree $n$, then the embedding
 $$
  {\Bbb C}^n\; \hookrightarrow\; M_n({\Bbb C})\,,
   \hspace{1em}
  (\lambda_1,\,\cdots\,,\,\lambda_n)\;
   \mapsto\; J^{(\lambda_1,\,\cdots\,,\,\lambda_n)}_+
 $$
 descends to an embedding
 $$
  \Phi_{\scriptsizeHilb}: ({\Bbb A}^1)^{[n]}\;
    \longrightarrow\;
    \Map((\Space M_n({\Bbb C});{\Bbb C}^n), {\Bbb A}^1)\,,
  \hspace{1em}
  \mbox{$\prod_{i=1}^n(y-\lambda_i)$}\;
    \longmapsto\; \varphi_{J^{(\lambda_1,\,\cdots\,,\,\lambda_n)}_+}\,.
 $$

On the other hand, the Chow variety
 $\Chow^{(n)}_{{0, \scriptsizeBbb A}^1}$ of $n$ points on ${\Bbb A}^1$
 parameterizes $0$-cycles of order $n$ on ${\Bbb A}^1$ and
 is identical to the $n$-th symmetric product $S^n({\Bbb A}^1)$
 of ${\Bbb A}^1$.
Such a $0$-cycle on ${\Bbb A}^1$ happens to be represented uniquely
 by a monic polynomial in $y$ of degree $n$ as well.
Thus there is a canonical isomorphism
 $({\Bbb A}^1)^{[n]}\simeq S^n({\Bbb A}^1)$.
However, from the general ground of Chow groups, the support of
 a cycle is meant to be a reduced subscheme with each of
 its irreducible components marked with a multiplicity.
Thus, in terms of matrices, it is represented by an $m\in M_n({\Bbb C})$
 such that the minimal polynomial of $m$ has only simple roots.
Such matrices are exactly the diagonalizable matrices. Again, let
${\Bbb C}^n$ parameterizes the ordered tuples of roots of
 monic polynomial of degree $n$, then it follows that
 the embedding
 $$
  {\Bbb C}^n\; \hookrightarrow\; M_n({\Bbb C})\,,
   \hspace{1em}
  (\lambda_1,\,\cdots\,,\,\lambda_n)\;
   \mapsto\; \Diag(\lambda_1,\,\cdots\,,\,\lambda_n)
 $$
 descends to an embedding
 $$
  \Phi_{\scriptsizeChow}: S^n({\Bbb A}^1)\;
    \longrightarrow\;
    \Map((\Space M_n({\Bbb C});{\Bbb C}^n), {\Bbb A}^1)\,,
  \hspace{1em}
  \mbox{$\prod_{i=1}^n(y-\lambda_i)$}\; \longmapsto\;
    \varphi_{\scriptsizeDiag(\lambda_1,\,\cdots\,,\,\lambda_n)}\,.
 $$
In other words,
 $\image\Phi_{\scriptsizeHilb}$ parameterizes conjugacy classes of
  regular representations of ${\Bbb C}[y]$ in $M_n({\Bbb C})$
 while $\image\Phi_{\scriptsizeChow}$ parameterizes conjugacy classes
  of diagonal representations of ${\Bbb C}[y]$ in $M_n({\Bbb C})$.

Note that,
 under the isomorphism $({\Bbb A}^1)^{[n]}\simeq S^n({\Bbb A}^1)$,
 $\Phi_{\scriptsizeHilb}$ and $\Phi_{\scriptsizeChow}$
 coincide only on the open dense subset, points of which correspond
 to $0$-dimensional reduced subschemes of length $n$ on ${\Bbb A}^1$.
For all $p$ in the complement of this subset,
 $\Phi_{\scriptsizeChow}(p) \prec \Phi_{\scriptsizeHilb}(p)$
 by an isotopic decay.
In particular, $\Map((\Space M_n({\Bbb C});{\Bbb C}^n), {\Bbb A}^1)$
 contains $({\Bbb A}^1)^{[n]}$ and $S^n({\Bbb A}^1)$ distinctly and,
 for $n\ge 3$, has more points than
 $\Phi_{\scriptsizeHilb}(({\Bbb A}^1)^{[n]})
  \cup \Phi_{\scriptsizeChow}(S^n({\Bbb A}^1))$.
In the current case, it happens that
 $\Phi_{\scriptsizeHilb}$ and $\Phi_{\scriptsizeChow}$
 give rise to
 $$
  ({\Bbb A}^1)^{[n]}\; \stackrel{\sim}{\longrightarrow}\;
  \Map((\Space M_n({\Bbb C});{\Bbb C}^n), {\Bbb A}^1)/\!\approx \;\;
  \stackrel{\sim}{\longleftarrow}\; S^n({\Bbb A}^1)\,.
 $$
This is only accidental and does not generalize to $Y$ of
 dimension $\ge 2$.

Note also that, for all $p$,
 the Chan-Paton module at $\Psi_{\scriptsizeHilb}(p)$ gives exactly
  the structure sheaf ${\cal O}_{Z_p}$ of the subscheme $Z_p$
  $p$ represents
 while the Chan-Paton module at $\Psi_{\scriptsizeChow}(p)$ gives
  an association of ${\Bbb C}^{n_i}$ to each $p_i$
  (as an ${\cal O}_{p_i} (={\Bbb C})$-module),
  for $p=\sum_{i=1}^kp_i$ as a $0$-cycle.
Thus, Chan-Paton spaces/modules in the sense of Definition 2.2.3
 tells the difference of subschemes versus cycles as well.\footnote{Some
                  stringy comments follow.
                   When generalized to higher-dimensional D-branes,
                   these notions produce different notions of
                   ``wrappings" of a D-brane around
                   a submanifold/subvariety in the target space(-time)
                   of strings.
                  Such a subtlety, among other things,
                   was recognized seriously only by a smaller group
                   of string theorists, e.g.\ [G-S] and [H-S-T].
                  For most of the stringy literatures, the simpler
                   cycle-picture are more dominating (in the region of
                   the related Wilson's theory-space where
                   ``branes are really branes").
                  In the hind sight, there might be a reason for this:
                   Recall that an open string interacts with
                   D-branes via its end-points.
                   In most disscussions/literatures, these end-points
                   are only taken to be simple points (i.e.\
                   reduced points in the algebro-geometric language)
                   and hence, despite the fact that D-brane warpping can
                   be a more complicated notion than usually thought of,
                   open strings do not see anything beyond
                   the cycle picture with a gauge bundle supported thereon.
                  Should one remember that an end-point is attached to
                    the open string and there are jets
                    (in the sense of differential topology or,
                     in the open-string world-sheet picture,
                     in the sense of real algebraic geometry)
                    at the end-point,
                   then one may expect to draw out
                    some open-string-parameterization-invariant details of
                    such hidden ``thickened structure"
                    (e.g.\ non-reducedness of subschemes, embedded points,
                           torsion-subsheaves within a torsion sheaf,
                           ..., etc.).
                  (However, except in the elementary discussion of
                   momentum conservation of open strings, in which $1$-jet
                   is involved, we are not aware of any other use of jets
                   at the end-point of open string in string theory.)

                  On the other hand, since a D-brane
                    (again in the ``brane is really a brane" region)
                    is now taken as an extended dynamical object
                    in its own right and hence has its own definition
                    and deformation-obstruction theory,
                   while it must contain contents induced from open strings,
                    it is completely legitimate that
                    it could also have contents without contradictions
                    with open strings and yet open strings cannot see.
                  In the current example and in Polchinski's picture,
                   open strings do not ``see" the (more complicated)
                   structure sheaf ${\cal O}_{Z_p}$ when it is non-reduced
                   but, rather, only see the (simpler) cycle
                   $\Sigma_{i=1}^kn_ip_i$ with the sub-Chan-Paton space
                   ${\footnotesizeBbb C}^{n_i}$ attached to each $p_i$.
                  Furthermore, what open strings do not see is
                   nevertheless transformable to what open strings
                   do see via an isotopic decay.
                  This is actually a general feature.
                  In our setting,
                   we take both as different yet allowable
                   existences of D-branes on the target space(-time)
                   from deformations of D-branes in the sense
                   of deformations of morphisms from an Azumaya-type
                   noncommutative space to the open-string
                   target space(-time).
                  This explains also the term in Definition 4.1.4,
                   cf.\ footnote 30.}    
This is a general feature.

Finally, the map that sends $\varphi_m$ to the diagonal of $J_m$
 gives rise to a continuous map
 $\pi_{\scriptsizeChow}:
  \Map((\Space M_n({\Bbb C});{\Bbb C}^n),{\Bbb A}^1)
                                     \rightarrow S^n({\Bbb A}^1)$.
It has $\Phi_{\scriptsizeChow}$ as a section.

\bigskip

\begin{flushleft}
{\bf Associated quiver.}
\end{flushleft}
Given a finite-dimensional ${\Bbb C}$-algebra $R$, one can associate
 a {\it quiver}\footnote{A
                         few definitions/remarks for readers' reference
                         are put here to make precise of the discussion
                         while avoiding distractions.
                        A `{\it quiver}' is an oriented graph $\Gamma$
                         introduced in, e.g., the work of Gabriel in early
                         1970s to study representations of algebras.
                        A {\it representation of a quiver} $\Gamma$
                         over ${\footnotesizeBbb C}$ is an assignment
                         to each vertex $v_i\in\Gamma$
                          a ${\footnotesizeBbb C}$-vector space $V_i$ and
                         to each arrow (i.e.\ oriented edge) $\in \Gamma$
                           from $v_i$ to $v_j$
                          a ${\footnotesizeBbb C}$-linear homomorphism
                          $\varphi_{ij}:V_i\rightarrow V_j$.
                        Such representations have now become also
                         a standard tool for string theorists
                        to encode the field contents in a supersymmetric
                         gauge field theory coupled with matters.
                        Such field theories occur particularly
                         on (the world-volume of) D-branes.
                        Due to the rigidity of supersymmetric field theory,
                         a quiver representation pretty much fixes
                         the combinatorial type of the field theory
                         under investigation.

                        There are {\it different} quivers that can be
                         associated to a finite-dimensional
                         ${\footnotesizeBbb C}$-algebra $R$,
                         regarded as a (left) $R$-module from
                         the algebra multiplication.
                        The one we choose here encodes
                         the embedded dimension
                         (i.e.\ the dimension of the tangent space
                                when re-phrased in geometry) of
                         of the Artinian ${\footnotesizeBbb C}$-algebra
                         in our problem.
                        See, e.g., [A-R-S], [G-R], and [Jat]
                         for more discussions.}
 $\Gamma_R$ to $R$ as follows:
 \begin{itemize}
  \item[(1)]
   Let $\{e_1,\,\cdots\,,e_k\}$ be a complete set of
    primitive orthogonal idempotents in $R$.
   Then associate to each $e_i$ a vertex, denoted also by $e_i$.

  \item[(2)]
   Let $J(R)$ be the radical of $R$.
   Then, associate $\dimm_{\scriptsizeBbb C}\,e_i(J(R)/J(A)^2)e_j$-many
    arrows from $e_i$ to $e_j$.
 \end{itemize}

Applying this to $\varphi_m$, representing a point
 in $\Map((\Space M_n({\Bbb C});{\Bbb C}^n)\,,\,{\Bbb A}^1)$,
 by associating a graph to the Artinian ${\Bbb C}$-algebra
 ${\Bbb C}[y]/\Ker\varphi_m\simeq \langle{\mathbf 1},m\rangle$,
 following the rules above,
we obtain a quiver $\Gamma_{\varphi_m}$ that captures part of
 the geometry of the D0-brane on ${\Bbb A}^1$ associated to
 $\varphi_m$:
 \begin{itemize}
  \item[$\cdot$]
   a {\it vertex} $e_i$ for the connected component $Z_i$ of
   $\Spec {\Bbb C}[y]/\Ker\varphi_m
    = \image\hat{\varphi}_m=\amalg_{i=1}^k Z_i$
   of the D0-brane on ${\Bbb A}^1$;

  \item[$\cdot$]
   an {\it arrow} with both ends attached to $e_i$ if $Z_i$
   has the embedded dimension $1$
   (i.e.\ if $Z_i$ is a non-reduced point on ${\Bbb A}^1$);
   there are no other arrows for any pair $(e_i,e_j)$, $1\le i,j\le k$.
 \end{itemize}
The Chan-Paton module discussed in an earlier theme is realized now
 as a representation of $\Gamma_{\varphi_m}$:
  (without loss of generality, we take $m$ to be the Jordan form
   $J=J_m$ and adopt earlier notations)
 \begin{itemize}
  \item[$\cdot$]
   assign the ${\cal O}_{Z_i}$-module
    $(\hat{\varphi}_{m\ast}{\Bbb C}^n)|_{Z_i}=\hat{\varphi}_{m\ast}V_i$
    to vertex $e_i$ for $i=1,\,\ldots\,,k$;

  \item[$\cdot$]
   if there is an arrow on $e_i$, then
    assign to that arrow the nilpotent endomorphism on $\hat{\varphi}_mV_i$
    associated to the multiplication by $(y-\lambda_i)\,$
   (i.e.\ the push-forward of
          the endomorphism $A_i-\lambda_i{\mathbf 1}_{n_i}$ on $V_i$).
 \end{itemize}
The quiver $\Gamma_{\varphi_m}$ together with this representation
 now encodes the full geometry of the connected components of
 the D0-brane on ${\Bbb A}^1$ except their exact locations
 $y=\lambda_1,\,\cdots\,,\lambda_k$.
 %

\bigskip

\begin{flushleft}
{\bf Higgsing/un-Higgsing of D-branes
     via deformations of morphisms.}\footnote{Readers who already know
                         the stringy side of Polchinski's D-branes are
                         suggested to compare it with the mathematical
                         picture described in this theme.
                        The Higgsing/un-Higgsing phenomenon described
                         in this theme following Definition 2.2.3
                         is a general feature.}
\end{flushleft}
The important open-string-induced
 Higgsing (i.e.\ gauge symmetry-breaking)/un-Higgsing
 (i.e.\ gauge symmetry enhancement) behavior on D-branes
 can be reproduced in the current content as follows.
As any associative ${\Bbb C}$-algebra $R$ gives rise to a Lie algebra
 $(R, [\cdot,\cdot])$ over ${\Bbb C}$ by taking the Lie bracket
 to be  $[m_1,m_2]=m_1m_2-m_2m_1$,
 we can equivalently make the discussion directly for
 associative algebras in our problem.

Since on $\Space M_n({\Bbb C})$, $M_n({\Bbb C})$ acts on
 the Chan-Paton space ${\Bbb C}^n$ as the endomorphism algebra
 $\End({\Bbb C}^n)$ of the Chan-Paton space,
this is the counterpart of (the Lie algebra of) the gauge symmetry
 on a D-brane in physicists' picture.
Given a
  $[\varphi_m:{\Bbb C}[y]\rightarrow M_n({\Bbb C})] \in
    \Map((\Space M_n({\Bbb C});{\Bbb C}^n),{\Bbb A}^1)$,
 the Chan-Paton space ${\Bbb C}^n$ on $\Space M_n({\Bbb C})$
  is turned into the Chan-Paton module on $\image \hat{\varphi}_m$
  by taking ${\Bbb C}^n$ now as a (left)
  $\langle{\mathbf 1},m\rangle$-module, as discussed earlier.
To distinguish them, we will denote the latter by
 $_{\langle{\mathbf 1},m\rangle}{\Bbb C}^n$.
Let
 $$
  \Centralizer\langle{\mathbf 1}, m\rangle\;
  :=\; \{m^{\prime\prime}\in M_n({\Bbb C})\,:\,
           m^{\prime\prime}m^{\prime}=m^{\prime}m^{\prime\prime}\;\;
           \mbox{for all $m^{\prime}\in \langle{\mathbf 1},m\rangle$} \}
 $$
 be the centralizer of $\langle{\mathbf 1},m\rangle$ in $M_n({\Bbb C})$.
Then,

\bigskip

\noindent
{\bf Lemma 4.1.8 [centralizer vs.\ pushed-forward endomorphism].} {\it 
 A ${\Bbb C}$-vector-space endomorphism
   $m^{\prime\prime}\in M_n({\Bbb C})$ of ${\Bbb C}^n$
  can be pushed forward to
   a $\langle{\mathbf 1}, m \rangle$-module endomorphism
    on $_{\langle{\mathbf 1},m\rangle}{\Bbb C}^n$
 if and only if
  $m^{\prime\prime}\in\Centralizer\langle{\mathbf 1},m\rangle
   \subset M_n({\Bbb C})$.
} 

\bigskip

\noindent
This gives a correspondence:
 \begin{itemize}
  \item[] {\it
   $\Centralizer\langle{\mathbf 1}, m\rangle\subset M_n({\Bbb C})$
    $\hspace{1em}\Longleftrightarrow\hspace{1em}$
   gauge symmetry on the D0-brane $\image\hat{\varphi}_m$
    on ${\Bbb A}^1$. }
 \end{itemize}
Recall further from earlier discussions
 the connected-component-decomposition
  $\image\hat{\varphi}_m=:Z=\amalg_{i=1}^kZ_i$  and
 the $\langle{\mathbf 1},m\rangle$-module direct-sum decomposition
  $_{\langle{\mathbf 1},m\rangle}{\Bbb C}^n=\sum_{i=1}^kV_i$
 with $\hat{\varphi}_{m\ast}V_i$ supported on $Z_i$.
Then, there is a natural direct-product decomposition
 as ${\Bbb C}$-algebras:
 $$
  \Centralizer\langle{\mathbf 1},m\rangle\;
   =\; \prod_{i=1}^k \Centralizer\langle{\mathbf 1},m\rangle_{(i)}\;
       \subset\; \prod_{i=1}^k \End(V_i)\;
       \simeq\;  \prod_{i=1}^kM_{n_i}({\Bbb C})\,.
 $$
Up to conjugation, we may assume that $m=J_m=J$ a Jordan form, then
 $\Centralizer\langle{\mathbf 1},m\rangle_{(i)}
   \subset M_{n_i}({\Bbb C})$
  consists of $n_i\times n_i$-matrices is of the form $B_i$
  given in Remark 4.1.6.
Thus, each $Z_i$ can be regarded as a D0-brane on ${\Bbb A}^1$
 in its own right, associated to
 $[\varphi_{B_i}]
   \in \Map((\Space M_{n_i}({\Bbb C});{\Bbb C}^{n_i}),{\Bbb A}^1)$,
 with the Chan-Paton module $\hat{\varphi}_{B_i\ast}{\Bbb C}^{n_i}$ and
 the gauge symmetry associated to the endomorphism subalgebra
 $\Centralizer\langle{\mathbf 1}_{n_i},B_i\rangle$ in $M_{n_i}({\Bbb C})$.
When $\varphi_m$ varies, this gives rise to Higgsing/un-Higgsing
 of gauge symmetry of D0-branes on ${\Bbb A}^1$.

%
%

In particular, if we restrict $\varphi_m$ to vary in
 $\Phi_{\scriptsizeChow}(S^n({\Bbb A}^1))
   \subset \Map((\Space M_n({\Bbb C}); {\Bbb C}^n),{\Bbb A}^1)$,
 then the Higgsing/un-Higgsing pattern of $n$ D0-branes on ${\Bbb A}^1$
  is as follows:
 \begin{itemize}
  \item[(1)]
   For $\varphi_m$ in the stratum associated to the type
    $(n=d_1+\,\cdots\,+d_k)$:
   \begin{itemize}
    \item[$\cdot$] [{\it D-branes on ${\Bbb A}^1$}]\\[.6ex]
     $Z=\amalg_{i=1}^kZ_i\simeq\amalg_{i=1}^k\Spec{\Bbb C}$,
      (i.e.\ $k$-collection of stacked D0-branes on ${\Bbb A}^1$);

    \medskip

    \item[$\cdot$] [{\it the Chan-Paton space}]\\[.6ex]
     ${\Bbb C}^{n_i}$ supported at the D0-brane
     $Z_i$ on ${\Bbb A}^1$ for $i=1,\,\ldots\,k$;

    \medskip

    \item[$\cdot$] [{\it gauge symmetry}]\\[.6ex]
      a factor $M_{n_i}({\Bbb C})\simeq \End({\Bbb C}^{n_i})$
       on $Z_i$ for $i=1,\,\cdots\,,k$;
      the total gauge symmetry of the $k$-many D0-brane system
       is the Lie algebra associated to
       the product $\prod_{i=1}^kM_{n_i}({\Bbb C})$.
   \end{itemize}

  \item[(2)]
   As a consequence of Item (1) above, when we vary
    $[\varphi_m]\in\Phi_{\scriptsizeChow}(S^n({\Bbb A}^1))$
    so that, for example,
   \begin{itemize}
    \item[$\cdot$] [{\it Higgsing}]\\[.6ex]
     $Z_1$ splits to $j$-many separated D0-brane collections
      $Z^{\prime}_1$, $\cdots\,$, $Z^{\prime}_j$ on ${\Bbb A}^1$,
      governed by the partition
      $n_1=n^{\prime}_1+\,\cdots\,+n^{\prime}_j$.
     Then the Chan-Paton space ${\Bbb C}^{n_1}$ splits as well and
      turns into a Chan-Paton-space ${\Bbb C}^{n^{\prime}_i}$
      at $Z^{\prime}_i$ for $i=1,\,\ldots\,,j$.
     The gauge symmetry associated to $M_{n_1}({\Bbb C})$ is now broken
      to the one associated to the sub-endomorphism-algebra
      $\prod_{i=1}^j M_{n^{\prime}_i}({\Bbb C})$
      with the factor $M_{n^{\prime}_i}({\Bbb C})$ assigned to
      $(Z^{\prime}_i,{\Bbb C}^{n^{\prime}}_i)$ for $i=1,\,\ldots\,,j$.

    \medskip

    \item[$\cdot$] [{\it un-Higgsing}]\\[.6ex]
     $Z_1,\,\cdots\,, Z_j$ collide/merge to a new $Z^{\prime}_j$.
     Then there is now a D0-brane collection at $Z^{\prime}_j$
      with Chan-Paton space ${\Bbb C}^{n_1+\,\cdots\,+n_j}$.
     The original gauge symmetry for the collection
      $\{(Z_1,{\Bbb C}^{n_1}) \,\cdots\,,(Z_j,{\Bbb C}^{n_j})\}$,
      which is the one associated to
      $M_{n_1}({\Bbb C})\times\,\cdots\,\times M_{n_j}({\Bbb C})$,
      is now enhanced to the gauge symmetry associated to
      $M_{n_1+\,\cdots\,+n_j}({\Bbb C})$, acting on
      $(Z^{\prime}_j,{\Bbb C}^{n_1+\,\cdots\,+n_j})$.
   \end{itemize}
 \end{itemize}
Except that we have to use algebraic groups
 -- in particular the $\GL_{\,\bullet}({\Bbb C})$-series
   in the current content -- in the pure algebro-geometric setting,
this is exactly the pattern of the oriented-open-string-induced
 Higgsing/un-Higgsing of unitary gauge symmetry of D-branes that
 Polchinski concluded in [Pol3: Sec.~3.3 and Sec.~3.4]\footnote{E-print
                                        version: {\tt hep-th/9611050}:
                                        Sec.~2.3 and Sec.~2.4.}.
In summary:

\bigskip

\noindent
{\bf Proposition 4.1.9
     [Higgsing/un-Higgsing of D0-branes on ${\Bbb A}^1$].}\footnote{For
                           non-string-theorists: On the physics side,
                           the Higgsing of gauge symmetry on D-branes
                           in the sense of Polchinski is originated
                           from the induced stretching of open strings
                           whose end-points are attached to D-branes
                           that are originally stacked
                           and then are deformed and separated.
                          Such stretching turns part of the massless
                           spectrum of open strings that contribute to
                           the gauge fields on the D-branes into
                           massive spectrum and hence reduces
                           the gauge fields on the D-branes.
                          The fact that this crucial open-string-induced
                           behavior of D-branes can be reproduced
                           by following Definition 2.2.3
                           alone {\it without} resorting to open strings
                           is what convince us that it makes sense
                           to take Definition 2.2.3
                           as the prototype intrinsic mathematical
                           definition for Polchinski's D-branes.
                          Unfamiliar readers are encouraged to
                           study [P-S] and [Pol4] to get a feeling.}
{\it
 The pattern of open-string-induced Higgsing/un-Higgsing behavior of
  $n$ D0-branes on ${\Bbb A}^1$ can be reproduced
  in the current content via deformations of morphisms
  $[\varphi_m:{\Bbb C}[y]\rightarrow M_n({\Bbb C})]$
  in $\Phi_{\scriptsizeChow}(S^n({\Bbb A}^1))
      \subset \Map((\Space M_n({\Bbb C}); {\Bbb C}^n),{\Bbb A}^1)$.
} 

\bigskip

\begin{flushleft}
{\bf Comparison with the spectral cover construction and
     the Hitchin system.}
\end{flushleft}
Fix a complex line bundle $\pi_L: L\rightarrow \pt\,$ over a point $\pt$.
We will identify $\pt$ with the zero-section of $L$ whenever needed.
Let $\lambda$ be the tautological section of $\pi_L^{\ast}L$ over $L$.

\bigskip

\noindent
{\bf Definition 4.1.10 [semi-simple pair].}\footnote{The
                               adjoint action of
                               $\GL_n({\footnotesizeBbb C})$
                               on $M_n({\footnotesizeBbb C})$ does not
                               have stable points in the sense of Mumford
                               in [M-F-K].
                              With Polchinski's D-branes in mind, we
                               choose semi-simple pairs for the role
                               of stable pairs in [Hi].}
{\rm
 A pair $(E,\phi)$,
  where
   $\pi_E:E\rightarrow \pt\,$
    is a rank-$n$ complex vector bundle over $\pt$  and
   $\phi:E\rightarrow E\otimes L$
    a complex-vector-bundle-homomorphism over $\pt\,$
  is called {\it semi-simple}
 if $\phi$ is semi-simple (i.e.\ diagonalizable)
  with respect to a (hence any) trivialization
  $E\simeq {\Bbb C}^n$ and $L\simeq {\Bbb C}$.
} 

\bigskip

\noindent
Associated to a semi-simple pair $(E,\phi)$, with the $\phi$ of type
 $(n=n_1+\,\cdots\,+n_k)$, are the following objects:
 \begin{itemize}
  \item[(1)]
   the reduced zero-locus $Z_{\phi}=\amalg_{i=1}^kZ_{\phi;i}$
    of the section $\det(\pi_L^{\ast}\phi-{\mathbf 1}\otimes\lambda)$
    of $\det(\pi_L^{\ast}E)\otimes (\pi_L^{\ast}L)^{\otimes n}$;

  \item[(2)]
   a direct-sum decomposition $E=\sum_{i=1}^k V_i$
    of bundles over $\pt\,$ so that\\
    $\hat{V}_i
     := (\pi_L^{\ast}V_i)|_{Z_{\phi,i}}
      = (\Ker(\pi_L^{\ast}\phi-{\mathbf 1}\otimes \lambda))|
                                                         _{Z_{\phi;i}}$
    for $i=1,\,\ldots\,,k$;

  \item[(3)]
   $\prod_{i=1}^k\End(V_i)\subset \End(E)\simeq M_n({\Bbb C})$
    acting on $E$ leaving each $V_i$ invariant for $i=1,\,\ldots\,,k$.
 \end{itemize}
This is the $0$-dimensional spectral cover construction
 in the sense of [Hi]; see also [B-N-R], [Don1], and [Ox].
The Hitchin system in this content takes the form of the isomorphism
  $S^n{\Bbb C}\stackrel{\sim}{\longrightarrow} {\Bbb C}^n$
 that sends $[\lambda_1, \,\cdots\,,\lambda_n]$
  to the monic polynomial $\prod_{i=1}^n(\lambda-\lambda_i)$
  of degree $n$ in $\lambda$.

Now identify
 $L$ with ${\Bbb A}^1$ by $y\mapsto \lambda$ and
 $E$ with the Chan-Paton space of $n$ D0-branes stacked
  at the origin $y=0$.
Then $\phi$ corresponds to a D0-brane configuration
  supported at $Z_{\phi}$,
 with the Chan-Paton space $\hat{V}_i$  and
 endomorphism algebra $\End(\hat{V}_i)=\End(V_i)$ at $Z_{\phi,i}$.
One may regard $Z_{\phi}$ as a deformation of the stacked D0-branes
 at $y=0$ (which corresponds to $\phi=0$).
{\it This reproduces also the Higgsing/un-Higgsing behavior
 of Polchinski's D-branes.}
Note that D0-branes on ${\Bbb A}^1$ described through this
 construction correponds to the locus $\image\Phi_{\scriptsizeChow}$
 in $\Map((\Space M_n({\Bbb C});{\Bbb C}^n),{\Bbb A}^1)$.

This spectral cover picture of D-branes is particularly fascinating
 when one recalls the Seiberg-Witten integrable system and
 the associated gauge-symmetry-breaking pattern revealed there;
 cf.\ [S-W1] and [Don2], [D-W], [Le].\footnote{However,
                              this setting has two drawbacks
                              one should be aware of:
                             (1) it obscures the important noncommutative
                                 nature of D-branes for it treats D-branes
                                 (of B-type) only as coherent torsion
                                 sheaves with a gauge symmetry, which we
                                 know now is not a complete picture,
                                (see also [Di-M] for subtleties in the
                                 case of D-brane bound-state systems),
                                  and
                             (2) while this construction is immediately
                                  generalizable to D-branes of complex
                                  codimension-$1$ in a complex target
                                  space, the further extension to
                                  describe higher-{\it co}dimensional
                                  D-branes becomes cumbersome.
                             These indicate that the spectral cover setting
                              might be just accidental for the cases it is
                              applicable and is overall not most natural
                              for D-branes. Cf.\ [Liu1].}

\bigskip

For the rest of this section, we will focus mainly
 on the moduli problem.

\bigskip

\subsection{D0-branes on the complex projective line ${\Bbb P}^1$.}

Let $Y$ be the projective line over ${\Bbb C}$:
 $$
  Y\; =\; {\Bbb P}^1\;
  =\; U_0\cup_{U_0\cap U_{\infty}} U_{\infty}\;
  =\; \Spec {\Bbb C}[y_0]
       \cup_{\scriptsizeSpec {\scriptsizeBbb C}[y_0\,,\,1/y_o] \simeq
         \scriptsizeSpec {\scriptsizeBbb C}[1/y_{\infty}\,,\,y_{\infty}]}
       \Spec {\Bbb C}[y_{\infty}]\,,
 $$
where
 $\Spec{\Bbb C}[y_0\,,\,1/y_o] \stackrel{\sim}{\rightarrow}
  \Spec{\Bbb C}[1/y_{\infty}\,,\,y_{\infty}]$
 is given by $y_{\infty}\mapsto 1/y_0$.
Having discussed the details of D0-branes on ${\Bbb A}^1$ in Sec.~4.1,
 we focus now on the issue of gluings for D0-branes on ${\Bbb P}^1$.

Recall
  the Grassmannian-like manifold $\Gr^{(2)}(n; d,n-d)$;
  the idempotents ${\mathbf 1}_d$, $d=0,\,\ldots\,,n$,
   in $M_n({\Bbb C})$;  and
  the notation `$m_1\sim m_2$' for similar matrices in $M_n({\Bbb C})$
 from Sec.~3.2.
Then, the ring-set representation variety
 $$
  \begin{array}{rcl}
   \Rep^{\ringsetscriptsize}({\Bbb C}[y], M_n({\Bbb C}))
    & =  & \{(e,m)\in M_n({\Bbb C})\times M_n({\Bbb C})\,:\,
                                          e^2=e\,,\; em=me=m\} \\[.6ex]
   & \subset  & {\Bbb A}^{n^2}\times {\Bbb A}^{n^2}
  \end{array}
 $$
 has $(n+1)$-many connected components, given by
 $$
  \Rep^{\ringsetscriptsize}({\Bbb C}[y], M_n({\Bbb C}))_{(d)}
   \; :=\; \{
    (e,m)\in\Rep^{\ringsetscriptsize}({\Bbb C}[y], M_n({\Bbb C}))\;:\;
      e\sim {\mathbf 1}_d
           \}\,,
 $$
 $d=0,\,\ldots\,,n$.
(Here we identify the pair $(e,m)$ with the ring-set-homomorphism
  $$
   \varphi_{(e,m)}:{\Bbb C}[y]\rightarrow M_n({\Bbb C})
    \hspace{2em}\mbox{with $1\mapsto e$ and $y\mapsto m$.})
  $$
$\Rep^{\ringsetscriptsize}({\Bbb C}[y], M_n({\Bbb C}))_{(d)}$
 is a $\GL_n({\Bbb C})$-manifold that goes with
 a natural $\GL_n({\Bbb C})$-equivariant bundle map
 $\Rep^{\ringsetscriptsize}({\Bbb C}[y], M_n({\Bbb C}))_{(d)}
  \rightarrow \Gr^{(2)}(n; d, n-d)$
 with fiber $\simeq M_d({\Bbb C})$.
In particular,
 $\dimm_{\scriptsizeBbb C}
   \Rep^{\ringsetscriptsize}({\Bbb C}[y], M_n({\Bbb C}))_{(d)}
   = d^2+2d(n-d) = n^2 -(n-d)^2$,
 which increases strictly when $d$ goes from $0$ to $n$.
The space $\Mor^{\ringsetscriptsize}({\Bbb C}[y], M_n({\Bbb C}))$
 of ring-set-homomorphisms from ${\Bbb C}[y]$ to $M_n({\Bbb C})$
 can be thought of as the $\GL_n({\Bbb C})$-space
 $\amalg_{d=0}^n
   \Rep^{\ringsetscriptsize}({\Bbb C}[y], M_n({\Bbb C}))_{(d)}$,
 but with the topology ${\cal T}$ in Definition 3.2.6.
It has the following properties:
\begin{itemize}
 \item[$\cdot$]
  $\Rep^{\ringsetscriptsize}({\Bbb C}[y], M_n({\Bbb C}))_{(n)}
   = \Rep({\Bbb C}[y], M_n({\Bbb C}))$
  is an open dense subset of\\
  $\Mor^{\ringsetscriptsize}({\Bbb C}[y], M_n({\Bbb C}))$.

 \item[$\cdot$]
  A neighborhood of $(e,m)$ with $e\sim {\mathbf 1}_d$
   consists of all
   $(e^{\prime},m^{\prime}) \in
      \Rep^{\ringsetscriptsize}({\Bbb C}[y], M_n({\Bbb C}))$
   such that
   \begin{itemize}
    \item[$\cdot$]
     $e^{\prime}\sim {\mathbf 1}_{d^{\prime}}$
     for some $d^{\prime}\ge d$;

    \item[$\cdot$]
     there is an idempotent $e^{\prime\prime}$
       in $Z(\langle e^{\prime},m^{\prime}\rangle)$ with the properties:
     \begin{itemize}
      \item[-]
       $e^{\prime\prime}\sim {\mathbf 1}_d$ and
       is in a neighborhood of $e$,

      \item[-]
       $e^{\prime\prime}m^{\prime}$ is in a neighborhood of $m$
       in $M_n({\Bbb C})$,

      \item[-]
       besides the characteristic value $0$ of multiplicity
        $d+(n-d^{\prime})$,
       the matrix
        $$
         (e^{\prime}-e^{\prime\prime})m^{\prime}\;
         =\; (e^{\prime}-e^{\prime\prime})m^{\prime}
                   (e^{\prime}-e^{\prime\prime}) \in M_n({\Bbb C})
        $$
        has all the remaining
        $(d^{\prime}-d)$-many characteristic values
        in a neighborhood of $\infty$ in ${\Bbb C}\cup \{\infty\}$.
     \end{itemize}
   \end{itemize}
\end{itemize}

The space $\Mor(\Space M_n({\Bbb C}),{\Bbb P}^1)$ of morphisms
 from $\Space M_n({\Bbb C})$ to ${\Bbb P}^1$
 is given by the locus in
 $\Mor^{\ringsetscriptsize}({\Bbb C}[y_0], M_n({\Bbb C}))
  \times \Mor^{\ringsetscriptsize}({\Bbb C}[y_{\infty}], M_n({\Bbb C}))$
 described by the following conditions:
\begin{itemize}
 \item[]\hspace{-1em}
  $( \varphi_{(e_{(0)}, m_{(0)})},
    \varphi_{(e_{(\infty)}, m_{(\infty)})} ) \in
  \Mor^{\ringsetscriptsize}({\Bbb C}[y_0], M_n({\Bbb C}))
  \times \Mor^{\ringsetscriptsize}({\Bbb C}[y_{\infty}], M_n({\Bbb C}))$,

 \item[(1)]
  $e_{(0)}e_{(\infty)}=e_{(\infty)}e_{(0)}$,
  $\,{\mathbf 1}=e_{(0)}+e_{(\infty)}-e_{(0)}e_{(\infty)}$;

 \item[(2)]
  $e_{(0)}m_{(\infty)}=m_{(\infty)}e_{(0)}$,
  $\,e_{(\infty)}m_{(0)}=m_{(0)}e_{(\infty)}$;

 \item[(3)]
  $e_{(\infty)}\langle e_{(0)}, m_{(0)}\rangle
   = e_{(0)}\langle e_{(\infty)}, m_{(\infty)} \rangle$ in $M_n({\Bbb C})$,
  (note that under Condition (2), \\
   $e_{(\infty)}\langle e_{(0)}, m_{(0)}\rangle
    =\langle e_{(\infty)}e_{(0)}, e_{(\infty)}m_{(0)}\rangle$ and
   $\,e_{(0)}\langle e_{(\infty)}, m_{(\infty)} \rangle
    =\langle e_{(0)}e_{(\infty)}, e_{(0)}m_{(\infty)}\rangle$);

 \item[(4)]
  $e_{(\infty)}m_{(0)}$ is invertible in
   $\langle e_{(\infty)}e_{(0)}, e_{(\infty)}m_{(0)}\rangle$,
  $\,e_{(0)}m_{(\infty)}$ is invertible in
    $\langle e_{(0)}e_{(\infty)}, e_{(0)}m_{(\infty)}\rangle$;

 \item[(5)]
 The identity in Condition (3) takes
  $e_{(\infty)}m_{(0)}$ to the inverse of $e_{(0)}m_{(\infty)}$ and
  $e_{(0)}m_{(\infty)}$ to the inverse of $e_{(\infty)}m_{(0)}$.
\end{itemize}
Note that Conditions (1) and (2) says that
 $$
  {\mathbf 1}\;\in\; \langle e_{(0)}, e_{(\infty)}\rangle\; \subset\;
  Z(\langle e_{(0)}, e_{(\infty)}, m_{(0)}, m_{(\infty)} \rangle)\;
  \subset\; M_n({\Bbb C})\,.
 $$
Conditions (3), (4), and (5) are
 the descendability to localizations and
 the gluability of pairs of ring-set-morphisms in
 $\Mor^{\ringsetscriptsize}({\Bbb C}[y_0], M_n({\Bbb C})) \times
  \Mor^{\ringsetscriptsize}({\Bbb C}[y_{\infty}], M_n({\Bbb C}))$.
$\GL_n({\Bbb C})$ acts diagonally on
 $\Mor^{\ringsetscriptsize}({\Bbb C}[y_0], M_n({\Bbb C})) \times
  \Mor^{\ringsetscriptsize}({\Bbb C}[y_{\infty}], M_n({\Bbb C}))$,
 leaving Conditions (1) - (5) invariant.

\bigskip

\noindent
{\bf Lemma 4.2.1 [closed condition].} {\it
 Assuming Conditions $(1)$ and $(2)$,
 then Conditions $(3)$, $(4)$, and $(5)$ together are equivalent to
 \begin{itemize}
  \item[$(3^{\prime})$]
   $e_{(0)}e_{(\infty)}m_{(0)}m_{(\infty)}=e_{(0)}e_{(\infty)}$.
 \end{itemize}
 In particular, the system $\{(1), (2), (3), (4), (5)\}$ realizes
  $\Mor(\Space M_n({\Bbb C}),{\Bbb P}^1)$
  as a $\GL_n({\Bbb C})$-invariant closed subset in
  $\Mor^{\ringsetscriptsize}({\Bbb C}[y_0], M_n({\Bbb C})) \times
   \Mor^{\ringsetscriptsize}({\Bbb C}[y_{\infty}], M_n({\Bbb C}))$.
} 

\bigskip

\noindent
$\Map((\Space M_n({\Bbb C});{\Bbb C}^n),{\Bbb P}^1)
 =\Mor(\Space M_n({\Bbb C}),{\Bbb P}^1)/\!\sim\,$
 is now given by the orbit-space of the $\GL_n({\Bbb C})$-action
 on the above locus in
 $\Mor^{\ringsetscriptsize}({\Bbb C}[y_0], M_n({\Bbb C})) \times
  \Mor^{\ringsetscriptsize}({\Bbb C}[y_{\infty}], M_n({\Bbb C}))$.

For
 ${\cal R}=
  (\varphi_{(e_{(0)},m_{(0)})},\varphi_{(e_{(\infty)},m_{(\infty)})})
  \in \Map((\Space M_n({\Bbb C});{\Bbb C}^n),{\Bbb P}^1)$,
the Chan-Paton module on each local chart $U$,
  where $U=U_0$ or $U_{\infty}$,
 is given by the $(e,m)$-module $e\cdot{\Bbb C}$ but now regarded
 as a ${\Bbb C}[y]$-module
 $_{{\scriptsizeBbb C}[y]}(e\cdot{\Bbb C})$ via $\varphi_{(e,m)}$.
We will denote this ${\cal O}_U$-module on $U$ by
 $\hat{\varphi}_{(e,m)\ast}(e\cdot{\Bbb C}^n)$.
It is supported on the image scheme $\image\hat{\varphi}$ on $U$
 associated to the ideal $\Ker\varphi_{(e,m)}$ in ${\Bbb C}[y]$.
Here,
 $(e,m)=(e_{(0)},m_{(0)})$ or $(e_{(\infty)},m_{(\infty)})$ respectively
  and
 ${\Bbb C}[y]={\Bbb C}[y_0]$ or ${\Bbb C}[y_{\infty}]$ respectively.
Except that $e\cdot{\Bbb C}^n$ now replaces ${\Bbb C}^n$,
 all the local details of $\hat{\varphi}_{(e,m)\ast}(e\cdot{\Bbb C}^n)$
 are the same as those in the case $Y={\Bbb A}^1$.
The total length of $\hat{\varphi}_{(e,m)\ast}(e\cdot{\Bbb C}^n)$ is
 $\dimm_{\scriptsizeBbb C}(e\cdot{\Bbb C}^n)$,
 ($=d$ for $e\sim {\mathbf 1}_d$).
The pair
 $\{ \image\hat{\varphi}_{(e_{(0)},m_{(0)})}\,,\,
     \image\hat{\varphi}_{(e_{(\infty)},m_{(\infty)})} \}$
 of local image schemes glue to a $0$-dimensional subscheme,
  denoted $\image\hat{\varphi}_{\cal R}$ or
  $\hat{\varphi}_{\cal R}(\Space M_n({\Bbb C}))$,
 of length $\le n$ on ${\Bbb P}^1$.
Idempotency of $e_{\bullet}$ and Conditions (1) and (2) imply that
 $\{ \hat{\varphi}_{(e_{(0)},m_{(0)})\ast}(e_{(0)}\cdot{\Bbb C}^n)\,,\,
     \hat{\varphi}_{(e_{(\infty)},m_{(\infty)})\ast}
                                    (e_{(\infty)}\cdot{\Bbb C}^n) \}$
 glues to a (torsion) ${\cal O}_{{\scriptsizeBbb P}^1}$-module
 on ${\Bbb P}^1$.
This is the push-forward $\hat{\phi}_{{\cal R}\ast}{\Bbb C}^n$ of
 ${\Bbb C}^n$ on $\Space M_n({\Bbb C})$ to ${\Bbb P}^1$
 under $\hat{\phi}_{\cal R}$;
cf.\ footnote 32.
It is the Chan-Paton module of the D0-branes
 $\hat{\varphi}(\Space M_n({\Bbb C}))$ on ${\Bbb P}^1$
 in the current setting.
Note that the total length of $\hat{\phi}_{{\cal R}\,\ast}{\Bbb C}^n$
 on ${\Bbb P}^1$ remains $n$.
The Higgsing/un-Higgsing behavior of Chan-Paton modules of D0-branes
 on any target $Y$ is a local issue and hence, for $Y={\Bbb P}^1$,
 is the same as that for $Y={\Bbb A}^1$ in Sec.~4.1.

The local discussions in Sec.~4.1 can be glued to global statements.
In particular,

\bigskip

\noindent
{\bf Proposition 4.2.2 [D0-branes on ${\Bbb P}^1$].} {\it
 There is an embedding
  $\Phi_{\scriptsizeHilb}:
   \Hilb^n_{{\scriptsizeBbb P}^1}=:({\Bbb P}^1)^{[n]}
    \rightarrow \Map((\Space M_n({\Bbb C});{\Bbb C}^n),{\Bbb P}^1)$,
  whose image is characterized by $\varphi_{\cal R}$ such that
   $\image\hat{\varphi}_{\cal R}$ is a subscheme of length $n$
   on ${\Bbb P}^1$.
 There is an embedding
  $\Phi_{\scriptsizeChow}:
   S^n({\Bbb P}^1)
    \rightarrow \Map((\Space M_n({\Bbb C});{\Bbb C}^n),{\Bbb P}^1)$,
  whose image is characterized by $\varphi_{\cal R}$ such that
   $\image\hat{\varphi}_{\cal R}$ is a reduced subscheme
   (of length $\le n$) on ${\Bbb P}^1$.
 There is a map
  $\Map((\Space M_n({\Bbb C});{\Bbb C}^n),{\Bbb P}^1)
   \rightarrow S^n({\Bbb P}^1)$
  that has $\Phi_{\scriptsizeChow}$ as a section.
 The pattern of open-string-induced Higgsing/un-Higgsing behavior of
  $n$ D0-branes on ${\Bbb P}^1$ can be reproduced
  in the current content via deformations of morphisms
  $[\varphi_{\cal R}]$
  in $\Phi_{\scriptsizeChow}(S^n({\Bbb P}^1))
      \subset \Map((\Space M_n({\Bbb C}); {\Bbb C}^n),{\Bbb P}^1)$.
} 

\bigskip

\noindent
{\it Remark 4.2.3 $[$strict morphism$]$.}
 A strict morphism from $\Space M_n({\Bbb C})$ to ${\Bbb P}^1$
  is given by a strict morphism
  (cf.\ Definition 1.2.11 and Definition 1.1.1)
  from $[( \{ {\Bbb C}[y_0], {\Bbb C}[y_{\infty}] \}
          \doublearrow \{ {\Bbb C}[y, 1/y]\} )]$
  to $[\{M_n({\Bbb C})\}]$.
 Since $Z(M_n({\Bbb C}))={\Bbb C}$,
 such a morphism factors as
  $$
   [( \{ {\Bbb C}[y_0], {\Bbb C}[y_{\infty}] \}
            \doublearrow \{ {\Bbb C}[y, 1/y]\} )]\;
    \longrightarrow\; [\{{\Bbb C}\}]\;
    \longrightarrow\; [\{M_n({\Bbb C})\}]
  $$
  and, hence, corresponds to a morphism
  $\Spec{\Bbb C}\rightarrow {\Bbb P}^1$.
 The corresponding D0-brane on ${\Bbb P}^1$ is supported at
  a reduced ${\Bbb C}$-point on ${\Bbb P}^1$
  with the Chan-Paton module ${\Bbb C}^n$,
  i.e.\ $n$-many coincident D0-branes on ${\Bbb P}^1$
        in the picture of Polchinski.
 The moduli space of such morphisms (i.e.\ coincident D0-branes)
  is ${\Bbb P}^1$.
 Thus, we see that the inclusion of general morphisms
   (cf.\ Definition 1.2.14 and Definition 1.1.1)
   in the definition of $\Mor(\Space M_n({\Bbb C}), {\Bbb P}^1)$
    and, hence,
   in the definition of $\Map((\Space M_n({\Bbb C});{\Bbb C}^n),{\Bbb P}^1)$
  is also required
 if one wants to incorporate the Higgsing/un-Higgsing behavior of,
  in this case, D0-branes on ${\Bbb P}^1$.
 Similar phenomenon occurs for other projective target spaces as well.
 This is another incident of the mysterious harmony
  between stringy requirement and mathematical naturality
  for a string-theory-related mathematical object.

\bigskip

\subsection{D0-branes on the complex affine plane ${\Bbb A}^2$.}

For a commutative $Y$ of dimension $\ge 2$, an additional ingredient
 than those in Sec.~4.1 and Sec.~4.2 is
                            commuting schemes/varieties\footnote{For
                                pure algebraic geometers:
                               Moduli problems in commutative algebraic
                                geometry tends to boil down to
                                Hilbert schemes, which in projective cases
                                are realized as a locus in an appropriate
                                Grassmannian variety.
                               In that sense, commuting schemes/varieties
                                play the same fundamental role
                                as Grassmannian varieties do
                                for the moduli problem of morphisms
                                from an Azumaya-type noncommutative space
                                to a commutative variety.
                               We hope this gives further motivation
                                to study commuting schemes/varieties.
                               See, e.g., [Bas], [Ge], [Ri], [S-T],
                                [Vac1], [Vac2].}.
We discuss in this subsection the case $Y={\Bbb A}^2$,
 for which the commuting variety that occurs is known slightly better.

Let $Y={\Bbb A}^2=\Spec {\Bbb C}[y_1,y_2]$
 be the affine plane over ${\Bbb C}$.
Then
 $\Mor(\Space M_n({\Bbb C}), Y) =\Mor({\Bbb C}[y_1,y_2], M_n({\Bbb C}))$
 is the variety\footnote{Throughout, we only consider
                                     the reduced scheme structure
                                     on a commuting scheme or
                                     a representation scheme
                                     that occurs in the problem.}
 what parameterizes the elements in the set
 $$
  C_2M_n({\Bbb C})\; :=\;
   \{(m_1,m_2)\in M_n({\Bbb C})\times M_n({\Bbb C})\,:\, m_1m_2=m_2m_1\}
 $$
 of pairs of commuting matrices in $M_n({\Bbb C})$.
This variety is identical with ${\Bbb C}$-algebra representation variety
 $\Rep({\Bbb C}[y_1,y_2], M_n({\Bbb C}))$
 with a point represented by
  $(m_1,m_2) \in M_n({\Bbb C})\times M_n({\Bbb C})$
  corresponding to the ${\Bbb C}$-algebra-homomorphism
  $$
   \varphi_{(m_1,m_2)}: {\Bbb C}[y_1,y_2]\rightarrow M_n({\Bbb C})\,,
   \hspace{1ex}
   \mbox{generated by $1\mapsto{\mathbf 1}$,
         $y_1\mapsto m_1$, and $y_2\mapsto m_2$}\,.
  $$

\bigskip

\noindent
{\bf Proposition 4.3.1 [irreducibility].} ([Ge], [Bas], and [Vac2].) {\it
 $\Rep({\Bbb C}[y_1,y_2], M_n({\Bbb C}))$ is an irreducible variety
  of dimension $n^2+n$ in ${\Bbb A}^{n^2}\times {\Bbb A}^{n^2}$.
 The $\GL_n({\Bbb C})$-action on $\Rep({\Bbb C}[y_1,y_2], M_n({\Bbb C}))$
  has stabilizer subgroups of minimal dimension $n$.
 A generic $\GL_n({\Bbb C})$-orbit thus has dimension $n^2-n$,
  that achieves the maximum orbit-dimension  and
 the subset that consists of $\varphi_{(m_1,m_2)}$,
  where $(m_1,m_2)$ is a diagonalizable commuting pair
   with both $m_1$ and $m_2$ having distinct characteristic values, is
   a smooth open dense subset in $\Rep({\Bbb C}[y_1,y_2], M_n({\Bbb C}))$.
} 

\bigskip

\noindent
It follows that
 $$
  \Map((\Space M_n({\Bbb C});{\Bbb C}^n), {\Bbb A}^2)\;
  \simeq\; \Map(\Space M_n({\Bbb C}), {\Bbb A}^2)\;
  =\; \Rep({\Bbb C}[y_1,y_2],M_n({\Bbb C}))/\!\sim\,,
 $$
 the orbit-space of the $\GL_n({\Bbb C})$-action with the quotient
 topology, is a connected non-Hausdorff topological space
 that contains a connected smooth open dense Hausdorff subset of
 dimension $2n$, namely the subset of $S^n({\Bbb A}^2)$
 that consists of $[(\lambda_1,\mu_1),\,\cdots\,,(\lambda_n,\mu_n)]$
 such that
  $\lambda_i$, $i=1,\,\ldots\,,n$, are all distinct from each other
  and so are $\mu_i$, $i=1,\,\ldots\,,n$.
Here $S^n({\Bbb A}^2) := ({\Bbb A}^2)^{n}/Sym_n$ is the $n$-th symmetric
 product of ${\Bbb A}^2$.

The complete set of dominance relations of the $\GL_n({\Bbb C})$-orbits
 in $\Rep({\Bbb C}[y_1,y_2], M_n({\Bbb C}))$,
 which generalizes [Ge], are not known.
However, there are two distinguished Hausdorff subspace in
 $\Map((\Space M_n({\Bbb C});{\Bbb C}^n), {\Bbb A}^2)$
 that can be understood
 through the work of Nakajima [Na] and of Vaccarino [Vac2]:
 \begin{itemize}
  \item[(1)]
   the naturally embedded image of the Hilbert scheme
   $({\Bbb A}^2)^{[n]}:= \Hilb^n_{{\scriptsizeBbb A}^2}$
   (with the reduced scheme structure) of $0$-dimensional subschemes
   of length $n$ on ${\Bbb A}^2$;

  \item[(2)]
   the naturally embedded image of the Chow variety
   $\Chow_{{0, \scriptsizeBbb A}^2}^{(n)} = S^n({\Bbb A}^2)$
   of $0$-cycles of order $n$ on ${\Bbb A}^2$.
 \end{itemize}
We now explain the details.

\bigskip

\noindent
{\bf Proposition 4.3.2 [regular representation].} {\it
 Let $R$ be a commutative Artinian algebra over ${\Bbb C}$
  of dimension $n$.
 Then, the regular representation\footnote{Recall that
                                 a {\it regular representation} of
                                 an algebra $R$ is the representation
                                 of $R$ on $R$ itself by,
                                 in our convention, left multiplications;
                                 i.e.\ $R$ as a (left) $R$-module.}
  of $R$ realizes $R$ as a maximal commutative subalgebra $R^{\prime}$
  of $M_n({\Bbb C})$.
 Furthermore, as an $R^{\prime}$-module,
  $_{R^{\prime}}{\Bbb C}^n\simeq R^{\prime}$.
} 

\bigskip

\noindent
{\it Proof.}
 This is an immediate corollary of [S-T: Sec.2.7, Theorem 11].
 When $R$ is generated by two commuting elements and the identity,
  as is in our case, there are two other independent proofs:
 (1) The first part of the proof of [Na: Sec.~1.2, Theorem 1.9]
     can be adapted directly to give another more analytic proof
     of the statement, cf.\ proof of Proposition 4.3.3 below.
 (2) This is a corollary of [Ge], which says that the
     maximum dimension of a commutative subalgebra in $M_n({\Bbb C})$
     generated by two commuting matrices and the identity is $n$.

\noindent\hspace{15cm}$\Box$

\bigskip

\noindent
Note that, in the above statement, different choices of
 $R\simeq {\Bbb C}^n$ as ${\Bbb C}$-vector spaces give rise to
 $R^{\prime}$'s in the same adjoint $\GL_n({\Bbb C})$-orbit.
It follows that there is an embedding of sets
 $$
  \Phi_{\scriptsizeHilb}: ({\Bbb A}^2)^{[n]}\;
    \longrightarrow\;
    \Map((\Space M_n({\Bbb C});{\Bbb C}^n), {\Bbb A}^2)\,,
  \hspace{1em}
   {\Bbb C}[y_1,y_2]/I\;
      \longmapsto\; \varphi_{(m_1,m_2)}\,.
 $$
Here,
 $I$ is an ideal of ${\Bbb C}[y_1,y_2]$ so that
  $\dimm_{\scriptsizeBbb C}({\Bbb C}[y_1,y_2]/I)=n$;
 it gives then the subalgebra
  $({\Bbb C}[y_1,y_2]/I)^{\prime}\subset M_n({\Bbb C})$
  as in Proposition 4.3.2, unique up conjugation;
 the corresponding matrix $m_i$ for $y_i$, $i=1,\,2$.
  under the built-in ${\Bbb C}$-algebra-isomorphism
  ${\Bbb C}[y_1,y_2]/I
   \stackrel{\sim}{\rightarrow} ({\Bbb C}[y_1,y_2]/I)^{\prime}$
 determines then $\varphi_{(m_1,m_2)}$.

\bigskip

\noindent
{\bf Proposition 4.3.3 [stable subset].} (Cf.\ [Na: Theorem 1.9].)
{\it
 Let $\Rep({\Bbb C}[y_1,y_2], M_n({\Bbb C}))^{\scriptsizest}$
  be the subset of $\Rep({\Bbb C}[y_1,y_2], M_n({\Bbb C}))$
  that consists of $\varphi_{(m_1,m_2)}$ such that
  $_{\langle{\mathbf 1}, m_1, m_2\rangle}{\Bbb C}^n
   \simeq \langle{\mathbf 1}, m_1, m_2\rangle$
  as $\langle{\mathbf 1}, m_1, m_2\rangle$-modules.
 Then
  $\Rep({\Bbb C}[y_1,y_2], M_n({\Bbb C}))^{\scriptsizest}$
  is smooth and $\GL_n({\Bbb C})$-invariant
 with stabilizers all of the same dimension $n$.
} 

\bigskip

\noindent
{\it Proof.}
 This is actually [Na: Theorem 1.9] in disguise.
 Note that the stability condition in the defining condition
  of the set $\tilde{H}$ in ibidem is precisely the condition
  ``$_{\langle{\mathbf 1}, m_1, m_2\rangle}{\Bbb C}^n
   \simeq \langle{\mathbf 1}, m_1, m_2\rangle$
    as $\langle{\mathbf 1}, m_1, m_2\rangle$-modules"
  in the statement here.
 Having said so, let us give a sketch of the proof in terms of
  the current setting.

 Using the trace map $M_n({\Bbb C})\rightarrow {\Bbb C}$
  as a complex bilinear inner product on the ${\Bbb C}$-vector space
  $M_n({\Bbb C})$,
 one can show that the (analytic quadric) commutator map
  (on analytic spaces)
  $$
   c\;:\;
    M_n({\Bbb C})\times M_n({\Bbb C})\;\longrightarrow\; M_n({\Bbb C})\,,
   \hspace{1em} (m_1,m_2)\;\longmapsto\; [m_1,m_2]:= m_1m_2-m_2m_1
  $$
  has cokernel $\coker dc_{(m_1,m_2)}$ at $(m_1,m_2)$ being
  $\{\xi\in M_n({\Bbb C}):[\xi,m_1]=[\xi,m_2]=0\}$,
  i.e.\ the centralizer
   $\Centralizer\langle{\mathbf 1}, m_1, m_2\rangle$ of the subalgebra
   $\langle{\mathbf 1}, m_1, m_2\rangle$ in the algebra $M_n({\Bbb C})$.
 Note that for $(m_1,m_2)\in C_2M_n({\Bbb C})$,
  $\langle{\mathbf 1}, m_1, m_2\rangle \subset
      \Centralizer\langle{\mathbf 1}, m_1, m_2\rangle$.

 If, furthermore,
   $\varphi_{(m_1,m_2)} \in
             \Rep({\Bbb C}[y_1,y_2], M_n({\Bbb C}))^{\scriptsizest}$,
  then
   $_{\langle{\mathbf 1}, m_1, m_2\rangle}{\Bbb C}^n
    = \langle{\mathbf 1}, m_1, m_2\rangle\cdot v_0$
   for some $v_0\in {\Bbb C}^n$.
 The ${\Bbb C}$-linear map
  $\Centralizer\langle{\mathbf 1}, m_1, m_2\rangle\rightarrow {\Bbb C}^n$,
   defined by $\xi\mapsto \xi\cdot v_0$, is then invertible
  and hence a ${\Bbb C}$-vector-space-isomorphism.
 It follows that
  $\langle{\mathbf 1}, m_1, m_2\rangle =
       \Centralizer\langle{\mathbf 1}, m_1, m_2\rangle$
  and $\dim_{\scriptsizeBbb C}\coker dc_{(m_1,m_2)}=n$
  for $\varphi_{(m_1,m_2)}
         \in \Rep({\Bbb C}[y_1,y_2], M_n({\Bbb C}))^{\scriptsizest}$.
 This shows that
  $\Rep({\Bbb C}[y_1,y_2], M_n({\Bbb C}))^{\scriptsizest}$
  is smooth.

 Finally, note that
 $\Stab(\varphi_{(m_1,m_2)})
  =\GL_n({\Bbb C})\cap\Centralizer\langle{\mathbf 1}, m_1, m_2\rangle$,
 which has the same dimension as
  $\Centralizer\langle{\mathbf 1}, m_1, m_2\rangle$.
 The proposition follows.

\noindent\hspace{15cm}$\Box$

\bigskip

Since the closure of $\overline{O}$ of a $G$-orbit $O$ of an action
 of a reductive algebraic group $G$ on an affine variety $V$
 (both over ${\Bbb C}$) is a union of $O$ with $G$-orbits of
 strictly smaller dimension, one has:

\bigskip

\noindent
{\bf Corollary 4.3.4 [good quotient].} {\it
 All the $\GL_n({\Bbb C})$-orbits are closed in
  $\Rep({\Bbb C}[y_1,y_2], M_n({\Bbb C}))^{\scriptsizest}$  and
 the map
  $\Rep({\Bbb C}[y_1,y_2], M_n({\Bbb C}))^{\scriptsizest}
   \rightarrow \image\Phi_{\scriptsizeHilb}$
  to the orbit-space is a good quotient.
} 

\bigskip

\noindent
This realizes the map
 $\Phi_{\scriptsizeHilb}: ({\Bbb A}^2)^{[n]}
    \rightarrow \Map((\Space M_n({\Bbb C});{\Bbb C}^n), {\Bbb A}^2)$
 as an embedding of the (reduced) Hilbert scheme
 as a variety/analytic space.

Let ${\Bbb C}^n$ parameterizes the diagonal matrices in $M_n({\Bbb C})$.
Then, the embedding
 $$
  \begin{array}{ccc}
   {\Bbb C}^n\times {\Bbb C^n}=({\Bbb C}^2)^n
    & \hookrightarrow  & M_n({\Bbb C})\times M_n({\Bbb C}) \\[.6ex]
   ((\lambda_1,\mu_1),\,\cdots\,,\,(\lambda_n,\mu_n))
    & \mapsto
    & (\Diag(\lambda_1,\,\cdots\,,\,\lambda_n)\,,\,
       \Diag(\mu_1,\,\cdots\,,\mu_n))
  \end{array}
 $$
 descends to an embedding
 $$
  \begin{array}{ccccc}
   \Phi_{\scriptsizeChow}  & :   & S^n({\Bbb A}^2)
    & \longrightarrow
    & \Map((\Space M_n({\Bbb C});{\Bbb C}^n), {\Bbb A}^2) \\[.6ex]
   && [(\lambda_1,\mu_1),\,\cdots\,,(\lambda_n,\mu_n)]  & \longmapsto\;
    & \varphi_{( \scriptsizeDiag(\lambda_1,\,\cdots\,,\,\lambda_n)\,,\,
                 \scriptsizeDiag(\mu_1,\,\cdots\,,\mu_n) )}
  \end{array}
 $$
 of the Chow variety.

$S^n({\Bbb A}^2)$ is the categorical quotient of
 $\Rep({\Bbb C}[y_1,y_2], M_n({\Bbb C}))$
 under the adjoint $\GL_n({\Bbb C})$-action.
The affine morphism
 $\Rep({\Bbb C}[y_1,y_2], M_n({\Bbb C}))\rightarrow S^n({\Bbb A}^2)$
 induced by the $\GL_n({\Bbb C})$-invariant function ring
 on $\Rep({\Bbb C}[y_1,y_2], M_n({\Bbb C}))$ descends to
 a morphism
 $\image\Phi_{\scriptsizeHilb}\rightarrow \image\Phi_{\scriptsizeChow}$
 of varieties that realizes
 $({\Bbb A}^2)^{[n]}$ as a desingularization of $S^n({\Bbb A}^2)$.
$\image\Phi_{\scriptsizeChow}$ consists of all the closed points
 in $\Map((\Space M_n({\Bbb C});{\Bbb C}^n), {\Bbb A}^2)$  and
the closure of any point in
 $\Map((\Space M_n({\Bbb C});{\Bbb C}^n), {\Bbb A}^2)$
 contains a unique point in $\image\Phi_{\scriptsizeChow}$.
Cf.\ [Na], [Pro], [Ri], and [Vac2].

Note that, for $(m_1,m_2)\in C_2M_n({\Bbb C})$,
 as $m_1$ and $m_2$ commute, they can be simultaneously triangularized.
If they have a simultaneous triangularization with the diagonal entries
 $(\lambda_1,\,\cdots\,,\lambda_n)$ and $(\mu_1,\,\cdots\,,\mu_n)$
 respectively,
let
 $I_{\{(\lambda_1,\mu_1),\,\cdots\,,(\lambda_n,\mu_n)\}}\;
  :=\; (y_1-\lambda_1, y_2-\mu_1)
        \cap\,\cdots\,\cap (y_1-\lambda_n,y_2-\mu_n)$
 be the ideal in ${\Bbb C}[y_1,y_2]$
 for the set of closed points
 $\{(\lambda_1,\mu_1),\,\cdots\,, (\lambda_n, \mu_n)\}$
 (as points on the analytic space ${\Bbb C}^2$
  with repeated points dropped).
Then,
 $$
  I_{\{(\lambda_1,\mu_1),\,\cdots\,,(\lambda_n,\mu_n)\}}^{\;n}\;
   \subset\; \Ker\varphi_{(m_1,m_2)}\;
   \subset\; I_{\{(\lambda_1,\mu_1),\,\cdots\,,(\lambda_n,\mu_n)\}}\,.
 $$
In particular,
 the image scheme
 $\image\hat{\varphi}_{(m_1,m_2)}\simeq
  \Spec({\Bbb C}[y_1,y_2]/\Ker\varphi_{(m_1,m_2)})$ on ${\Bbb A}^2$
 has the reduced scheme structure exactly the set
 $\{(\lambda_1,\mu_1),\,\cdots\,, (\lambda_n, \mu_n)\}$ above.

For $\varphi_{(m_1,m_2)}\in\image\Phi_{\scriptsizeHilb}$,
 the Chan-Paton module $\hat{\varphi}_{(m_1,m_2)\ast}{\Bbb C}^n$,
 as a ${\cal O}_{\scriptsizeimage\hat{\varphi}_{(m_1,m_2)}}$-module,
 is isomorphic to the structure sheaf
 ${\cal O}_{\scriptsizeimage\hat{\varphi}_{(m_1,m_2)}}$
while for $\varphi_{(m_1,m_2)}\in\image\Phi_{\scriptsizeChow}$,
 the Chan-Paton module $\hat{\varphi}_{(m_1,m_2)\ast}{\Bbb C}^n$,
 as a ${\cal O}_{\scriptsizeimage\hat{\varphi}_{(m_1,m_2)}}$-module,
 is isomorphic to
 $\oplus_{i=1}^n{\cal O}_{(\lambda_i,\mu_i)}$.
Here, $(\lambda_i,\mu_i)$, $i=1,\,\ldots\,,n$,
 are the image point from earlier notations with repeated
 $(\lambda_i,\mu_i)$ kept to contribute to the direct sum.
Behavior of Higgsing/un-Higgsing follows similar pattern
 as in Sec.~4.1.

\bigskip

\subsection{D0-branes on a complex quasi-projective variety.}

A picture of D0-branes on a (commutative) complex quasi-projective
 variety that follows from a combination and an immediate generalization
 of Sec.~4.1 - Sec.~4.3 is given in this subsection.
A comparison with gas of D0-branes in [Vafa1] of Vafa is given in the end.

\bigskip

\begin{flushleft}
{\bf D0-branes on ${\Bbb P}^r$.}
\end{flushleft}
Let $Y$ be the projective space over ${\Bbb C}$:
 $$
  Y\;=\; {\Bbb P}^r\; =\; \Proj {\Bbb C}[y_0,y_1,\,\cdots\,, y_r]\;
   =\; \cup_{i=0}^r\, U_i\;
   =\; \cup_{i=0}^r\,
         \Spec {\Bbb C}[\mbox{$\frac{y_0}{y_i}$}\,,\,\cdots\,,\,
                        \mbox{$\frac{y_r}{y_i}$}]\,.
 $$
Here
 $y_{\bullet}/y_i$ are treated as formal variables
  with $y_i/y_i =$ the identity $1$ of the ring
  ${\Bbb C}[\mbox{$\frac{y_0}{y_i}$}\,,\,\cdots\,,\,
            \mbox{$\frac{y_r}{y_i}$}]$;
 the gluings
  $U_i \supset U_{ij}:= U_i\cap U_j\stackrel{\sim}{\leftarrow}
       U_{ji}:= U_j\cap U_i \subset U_j$
  of local affine charts are given by
 $$
  \begin{array}{cccccccl}
   {\Bbb C}[\frac{y_0}{y_i}\,,\,\cdots\,,\,\frac{y_r}{y_i}]
    & \hookrightarrow
    & \frac{{\smallBbb C}[\frac{y_0}{y_i}\,,\,\cdots\,,\,
                          \frac{y_r}{y_i}\,,\,\frac{y_i}{y_j}]}
           {\left( \frac{y_j}{y_i}\cdot\frac{y_i}{y_j}-1 \right)}
    & \stackrel{\sim}{\longrightarrow}
    & \frac{{\smallBbb C}[\frac{y_0}{y_j}\,,\,\cdots\,,\,
                          \frac{y_r}{y_j}\,,\,\frac{y_j}{y_i}]}
           {\left( \frac{y_i}{y_j}\cdot\frac{y_j}{y_i}-1 \right)}
    & \hookleftarrow
    & {\Bbb C}[\frac{y_0}{y_j}\,,\,\cdots\,,\,\frac{y_n}{y_j}]   \\[3ex]
   && \frac{y_{\bullet}}{y_i}
    & \longmapsto  & \frac{y_{\bullet}}{y_j}\cdot\frac{y_j}{y_i} \\[1ex]
   && \frac{y_i}{y_j}  & \longmapsto  & \frac{y_i}{y_j}   &&&.
  \end{array}
 $$

Let
 $$
  C_{r+1}M_n({\Bbb C})\; :=\;
   \{ (m_0,\,\cdots\,,m_r)\in M_n({\Bbb C})^{r+1}\,:\, m_im_j=m_jm_i\,,
                              i,j = 0,\,\ldots\,, r \}\,.
 $$
The ring-set representation variety
 \begin{eqnarray*}
  \lefteqn{
   \Rep^{\ringsetscriptsize}(
   {\Bbb C}[\mbox{$\frac{y_0}{y_i}$}\,,\,\cdots\,,\,
            \mbox{$\frac{y_n}{y_i}$}], M_n({\Bbb C}))  }\\[.6ex]
  && =\; \{(m_{(i),0},\,\cdots\,,m_{(i),r})
          \in C_{r+1}M_n({\Bbb C})\,:\,
             m_{(i),i}m_{(i),i^{\prime}}=m_{(i),i^{\prime}}m_{(i),i}
             =m_{(i),i^{\prime}},\, i^{\prime}=0,\,\ldots\,,r \} \\[.6ex]
  && \subset\; \prod_{r+1}{\Bbb A}^{n^2}\;=\; {\Bbb A}^{n^2(r+1)}\,,
 \end{eqnarray*}
 (in particular, $e_{(i)} := m_{(i),i}$ is an idempotent),
 is a disjoint union of
 \begin{eqnarray*}
  \lefteqn{
   \Rep^{\ringsetscriptsize}(
    {\Bbb C}[\mbox{$\frac{y_0}{y_i}$}\,,\,\cdots\,,\,
             \mbox{$\frac{y_n}{y_i}$}],M_n({\Bbb C}) )_{(d)}  }\\[.6ex]
  && :=\; \{
    (m_{(i),\bullet})_{\bullet} \in
      \Rep^{\ringsetscriptsize}(
      {\Bbb C}[\mbox{$\frac{y_0}{y_i}$}\,,\,\cdots\,,\,
               \mbox{$\frac{y_n}{y_i}$}], M_n({\Bbb C}) )\;:\;
     m_{(i),i}\sim {\mathbf 1}_d
           \}\,,
  \hspace{1em} d\,=\,0,\,\ldots\,,n\,.
 \end{eqnarray*}
Here, again, we identify
 the ring-set-homomorphism
  $\varphi_{(m_{(i),0},\,\cdots\,,m_{(i),r})}:
   {\Bbb C}[\frac{y_0}{y_i}\,,\,\cdots\,,\,\frac{y_r}{y_i}]
   \rightarrow M_n({\Bbb C})$
  that sends $y_{\bullet}/y_i$ to $m_{(i),\bullet}$
 with $(m_{(i),0},\,\cdots\,,m_{(i),r})\in C_{r+1}M_n({\Bbb C})$.

Similar to the case $Y={\Bbb P}^1$ in Sec.~4.2,
the space $\Mor(\Space M_n({\Bbb C}),{\Bbb P}^r)$ of morphisms
 from $\Space M_n({\Bbb C})$ to ${\Bbb P}^r$
 is given by the locus on
 $\prod_{i=0}^r\,
  \Mor^{\ringsetscriptsize}({\Bbb C}[\frac{y_0}{y_i}\,,\,
                          \cdots\,,\,\frac{y_r}{y_i}], M_n({\Bbb C}))$
 described by the following conditions:\footnote{For
                                            readers who are familiar
                                            with toric geometry:
                                           Such system of conditions
                                            can be formally associated
                                            to the fan (or polytope
                                            in the projective case)
                                            of a toric variety.}
\begin{itemize}
\item[]\hspace{-1em}
 $( \varphi_{(m_{(0),\bullet})_{\bullet}}\,,\,
    \cdots\,,\,\varphi_{(m_{(r),\bullet})_{\bullet}} )  \in
  \prod_{i=0}^r\,
  \Mor^{\ringsetscriptsize}({\Bbb C}[\frac{y_0}{y_i}\,,\,
                        \cdots\,,\,\frac{y_r}{y_i}], M_n({\Bbb C}))$,

 \item[(1)]
  $m_{(i), i}m_{(j),j}=m_{(j),j}m_{(i),i}\,$,\\[1ex]
  ${\mathbf 1}
    =\sum_i m_{(i),i} -\sum_{i_1<i_2} m_{(i_1),i_1}m_{(i_2),i_2}
     + \,\cdots\,\\[.6ex]
     \mbox{\hspace{6em}}
      + (-1)^{k-1}\sum_{i_1<\,\cdots\,<i_k}
                       m_{(i_1),i_1}\,\cdots\,m_{(i_k),i_k}
      +\,\cdots\,+(-1)^r m_{(0),0}\,\cdots\,m_{(r),r}$;

 \item[(2)]
  $m_{(i),i}m_{(j),\bullet}=m_{(j),\bullet}m_{(i),i}$,
  $\,i,j, \bullet =0,\,\ldots\,,r$;

 \item[(3)]
  $( m_{(j),j} m_{(i),j} )\,( m_{(j),j} m_{(j),i})\,
   =\, m_{(i),i} m_{(j),j}\,$,
  $\,i,j = 0,\,\ldots\,,r$; \hspace{1ex}
  cf.\ Lemma 4.2.1;  

 \item[(4)]
  $m_{(j),j} m_{(i),\bullet}\,
   =\, m_{(j),\bullet} \cdot ( m_{(j),j} m_{(i),j} )\,$
  $\,i,j, \bullet =0,\,\ldots\,,r$; \hspace{1ex}
  cf.\ the gluing $U_{ij}\stackrel{\sim}{\leftarrow}U_{ji}$.
\end{itemize}
$\GL_n({\Bbb C})$ acts diagonally on
 $\prod_{i=0}^r\,
  \Mor^{\ringsetscriptsize}({\Bbb C}[\frac{y_0}{y_i}\,,\,
                     \cdots\,,\,\frac{y_r}{y_i}], M_n({\Bbb C}))$,
  via the post-composition with the adjoint $\GL_n({\Bbb C})$-action
  on $M_n({\Bbb C})$,  and
the above system of conditions describes a $\GL_n({\Bbb C})$-invariant
 closed subset therein.
The space of D0-branes on ${\Bbb P}^r$ is given by\\
 $\Map((\Space M_n({\Bbb C});{\Bbb C}^n),{\Bbb P}^r)
  =\Mor(\Space M_n({\Bbb C}),{\Bbb P}^r)/\!\sim\,$,
 described by the orbit-space of the $\GL_n({\Bbb C})$-action
 on the above subset in
 $\prod_{i=0}^r\,
  \Mor^{\ringsetscriptsize}({\Bbb C}[\frac{y_0}{y_i}\,,\,
                     \cdots\,,\,\frac{y_r}{y_i}], M_n({\Bbb C}))$.

The Chan-Paton modules of D0-branes on ${\Bbb P}^r$ and
 their Higgsing/un-Higgsing behavior follow the reasoning
 that combines the cases $Y={\Bbb P}^1$ and $Y={\Bbb A}^2$.
Together with
  the simultaneous triangularizability of any family of commuting matrices
   and
  the map that takes a tuple of triangularized matrices to the tuple of
   the respective diagonal,
 one has: (cf.\ Proposition 4.2.2)

\bigskip

\noindent
{\bf Proposition 4.4.1 [D0-branes on ${\Bbb P}^r$].} {\it
 There is an embedding
  $\Phi_{\scriptsizeHilb}:
   \Hilb^n_{{\scriptsizeBbb P}^r}=:({\Bbb P}^r)^{[n]}
    \rightarrow \Map((\Space M_n({\Bbb C});{\Bbb C}^n),{\Bbb P}^r)$.
 $\varphi_{\cal R}\in\Phi_{\scriptsizeHilb}(({\Bbb P}^r)^{[n]})$
  has the property that
   $\image\hat{\varphi}_{\cal R}$ is a subscheme of length $n$
   on ${\Bbb P}^r$.
 There is an embedding
  $\Phi_{\scriptsizeChow}:
   S^n({\Bbb P}^r)
    \rightarrow \Map((\Space M_n({\Bbb C});{\Bbb C}^n),{\Bbb P}^r)$,
  whose image is characterized by $\varphi_{\cal R}$ associated to
   a system of commuting diagonalizable matrices.
 (In particular, $\image\hat{\varphi}_{\cal R}$ is a reduced subscheme
   of length $\le n$ on ${\Bbb P}^r$.)
 There is a map
  $\Map((\Space M_n({\Bbb C});{\Bbb C}^n),{\Bbb P}^r)
   \rightarrow S^n({\Bbb P}^r)$
  that has $\Phi_{\scriptsizeChow}$ as a section.
 The pattern of open-string-induced Higgsing/un-Higgsing behavior of
  $n$ D0-branes on ${\Bbb P}^r$ can be reproduced
  in the current content via deformations of morphisms
  $[\varphi_{\cal R}]$
  in $\Phi_{\scriptsizeChow}(S^n({\Bbb P}^r))
      \subset \Map((\Space M_n({\Bbb C}); {\Bbb C}^n),{\Bbb P}^r)$.
} 

\bigskip

\begin{flushleft}
{\bf D0-branes on a quasi-projective variety.}
\end{flushleft}
Let $Y$ be a quasi-projective variety  and
suppose that $Y$ is embedded in ${\Bbb P}^r$ as $Y_1-Y_2$,
 where both $Y_1$ and $Y_2$ are closed subschemes of ${\Bbb P}^r$.
Let $I_1=\langle f_{11},\,\cdots\,,f_{1l_1}\rangle$
  (resp.\ $I_2=\langle f_{21},\,\cdots\,,f_{2l_2}\rangle$)
 be the homogeneous ideal in ${\Bbb C}[y_0,\,\cdots\,,y_r]$
 associated to $Y_1$ (resp.\ $Y_2$) in ${\Bbb P}^r$.
Recall the local affine charts $\cup_{i=0}^rU_i$ of ${\Bbb P}^r$.
Consider the (in general only quasi-affine) open cover
 $\cup_{i=0}^r((Y_1-Y_2)\cap U_i)$ of $Y$.
Then, the pair $(I_1, I_2)$ gives rise to a pair
 $$
  \left(\,
   I_{1,(i)}=( f_{11,(i)},\,\cdots\,,f_{1l_1,(i)} )\,,\,
   I_{2,(i)}=( f_{21,(i)},\,\cdots\,,f_{2l_2,(i)} )\,
  \right)
 $$
 of ideals in ${\Bbb C}[\frac{y_0}{y_i}\,,\,\cdots\,,\,\frac{y_r}{y_i}]$
 via the dehomogenization of $(I_1, I_2)$ on the affine chart
 $U_i$ of ${\Bbb P}^r$ for $i=0,\,\ldots\,,r$.
The space $\Mor(\Space M_n({\Bbb C}),Y)$ of morphisms from
 $\Space M_n({\Bbb C})$ to $Y$ is given by further restricting
 the locus $\Mor(\Space M_n({\Bbb C}),{\Bbb P}^r)$ in
 $\prod_{i=0}^r\,
  \Mor^{\ringsetscriptsize}({\Bbb C}[\frac{y_0}{y_i}\,,\,
                          \cdots\,,\,\frac{y_r}{y_i}], M_n({\Bbb C}))$,
 described by Conditions (1) - (4) in the previous theme,
 to the following system of incidence relation from $I_1$ and
 exclusion relations from $I_2$:
 \begin{itemize}
  \item[(5)]
   [({\it closed}) {\it incidence conditions from} $I_1$]$\,$:\\[.6ex]
   $\mbox{\hspace{1em}}
    f_{1\bullet, (i)}(m_{(i),0}\,,\,\cdots\,,\,m_{(i),r})\,=\,0\,
    \in M_n({\Bbb C})\,$, $\,\bullet=1,\,\ldots\,, l_1\,$,
   $\,i=0,\,\ldots,\,r\,$;

  \item[(6)]
   [({\it open}) {\it exclusion conditions from} $I_2$]$\,$:\\[.6ex]
   $\mbox{\hspace{1em}}
    m_{(i),i}\;\in\;
    \left\langle\,
     f_{2\bullet, (i)}(m_{(i),0}\,,\,\cdots\,,\,m_{(i),r})\,
    \right\rangle_{\bullet=1}^{l_2}\; \subset\; M_n({\Bbb C})\,$,
   $\,i=0,\,\ldots,\,r\,$.
 \end{itemize}
The diagonal $\GL_n({\Bbb C})$-action on
 $\prod_{i=0}^r\,
  \Mor^{\ringsetscriptsize}({\Bbb C}[\frac{y_0}{y_i}\,,\,
                     \cdots\,,\,\frac{y_r}{y_i}], M_n({\Bbb C}))$
 leaves the locally-closed subset that satisfies
 Conditions (1) - (6) invariant.
The space of D0-branes on $Y$ is given then by
 $\Map((\Space M_n({\Bbb C});{\Bbb C}^n),Y)
  =\Mor(\Space M_n({\Bbb C}),Y)/\!\sim\,$,
 described by the orbit-space of the $\GL_n({\Bbb C})$-action
 on the above locally-closed subset in
 $\Mor(\Space M_n({\Bbb C}),{\Bbb P}^r)$.

\bigskip

\noindent
{\it Remark 4.4.2 $[$Independence of embedding$]$.} {\rm
 The open cover $\cup_{i=0}^r((Y_1-Y_2)\cap U_i)$ of $Y$
  can be refined to an affine open cover of $Y$, which realizes
  $Y$ as a gluing system of rings.
 Different embeddings of $Y$ in projective spaces realizes $Y$
  as different gluing systems of rings that have a common refinement.
 It follows then from Sec.~1.2
  that $\Map((\Space M_n({\Bbb C});{\Bbb C}^n), Y)$ thus constructed
  is independent of the embedding of $Y$ in a projective space.
} 

\bigskip

Proposition 4.4.1  
 implies then:

\bigskip

\noindent
{\bf Theorem 4.4.3 [D0-branes on quasi-projective variety].} {\it
 Let $Y$ be a quasi-projective variety over ${\Bbb C}$.
 $\;(1)$
 There is an embedding
  $\Phi_{\scriptsizeHilb}:
   \Hilb^n_Y=:Y^{[n]}
    \rightarrow \Map((\Space M_n({\Bbb C});{\Bbb C}^n),Y)$.
 $\varphi_{\cal R}\in\Phi_{\scriptsizeHilb}(Y^{[n]})$
  has the property that
   $\image\hat{\varphi}_{\cal R}$ is a subscheme of length $n$ on $Y$.
 $\;(2)$
 There is an embedding
  $\Phi_{\scriptsizeChow}:
   S^nY\rightarrow \Map((\Space M_n({\Bbb C});{\Bbb C}^n),Y)$,
  whose image is characterized by $\varphi_{\cal R}$ associated to
   a system of commuting diagonalizable matrices.
 (In particular, $\image\hat{\varphi}_{\cal R}$ is a reduced subscheme
   of length $\le n$ on $Y$.)
 $\;(3)$
 There is a map
  $\Map((\Space M_n({\Bbb C});{\Bbb C}^n),Y)\rightarrow S^nY$
  that has $\Phi_{\scriptsizeChow}$ as a section.
 $\;(4)$
 The pattern of open-string-induced Higgsing/un-Higgsing behavior of
  $n$ D0-branes on $Y$ can be reproduced
  in the current content via deformations of morphisms
  $[\varphi_{\cal R}]$
  in $\Phi_{\scriptsizeChow}(S^nY)
           \subset \Map((\Space M_n({\Bbb C}); {\Bbb C}^n),Y)$.
} 

\bigskip

\noindent
{\it Remark 4.4.4 $[$toric variety$]$.} {\rm
 The discussions for D0-branes on ${\Bbb P}^r$
  (resp.\ a quasi-projective variety)
  generalize immediately to D0-branes
  on a toric variety (resp.\ a subscheme of a toric variety).
} 

\bigskip

\begin{flushleft}
{\bf D0-branes, gauged matrix models, and quantum moduli spaces.}
\end{flushleft}
When $Y$ is a closed subvariety of a toric variety/${\Bbb C}$,
the space $\Map((\Space M_n({\Bbb C}); {\Bbb C}^n),Y)$ are described
 by a system of noncommutative-polynomial-like algebraic equations
 that give only closed conditions.
In this case, $\Map((\Space M_n({\Bbb C}); {\Bbb C}^n),Y)$ is
 realizable as the classical moduli space of vacua
 (also known as vacuum manifold/variety) of a gauged matrix model.
The construction is similar to that of [Wi1] but adjusted to
 $d=0+1$ matrix models.
See also the discussions in [D-G-M], [Do-M], and [G-L-R]
 for related situations and [L-Y5] for further discussions.
The real issue, particularly from the mathematical/geometric aspect,
 is whether there is or needs to be also a good/mathematical notion
 of {\it quantum moduli space} in this case to incorporate more physics
 into the current mathematical setting.
In the next theme, we will see an example from string theory
 in which $\Map((\Space M_n({\Bbb C}); {\Bbb C}^n),Y)$
 already contains both a classical and a quantum moduli space
 of D0-branes on $Y$ in the sense of [Vafa1].

\bigskip

\begin{flushleft}
{\bf A comparison with the moduli problem of gas of D0-branes
     in {\rm [Vafa1]} of Vafa.}
\end{flushleft}
In [Vafa1], Vafa studied, among other things, the physics of
 finitely many D0-branes and D4-branes.
In particular, for a gas of $n$-many identical D0-branes
  on one D4-brane supported on a complex torus ${\Bbb T}^4$,
 except the additional $U(1)$-factor in the whole gauge group that
  comes from the simple D4-brane, the Higgsing/un-Higgsing behavior
  of such D0-D4 systems is the same as
  that for $n$-many D0-branes alone  and
 the classical moduli/configuration space of the $n$-many D0-branes
  on the ${\Bbb T}^4$ is given by $S^n({\Bbb T}^4)$,
  which is a singular complex space.
This moduli space is subject to a quantum correction to
 a quantum moduli space
 $\widetilde{S^n({\Bbb T}^4)}$, dictated by the requirement that
 the cohomology $H^{\ast}(\widetilde{S^n({\Bbb T}^4)},{\Bbb C})$
 should be the orbifold cohomology (e.g.\ [V-W1] and [V-W2])
 of $S^n({\Bbb T}^4)$ from string theory.
It is also anticipated that $\widetilde{S^n({\Bbb T}^4)}$ should be
 a hyperk\"{a}hler resolution of $S^n({\Bbb T}^4)$.
See also related discussions in [B-V-S1], [B-V-S2], and [Vafa2].

The related orbifold cohomology was later constructed mathematically
 by Chen and Ruan in [C-R1] and [C-R2].
In [Ru: Conjecture 6.3], Ruan conjectured in particular that,
 for $Y$ a smooth projective surface over ${\Bbb C}$
 such that $Y^{[n]}$ has a hyperk\"{a}hler structure,
 the orbifold cohomology ring $H^{\ast}_{\scriptsizeorb}(S^nY,{\Bbb C})$
 of $S^nY$ is isomorphic to the (ordinary) cohomology ring
 $H^{\ast}(X^{[n]};{\Bbb C})$ of $X^{[n]}$.
For the case $Y$ is a smooth projective surface$/{\Bbb C}$
  with trivial canonical line bundle,
 this was proved by Uribe [Ur: Theorem 3.2.3] together with
 previous result of Lehn and Sorger in [L-S].
Thus, for $Y$ a smooth projective Calabi-Yau surface,
 the $\widetilde{S^nY}$ anticipated in [Vafa1] is $Y^{[n]}$.

In our current setting,
a gas of $n$-many D0-branes on a D4-brane\footnote{A
                                               complete treatment of this
                                               involves an intrinsic
                                               mathematical
                                               construction/definition of
                                               a bound system of D-branes.
                                              Here, we only consider
                                               the pure D0-brane
                                               sector/factor
                                               in such a system.},
 supported on a smooth projective surface $Y$, is regarded
 as the image of a morphism from $(\Space M_n({\Bbb C});{\Bbb C}^n)$
 to $Y$.
The moduli space $\Map((\Space M_n({\Bbb C});{\Bbb C}^n),Y)$
 of such morphisms contains both
 $Y^{[n]}\simeq \image\Phi_{\scriptsizeHilb}$ and
 $S^nY\simeq \image\Phi_{\scriptsizeChow}$,
and the restriction of
 $\pi_{\scriptsizeHilb}:
   \Map((\Space M_n({\Bbb C});{\Bbb C}^n),Y)\rightarrow S^nY$
 to $\image\Phi_{\scriptsizeHilb}$ realizes the resolution
 $Y^{[n]}\rightarrow S^nY$.
In the special case that $Y$ is in addition Calabi-Yau, we see that
 $\Map((\Space M_n({\Bbb C});{\Bbb C}^n),Y)$ contains
 {\it both} the classical and the quantum moduli space of
 D0-brane configurations on $Y$ in [Vafa1].

\bigskip

\subsection{A remark on D-branes and universal moduli space.}

In the previous subsections, we see an interesting feature of
 the moduli space of D0-branes
 on a (commutative) quasi-projective variety:
 namely, it incorporates both the Hilbert scheme and the Chow variety.
We also see in the end of Sec.~4.4 that in a special occasion
 this is interpreted as containing both the classical and
 the quantum moduli space of D0-branes in physics.

While the encompassing of both the classical and the quantum moduli
 space of a D-brane system on a string target space in general
 is an issue that will be subject to how we formulate
 the intrinsic definition of D-brane bound system,
the unifying feature of the moduli space of D-branes on a target space
 (in the sense of Definition 2.2.3
  and its extension/generalization to systems that contains
  NS-branes as well) for different moduli spaces
  (e.g.\ Hilbert schemes and Chow varieties in the above example)
  in commutative geometry should be an anticipated feature
  when the mathematical definition/formulation of D-branes is
  ``correct".
Indeed, since 1995 new stringy dualities have made predictions
 that relate invariants of different mathematical origins, e.g.\
 from the stable maps, the stable/torsion sheaves, and subschemes
 respectively (when put in the setting of algebraic geometry).
These stringy dualities involve D-branes at work.
It is thus natural to anticipate that all these standard moduli spaces
 that appear in the mathematical definition of these invariants
 should live in different, possibly partially-overlapped
 regions/corners of the moduli space of D-branes (or in general
  D-branes coupled with NS-branes) on a target space.
This anticipation is particularly compelling from the viewpoint of
 Wilson's theory-space underlying these stringy dualities;
cf.\ [Liu2] and [L-Y1: appendix A.1].

\newpage

\end{document}